%% file: AAMR_Subspaces.tex
\documentclass[11pt,a4paper]{article}

\usepackage{amsmath}
\usepackage{amsfonts}
\usepackage{amssymb}
\usepackage{centernot}
\usepackage[margin=2.5cm]{geometry}
\usepackage[colorlinks=true,linkcolor=blue,citecolor=blue,urlcolor=blue]{hyperref}
\usepackage{enumitem}
\usepackage{graphicx,fancybox}
\usepackage{algorithm2e}
\usepackage{tikz}
\usepackage{subfigure}
\usepackage{todonotes}
\usepackage[capitalise]{cleveref}
\usetikzlibrary{calc}
\usepackage{hhline}
\usepackage{multirow}
\usepackage{pgf}
\usepackage{bigstrut}
\usepackage{rotating}
\usepackage{sectsty}
\usepackage[font={small}]{caption}

\subsectionfont{\fontsize{11}{15}\selectfont}
\DeclareMathOperator*{\argmin}{argmin}

\renewcommand{\arraystretch}{1.2}

\newcommand{\Fix}{\operatorname{Fix}}

\newcommand{\diag}{\operatorname{diag}}

\newcommand{\Hi}{\mathcal{H}}
\newcommand{\R}{\mathbb{R}}

\newcommand{\Cnn}{\mathbb{C}^{n\times n}}

\Crefname{fact}{Fact}{Facts}
\usepackage{amsthm}
\newtheorem{theorem}{Theorem}[section]
\newtheorem{definition}{Definition}[section]
\newtheorem{proposition}{Proposition}[section]
\newtheorem{corollary}{Corollary}[section]

\newtheorem{fact}{Fact}[section]
\newtheorem{remark}{Remark}[section]

\title{Optimal rates of linear convergence of the averaged alternating modified reflections method for two subspaces}

\author{Francisco J. Arag\'on Artacho\thanks{Department of Mathematics,
University of Alicante, \textsc{Spain}. e-mail:~\url{francisco.aragon@ua.es}}
        \and Rub\'en Campoy\thanks{Department of Mathematics,
University of Alicante, \textsc{Spain}. e-mail:~\url{ruben.campoy@ua.es}}
}

\begin{document}
\maketitle

\begin{abstract}
The averaged alternating modified reflections (AAMR) method is a projection algorithm for finding the closest point in the intersection of convex sets to any arbitrary point in a Hilbert space. This method can be seen as an adequate modification of the Douglas--Rachford method that yields a solution to the best approximation problem.
In this paper we consider the particular case of two subspaces in a Euclidean space. We obtain the rate of linear convergence of the AAMR method in terms of the Friedrichs angle between the subspaces and the parameters defining the scheme, by studying the linear convergence rates of the powers of matrices. We further optimize the value of these parameters in order to get the minimal convergence rate, which turns out to be better than the one of other projection methods. Finally, we provide some numerical experiments that demonstrate the theoretical results.
\end{abstract}



\section{Introduction}\label{sec:intro}
The \emph{averaged alternating modified reflections (AAMR) algorithm}, introduced in~\cite{AAMR}, is a projection method for solving best approximation problems in the convex setting. A best approximation problem consists in finding the closest point to any given point in the intersection of a collection of sets. In this work we study problems involving two subspaces $U$ and $V$ in $\mathbb{R}^n$.

Given a point $z\in\R^n$, the corresponding best approximation problem is defined as
\begin{equation}\label{eq:bestaproxprob}
\text{Find } w\in U\cap V\text{ such that } \|w-z\|=\inf_{x\in U\cap V} \|x-z\|.
\end{equation}
For any initial point $x_0\in\R^n$, the AAMR algorithm is iteratively defined by
\begin{equation}\label{eq:aamr_intro}
x_{k+1}:=(1-\alpha)x_k+\alpha(2\beta P_{V-z}-I)(2\beta P_{U-z} - I)(x_k), \quad k=0,1,2\ldots.
\end{equation}
When $\alpha,\beta\in\,]0,1[$, the generated sequence $\{x_k\}_{k=0}^\infty$ converges to a point $x^\star$ such that
\begin{equation*}
P_U(x^\star+z)=P_{U\cap V}(z),
\end{equation*}
which solves problem~\eqref{eq:bestaproxprob}. Furthermore, the shadow sequence $\left\{P_U(x_k+z)\right\}_{k=0}^\infty$ is convergent to the solution $P_{U\cap V}(z)$ even if $\alpha=1$, see~\cite[Theorem~4.1]{AAMR}. In fact, when the sets involved are subspaces, we prove that the sequence $\{x_k\}_{k=0}^\infty$ is also convergent for $\alpha=1$, see~\Cref{cor:convergent}.

Several projection methods have been developed for solving convex feasibility problems in Hilbert spaces, see e.g.~\cite{BC11,C12,CC15,D01,ER11}. In the case where the sets are subspaces, some of these methods converge to the closest point in the the intersection to the starting point, providing thus a solution of the best approximation problem~\eqref{eq:bestaproxprob}. Among these schemes, probably the two most well-known are the method of alternating projections (AP), which was originally introduced by John von Neumann~\cite{VN50}, and the Douglas--Rachford method (DR)~\cite{DR56,LM79}, which is also referred as averaged alternating reflections. The rate of linear convergence of these methods is known to be the cosine of the Friedrichs angle between the subspaces for DR~\cite{BCNPW14}, and the squared cosine of this angle for AP~\cite{D95}. Several relaxations and generalizations of these methods have been proposed, such as the relaxed and the partial relaxed alternating projections (RAP, PRAP) \cite{Ag54,C08,MS54}, the generalized alternating projections (GAP)~\cite{FG16,GPR67}, the relaxed averaged alternating reflections (RAAR) \cite{L08}, and the generalized Douglas--Rachford (GDR) \cite{EB92}, among others. We note that AAMR can also be seen as a modified version of DR, since both methods coincide when $\alpha=\frac{1}{2}$ and $\beta=1$ in \eqref{eq:aamr_intro}.

Thanks to the linearity of the projector operator onto subspaces, projection methods reduce to matrix iterations. Taking advantage of this fact, optimal convergence rates have been obtained in~\cite{BCNPW15} for RAP, PRAP and GDR. By following an analogous matrix analysis, the rate of convergence with optimal parameters for GAP has been recently given in~\cite{FG17}.

In the current setting, the rate of convergence of the AAMR algorithm  was numerically analyzed in various computational experiments in~\cite[Section~7]{AAMR}. The goal of this work is to provide the theoretical results that substantiate the behavior of the algorithm that was numerically observed. By following the same approach as in~\cite{BCNPW15}, we analyze the linear rate of convergence of the AAMR method by studying the convergence rates of powers of matrices. The rate obtained depends on both the Friedrichs angle and the parameters defining the algorithm. In addition, we also obtain the optimal selection of the parameters according to the Friedrichs angle, so that the rate of convergence is minimized. This rate coincides with the one for GAP, which is the best among the rates of all the projection methods mentioned above. This is not just by chance: the shadow sequences of GAP and AAMR coincide for linear subspaces under some conditions (see \Cref{th:GAP_AAMR} and \Cref{fig:th_GAP_AAMR}).

The remaining of the paper is structured as follows. We present some definitions and preliminary results in~\Cref{sec:preliminaries}. In~\Cref{sec:rate} we collect our main results regarding the rate of convergence of the AAMR method. We compare the rate with optimal parameters of AAMR with the rate of various projection methods in~\Cref{sec:comparison}. In~\Cref{sec:experiments} we perform two computational experiments that validate the theoretical results obtained. We finish with some conclusions and future work in~\Cref{sec:conclusions}.

\section{Preliminaries}\label{sec:preliminaries}
In this work, our setting is the Euclidean space $\R^n$ with inner product $\langle\cdot ,\cdot\rangle$ and induced norm~$\|\cdot\|$. For a given set $C\subset\R^n$ we denote by $C^\bot=\{x\in\Hi : \langle c,x\rangle =0, \forall c\in C\}$ the \emph{orthogonal complement} of $C$. Given $x\in \R^n$, a point $p\in C$ is said to be a \emph{best approximation} to $x$ from $C$~if
\begin{equation*}
\|p-x\|=d(x,C):=\inf_{c\in C}\|c-x\|.
\end{equation*}
The operator $P_C(x):=\arg\!\min\left\{\|x-c\|, c\in C\right\}$ is called the \emph{projector} onto $C$. When $C$ is closed and convex, $P_C$ is single-valued. In the case when $C$ is a subspace, $P_C(x)$ is sometimes called the \emph{orthogonal projection} of $x$ to $C$, due to the fact that $x-P_C(x)\in C^\bot$.

Throughout this paper, we assume without lost of generality that $U$ and $V$ are two subspaces of $\R^n$ such that $1\leq p:=\dim U\leq\dim V=:q\leq n-1$, with $U\neq U\cap V$ and $U\cap V\neq\{0\}$  (otherwise, problem~\eqref{eq:bestaproxprob} would be trivial). We now recall the concept of principal angles and Friedrichs angle between a pair of subspaces, and a result relating both concepts.

\begin{definition}
The \emph{principal angles} between $U$ and $V$ are the angles $0\leq\theta_1\leq\theta_2\leq\cdots\leq~\theta_p\leq~\frac{\pi}{2}$ whose cosines are recursively defined by
\begin{align*}
\cos\theta_k:=&\langle u_k,v_k\rangle\\
=&\max\left\{\langle u,v\rangle : u\in U, v\in V, \|u\|=\|v\|=1, \langle u,u_j\rangle=\langle v,v_j\rangle=0 \text{ for } j=1,\ldots,k-1   \right\},
\end{align*}
with $u_0=v_0:=0$.
\end{definition}

\begin{definition}\label{def:F_angle}
	The \emph{Friedrichs angle} between $U$ and $V$ is the angle in $\theta_F\in\,]0,\frac{\pi}{2}]$ whose cosine~is
	$$
	c_F(U,V):=\sup\left\{\langle u, v \rangle \, : \, u\in U\cap(U\cap V)^\bot, v\in V\cap(U\cap V)^\bot, \|u\|\leq 1,\|v\|\leq 1 \right\}.
	$$
\end{definition}

\begin{fact}\label{fact:principal_angles}
Let $\theta_1,\theta_2,\ldots,\theta_p$ be the principal angles between $U$ and $V$, and let $s:=\dim (U\cap V)$. Then we have $\theta_k=0$ for $k=1,\ldots,s$ and $\theta_{s+1}=\theta_F>0$.
\end{fact}
\begin{proof}
See~\cite[Proposition~3.3]{BCNPW15}.
\end{proof}
\begin{remark}
By our standing assumption that $U\neq U\cap V$, we have $s=\dim (U\cap V)<p$.
\end{remark}

The projector operator onto subspaces is known to be a linear mapping. The following result provides a matrix representation of the projectors onto $U$ and $V$, according to their principal angles. We denote by $I_n$, $0_n$ and $0_{m\times n}$, the $n\times n$ identity matrix, the $n\times n$ zero matrix, and the $m\times n$ zero matrix, respectively. For simplicity, we shall omit the subindices when the size can be deduced.

\begin{fact}\label{fact:projections_principal_angles}
If $p+q<n$, we may find an orthogonal matrix $D\in\R^{n\times n}$ such that
\begin{equation}\label{eq:projections_principal_angles}
P_U=D\begin{pmatrix}
I_p & 0 & 0 & 0\\
0 & 0_p & 0 & 0\\
0 & 0 & 0_{q-p} & 0\\
0 & 0 & 0 & 0_{n-p-q}
\end{pmatrix}D^* \quad\text{and}\quad
P_V=D\begin{pmatrix}
C^2 & CS & 0 & 0\\
CS & S^2 & 0 & 0\\
0 & 0 & I_{q-p} & 0\\
0 & 0 & 0 & 0_{n-p-q}
\end{pmatrix}D^*,
\end{equation}
where $C$ and $S$ are two $p\times p$ diagonal matrices defined by
$$C:=\diag(\cos\theta_1,\ldots,\cos\theta_p)\quad\text{and}\quad S:=\diag(\sin\theta_1,\ldots,\sin\theta_p),$$
with $\theta_1,\ldots,\theta_p$ being the principal angles between $U$ and $V$.
\end{fact}
\begin{proof}
See~{\cite[Proposition~3.4]{BCNPW15}}.
\end{proof}

\subsection{The averaged alternating modified reflections operator for two subspaces}\label{subsec:aamr}

The AAMR operator was originally introduced for two arbitrary closed and convex sets \cite[Definition~3.2]{AAMR}. In this section, we present the scheme in the case of two subspaces, as well as some properties of the operator and its set of fixed points within this context.

\begin{definition}
	Given $\alpha\in\,]0,1]$ and $\beta\in\,]0,1[$,  the \emph{averaged alternating modified reflections (AAMR) operator} is the mapping $T_{U,V,\alpha,\beta}:\R^n \mapsto \R^n$ given by
	\begin{equation*}
	T_{U,V,\alpha,\beta}:=(1-\alpha)I+\alpha(2\beta P_V-I)(2\beta P_U - I).
	\end{equation*}
	Where there is no ambiguity, we shall abbreviate the notation of the operator $T_{U,V,\alpha,\beta}$ by~$T_{\alpha,\beta}$.	
\end{definition}

\begin{fact}\label{fact:nonexpansive}
Let $\alpha\in\,]0,1]$ and $\beta\in\,]0,1[$. Then, the AAMR operator
$T_{\alpha,\beta}$ is nonexpansive.
\end{fact}
\begin{proof}
See~\cite[Proposition~3.3]{AAMR}.
\end{proof}

\begin{proposition}\label{prop:FixT}
Let $\alpha\in\,]0,1]$ and $\beta\in\,]0,1[$. Then
\begin{equation*}
\Fix T_{U,V,\alpha,\beta}=U^\perp \cap V^\perp.
\end{equation*}
\end{proposition}
\begin{proof}
Observe that
\begin{equation}\label{eq:caract_fixedpoints}
x\in\Fix T_{U,V,\alpha,\beta} \Leftrightarrow P_V(2\beta P_U(x)-x)=P_U(x).
\end{equation}
Moreover, by~\cite[Proposition~3.4]{AAMR},
$$P_U(x)=P_{U\cap V}(0), \quad \text{ for all } x\in\Fix T_{U,V,\alpha,\beta}.$$
Therefore, if $x\in \Fix T_{U,V,\alpha,\beta}$, then $0=P_{U\cap V}(0) =P_U(x)$. Using this equality, together with \eqref{eq:caract_fixedpoints}, we deduce that $P_U(x)=P_V(x)=0$, which implies $x\in U^\perp \cap V^\perp$.

To prove the converse implication, pick any $x\in U^\perp \cap V^\perp$. Then we trivially have that $P_V(2\beta P_U(x)-x)=P_U(x)$, and thus $x\in\Fix T_{U,V,\alpha,\beta}$.
\end{proof}

The AAMR scheme is iteratively defined by~\eqref{eq:aamr_intro}. Using the linearity of the projector operator onto subspaces, we deduce that the iteration takes the form
\begin{align}\label{eq:AAMR_subspaces}
x_{k+1}&=(1-\alpha)x_k+\alpha(2\beta P_{V-z}-I)(2\beta P_{U-z} - I)(x_k)\nonumber\\
&=(1-\alpha)x_k+\alpha(2\beta P_{V-z}-I)\left(2\beta \left(P_{U}(x_k+z) -z\right) -x_k\right)\nonumber\\
&=(1-\alpha)x_k+\alpha\left(2\beta P_{V-z}\left(2\beta \left(P_{U}(x_k+z) -z\right) -x_k\right)-2\beta \left(P_{U}(x_k+z) -z\right) +x_k\right)\nonumber\\
&=x_k+2\alpha\beta\left( P_{V}\left(2\beta \left(P_{U}(x_k+z) -z\right) -x_k+z\right)-P_{U}(x_k+z)\right)\nonumber\\
&=x_k+2\alpha\beta\left( 2\beta P_{V}P_{U}(x_k+z) +(1-2\beta)P_V(z) -P_V(x_k)-P_{U}(x_k+z)\right).
\end{align}

\begin{fact}\label{fact:translation_aamr}
	Let $\alpha\in\,]0,1]$, $\beta\in\,]0,1[$ and  $z\in\R^n$. Then, one has $\Fix T_{U-z,V-z,\alpha,\beta}\neq\emptyset$ and
	\begin{equation*}\label{eq:fact_affine_traslation_FixT2}
	\Fix T_{U-z,V-z,\alpha,\beta}=x^*+ U^\perp \cap V^\perp, \quad \forall x^*\in\Fix T_{U-z,V-z,\alpha,\beta}.
	\end{equation*}
	Furthermore, for any $x\in\R^n$,
	\begin{equation*}\label{eq:fact_affine_traslation_FixT}
	T_{U-z,V-z,\alpha,\beta}(x)=T_{U,V,\alpha,\beta}(x-x^*)+x^*, \quad \forall x^*\in\Fix T_{U-z,V-z,\alpha,\beta}.
	\end{equation*}
\end{fact}
\begin{proof}
Since $U$ and $V$ are subspaces in a finite dimensional space, by~\cite[Fact~2.11 and Corollary~3.1]{AAMR} we get that $\Fix T_{U-z,V-z,\alpha,\beta}\neq \emptyset$. The remaining assertions are obtained by applying \Cref{prop:FixT} and \cite[Proposition~3.6]{AAMR} to $U-q$,$V-q$ and $-q\in(U-q)\cap(V-q)$, noting that \cite[Proposition~3.6]{AAMR} also holds for $\alpha=1$.
\end{proof}

\subsection{Optimal convergence rate of powers of matrices}\label{subsec:matrices}

We denote by $\Cnn$ ($\mathbb{R}^{n\times n}$), the space of $n\times n$ complex (real) matrices, equipped with the induced  matrix norm~$\|A\|:=\max\left\{\|Ax\| : x\in\mathbb{C}^{n}, \|x\|\leq 1\right\}$. The \emph{kernel} of a matrix $A\in\mathbb{C}^{n\times n}$ is denoted by $\ker A:=\left\{x\in\Cnn : Ax=0\right\}$ and the set of \emph{fixed points} of $A$ is denoted by $\Fix A:=\ker (A-I)$. We say $A$ is nonexpansive if $\|Ax-Ay\|\leq\|x-y\|$ for all $x,y\in\Cnn$.

\begin{definition}
A matrix $A\in\Cnn$ is said to be convergent to $A^\infty\in\Cnn$ if and only if
\begin{equation*}
\lim_{k\mapsto\infty}\|A^k-A^\infty\|=0.
\end{equation*}
We say $A$ is \emph{linearly} convergent to $A^\infty$ with rate $\mu\in[0,1[$ if there exist a positive integer $k_0$ and some $M>0$ such that
\begin{equation*}
\|A^k-A^\infty\|\leq M\mu^k, \quad \text{for all } k\geq k_0.
\end{equation*}
In this case, $\mu$ is called a \emph{linear convergence rate} of $A$. When the infimum of all the convergence rates is also a convergence rate, we say this minimum is the \emph{optimal linear convergence rate}.
\end{definition}

For any matrix $A\in\Cnn$, we denote by $\sigma(A)$ the \emph{spectrum} of $A$ (the set of all eigenvalues). An eigenvalue $\lambda\in\sigma(A)$ is said to be \emph{semisimple} if its algebraic multiplicity coincides with its geometric multiplicity (cf.~\cite[p. 510]{Meyer}), or, equivalently, if $\ker(A-\lambda I)=\ker \left((A-\lambda I)^2\right)$ (see~\cite[Fact~2.3]{BCNPW15}). The \emph{spectral radius} of $A$ is defined by
$$\rho(A) := \max\{|\lambda|: \lambda\in\sigma(A)\},$$
and the second-largest modulus of the eigenvalues of $A$ after $1$ is denoted by
$$
\gamma(A):=\max\left\{|\lambda| : \lambda\in\{0\}\cup\sigma(A)\setminus\{1\}\right\}.
$$
An eigenvalue $\lambda\in\sigma(A)$ with $|\lambda |=\gamma(A)$ is called a \emph{subdominant} eigenvalue.

\begin{fact}\label{fact:convergence}
Let $A\in\Cnn$. Then $A$ is convergent if and only if one of the following holds:
\begin{itemize}[nosep]
\item[(i)]$\rho(A) < 1$;
\item[(ii)]$\rho(A)=1$ and $\lambda=1$ is semisimple and is the only eigenvalue on the unit circle.
\end{itemize}
When this happens, $A$ is linearly convergent with any rate $\mu\in\,]\gamma(A),1[\,$, and $\gamma(A)$ is the optimal linear convergence rate of $A$ if and only if all the subdominant eigenvalues are semisimple.
Furthermore, if $A$ convergent and nonexpansive, then $\lim_{k\to\infty}A^k=P_{\Fix A}$.
\end{fact}
\begin{proof}
See~\cite[pp.~617--618,~630]{Meyer} and~\cite[Theorem~2.12, Theorem~2.15 and~Corollary~2.7(ii)]{BCNPW15}.
\end{proof}

\section{Convergence rate analysis}\label{sec:rate}

We begin this section with the following theorem that establishes the rate of convergence of the AAMR algorithm in terms of $\alpha$, $\beta$ and the Friedrichs angle between the subspaces. We denote the positive part of $x\in\R$ by $x^+:=\max\{0,x\}$.

\begin{theorem}\label{th:rate}
Let $\alpha\in\,]0,1]$ and $\beta\in\,]0,1[\,$. Then, the AAMR operator $$T_{\alpha,\beta}:=(1-\alpha)I_n+\alpha(2\beta P_V-I_n)(2\beta P_U-I_n)$$ is linearly convergent to $P_{U^\perp \cap V^\perp}$ with any rate $\mu\in\,]\gamma(T_{\alpha,\beta}),1[\,$, where
\begin{equation}\label{eq:rate}
\gamma(T_{\alpha,\beta})=\left\{\begin{array}{ll}
1-4\alpha\beta(1-\beta),&\text{ if }0\leq c_F< c(\alpha,\beta);\\
\sqrt{4(1-\alpha)\alpha\beta^2 c_F^2+(1-2\alpha\beta)^2},&\text{ if } c(\alpha,\beta)\leq c_F <\widehat{c}_{\beta};\\
1+2\alpha\beta\left(\beta c_F^2-1+c_F\sqrt{\beta^2 c_F^2-2\beta+1}\right),&\text{ if } \widehat{c}_{\beta}\leq c_F <1;
\end{array}\right.
\end{equation}
with $c_F:=\cos\theta_F$  and $\theta_F$ being the Friedrichs angle between $U$ and $V$,
\begin{equation}\label{eq:cs}
\widehat{c}_\beta :=\frac{\sqrt{(2\beta-1)^+}}{\beta}\quad\text{and}\quad c(\alpha,\beta):=\left\{\begin{array}{cl}
\sqrt{\frac{\left(\left(1-4\alpha\beta(1-\beta)\right)^2-(1-2\alpha\beta)^2\right)^+}{4(1-\alpha)\alpha\beta^2}},&\text{ if } \alpha<1;\\
0,&\text{ if }\alpha=1.\end{array}\right.
\end{equation}
Furthermore, $\gamma(T_{\alpha,\beta})$ is the optimal linear convergence rate if and only if $\beta\neq\frac{1}{1+\sin\theta_F}$ or $\theta_F=\frac{\pi}{2}$.
\end{theorem}
\begin{proof} To prove the result, we consider two main cases.\\
Case 1: $p+q<n$. By~\Cref{fact:projections_principal_angles}, we can find an orthogonal matrix $D\in\R^{n\times n}$ such that~\eqref{eq:projections_principal_angles} holds. After some calculations, we obtain
\begin{align}\label{eq:T_matrix}
T_{\alpha,\beta}&=D\left(\begin{array}{ccc}
M_{\alpha,\beta} & 0 & 0\\
0 &(1-2\alpha\beta)I_{q-p} &0\\
0 & 0 &I_{n-p-q}
\end{array}\right)D^*,
\end{align}
where $$M_{\alpha,\beta}:=\left( \begin{array}{cc}
2\alpha\beta(2\beta-1)C^2+(1-2\alpha\beta)I_p & -2\alpha\beta CS \\
2\alpha\beta(2\beta-1)CS& 2\alpha\beta C^2+(1-2\alpha\beta)I_p
\end{array}\right).$$
Let $s:=\dim(U\cap V)$ and let $1=c_1=\cdots=c_s>c_{s+1}=c_F\geq c_{s+2}\geq\cdots \geq c_p\geq 0$ be the cosine of the principal angles $0=\theta_1=\cdots=\theta_s<\theta_{s+1}=\theta_F\leq \theta_{s+2}\leq\cdots \leq \theta_p\leq \frac{\pi}{2}$ between $U$ and $V$ (see~\Cref{fact:principal_angles}). By the block determinant formula (see, e.g.,~\cite[(0.8.5.13)]{HJ13}), we deduce after some algebraic manipulation that the spectrum of $T_{\alpha,\beta}$ is given by
$$ \sigma(T_{\alpha,\beta})=\left\{\begin{array}{ll}
\displaystyle\bigcup_{k=1}^p\left\{1+2\alpha\beta\left(\beta c_k^2-1\pm c_k\sqrt{\beta^2c_k^2-2\beta+1}\right)\right\}\cup\{1\},& \text{if } q=p;\\
\displaystyle\bigcup_{k=1}^p\left\{1+2\alpha\beta\left(\beta c_k^2-1\pm c_k\sqrt{\beta^2c_k^2-2\beta+1}\right)\right\}\cup\{1\}
\cup\{1-2\alpha\beta\},& \text{if } q>p.
\end{array}\right.$$
Then $\lambda_{k,r}:=1+2\alpha\beta\left(\beta c_k^2-1+(-1)^rc_k\sqrt{\beta^2c_k^2-2\beta+1}\right)$ are eigenvalues of $T_{\alpha,\beta}$, with $r=1,2$ and $k=1,\ldots,p$. Observe that $\lambda_{k,r}\in\mathbb{R}$ if $c_k\geq \beta^{-1}\sqrt{(2\beta-1)^+}=:\widehat{c}_\beta$, while $\lambda_k^r\in\mathbb{C}$ otherwise. To study the modulus of the eigenvalues $\lambda_{k,r}$, consider the function $f_{\alpha,\beta,r}:[0,1]\to\mathbb{R}$ given by
$$f_{\alpha,\beta,r}(c):=\left\{\begin{array}{ll}
\left(1+2\alpha\beta(\beta c^2-1)\right)^2-4\alpha^2\beta^2 c^2\left(\beta^2 c^2-2\beta+1\right),&\text{ if } c<\widehat{c}_\beta;\\
\left(1+2\alpha\beta\left(\beta c^2-1+(-1)^rc\sqrt{\beta^2 c^2-2\beta+1}\right)\right)^2,&\text{ if } c\geq\widehat{c}_\beta.
\end{array}\right.$$
Hence, one has $|\lambda_{k,r}|^2=f_{\alpha,\beta,r}(c_k)$.

Let us analyze some properties of the function $f_{\alpha,\beta,r}$. When $\widehat{c}_\beta>0$, observe that $f_{\alpha,\beta,r}$ is continuous at $\widehat{c}_\beta$, since
$$
\lim_{c\to\widehat{c}_\beta^{\hspace{1pt}-}}f_{\alpha,\beta,r}(c)=\lim_{c\to\widehat{c}_\beta^{\hspace{1pt}+}}f_{\alpha,\beta,r}(c)=\left(1-2\alpha(1-\beta)\right)^2.
$$
Define the auxiliary function $g_{\beta,r}(c):=\beta c^2-1+(-1)^rc\sqrt{\beta^2 c^2-2\beta+1}$ for $c\geq\widehat{c}_\beta$. Then,
$$f_{\alpha,\beta,r}(c)=\left\{\begin{array}{ll}
4(1-\alpha)\alpha\beta^2 c^2+(1-2\alpha\beta)^2,&\text{ if } c<\widehat{c}_\beta;\\
\left(1+2\alpha\beta g_{\beta,r}(c)\right)^2,&\text{ if } c\geq\widehat{c}_\beta.\\
\end{array}\right.$$
The derivative of $f_{\alpha,\beta,r}$ is given for $c\neq\widehat{c}_\beta$ by
$$f'_{\alpha,\beta,r}(c)=\left\{\begin{array}{ll}
8(1-\alpha)\alpha\beta^2 c,&\text{ if } c<\widehat{c}_\beta;\\
4\alpha\beta\left(1+2\alpha\beta g_{\beta,r}(c)\right)(-1)^{r} \frac{\left(\sqrt{\beta^2 c^2-2\beta+1}+(-1)^r\beta c\right)^2}{\sqrt{\beta^2 c^2-2\beta+1}} ,&\text{ if } c>\widehat{c}_\beta.
\end{array}\right.$$
Further, we claim that $1+2\alpha\beta g_{\beta,2}(c)>1+2\alpha\beta g_{\beta,1}(c)\geq 0$ for all $c>\widehat{c}_\beta$. Indeed, since
$$(2\beta^2c^2-2\beta+1)^2=4\beta^2c^2(\beta^2 c^2-2\beta+1) + (2\beta-1)^2\geq (2\beta c)^2(\beta^2 c^2-2\beta+1),$$
we deduce, after taking square roots and reordering, that
$$-1\leq 2\beta\left(\beta c^2-1-c\sqrt{\beta^2 c^2-2\beta+1}\right)=2\beta g_{\beta,1}(c)< 2\beta g_{\beta,2}(c),$$
from where the assertion easily follows.

All the above properties of the function $f_{\alpha,\beta,r}$ can be summarized as follows:
\begin{itemize}
\item For all $0\leq c< d\leq\widehat{c_\beta}$,
$$F_{\alpha,\beta}^0:=(1-2\alpha\beta)^2\leq f_{\alpha,\beta,r}(c)\leq f_{\alpha,\beta,r}(d)\leq \left(1-2\alpha(1-\beta)\right)^2.$$
\item For all $\widehat{c}_\beta\leq c< d\leq 1$,
\begin{align*}
F_{\alpha,\beta}^1&:=\left(1-4\alpha\beta(1-\beta)\right)^2 \leq f_{\alpha,\beta,1}(d) \leq f_{\alpha,\beta,1}(c)\\
 &\leq \left(1-2\alpha(1-\beta)\right)^2 \leq f_{\alpha,\beta,2}(c)<f_{\alpha,\beta,2}(d)\leq 1.
\end{align*}
\end{itemize}

In view of \Cref{fact:convergence}, we have to show that the eigenvalue $\lambda=1$ is semisimple and the only eigenvalue in the unit circle. According to the monotonicity properties of $f_{\alpha,\beta,r}$, we have that $|\lambda_{k,r}|\leq 1$ for all $k=1,\ldots,p$ and $r=1,2$. Further,
$$
|\lambda_{k,r}|=1 \Leftrightarrow k\in\{1,2,\ldots,s\} \text{ and } r=2,
$$
in which case $\lambda_{k,2}=1$. Thus, we have shown that $\rho\left(T_{\alpha,\beta}\right)=1$ and $\lambda=1$ is the only eigenvalue in the unit circle.

Let us see now that $\lambda=1$ is semisimple. First observe that, for any $\lambda\in\mathbb{C}$, given the block diagonal structure of $T_{\alpha,\beta}$,  one has
$$
\ker\left(T_{\alpha,\beta}-\lambda I\right)=\ker\left(\left(T_{\alpha,\beta}-\lambda I\right)^2\right) \iff \ker\left(M_{\alpha,\beta}-\lambda I\right)=\ker\left(\left(M_{\alpha,\beta}-\lambda I\right)^2\right).
$$
Then, we can compute
\begin{equation*}
M_{\alpha,\beta}-I=2\alpha\beta\left(\begin{array}{cc}
(2\beta-1)C^2-I_p & -CS\\
(2\beta-1)CS& -S^2
\end{array}\right).
\end{equation*}
Observe that the matrices $C$ and $S$ can be decomposed as
\begin{equation}\label{decompCS}
C=\left( \begin{array}{cc} I_s & 0\\ 0 & \widetilde{C}\end{array}\right)\quad\text{and} \quad S=\left( \begin{array}{cc} 0_s & 0 \\ 0 & \widetilde{S}\end{array}\right),
\end{equation}
where both $\widetilde{C}$ and $\widetilde{S}$ are diagonal matrices and $\widetilde{S}$ has strictly positive entries. Hence,
\begin{equation*}
M_{\alpha,\beta}-I=2\alpha\beta\left( \begin {array}{cccc}
-2(1-\beta)I_s & 0  & 0 & 0\\
0 & (2\beta-1)\widetilde{C}^2-I_{p-s}  & 0 & -\widetilde{C}\widetilde{S}\\
0 & 0  & 0 & 0\\
0 & (2\beta-1)\widetilde{C}\widetilde{S}  & 0 & -\widetilde{S}^2\\
\end {array} \right),
\end{equation*}
and one has that
$\ker\left(M_{\alpha,\beta}-I\right)=\ker\left(\left(M_{\alpha,\beta}-I\right)^2\right)$ if and only if $\ker\left(M_0\right)=\ker\left(M_0^2\right)$, where $$M_0:=\left(\begin{array}{cc}
(2\beta-1)\widetilde{C}^2-I_{p-s} & -\widetilde{C}\widetilde{S}\\
(2\beta-1)\widetilde{C}\widetilde{S} & -\widetilde{S}^2
\end{array}\right).$$
Since $\det\left(M_0\right)=\det\left(\widetilde{S}^2\right)\neq 0$ (again, by the block determinant formula), we conclude that $\lambda=1$ is a semisimple eigenvalue. Then, since $T_{\alpha,\beta}$ is nonexpansive by \Cref{fact:nonexpansive}, we have by \Cref{fact:convergence} that $T_{\alpha,\beta}$ is linearly convergent to $P_{\Fix T_{\alpha,\beta}}$ with any rate $\mu\in\,]\gamma(T_{\alpha,\beta}),1[$, and $\Fix T_{\alpha,\beta}=U^\perp \cap V^\perp$ by \Cref{prop:FixT}.

Furthermore, we can also deduce from the monotonicity properties of $f_{\alpha,\beta,r}$ that the subdominant eigenvalues of $T_{\alpha,\beta}$ are determined by
\begin{align*}
\gamma(T_{\alpha,\beta})&=\max\left\{|\lambda_{s+1,2}|,|\lambda_{1,1}|\right\}\\
&=\max\left\{\left|1+2\alpha\beta\left(\beta c_F^2-1+c_F\sqrt{\beta^2c_F^2-2\beta+1}\right)\right|,1-4\alpha\beta(1-\beta)\right\}.
\end{align*}
To prove~\eqref{eq:rate}, let us compute the value of $\gamma(T_{\alpha,\beta})$. If $c_F> \widehat{c}_{\beta}$, then $|\lambda_{1,1}|<|\lambda_{s+1,2}|$. Otherwise,
$$
|\lambda_{s+1,2}|\leq |\lambda_{1,1}| \Leftrightarrow f_{\alpha,\beta,2}(c_F)\leq f_{\alpha,\beta,1}(1) \Leftrightarrow 4(1-\alpha)\alpha\beta^2c_F^2\leq F_{\alpha,\beta}^1-F_{\alpha,\beta}^0.
$$
Consequently, if we define
$$c(\alpha,\beta):=\left\{\begin{array}{cl}
\sqrt{\frac{F_{\alpha,\beta}^1-F_{\alpha,\beta}^0}{4(1-\alpha)\alpha\beta^2}},&\text{ if } F_{\alpha,\beta}^1> F_{\alpha,\beta}^0;\\
0,&\text{ otherwise;}\end{array}\right.$$
which is equivalent to the expression in~\eqref{eq:cs}, we obtain~\eqref{eq:rate}. Three possible scenarios for $f_{\alpha,\beta,r}$ and the constants $c(\alpha,\beta)$ and $\widehat{c}_{\beta}$ depending on the values of $\alpha$ and $\beta$ are shown in \Cref{fig:f_alpha_beta_r}.
\begin{figure}[ht!]
	\centering
	 \subfigure[$0=c(\alpha,\beta)=\widehat{c}_{\beta}<1$]{\scalebox{.4}{\input{Fig1a.pgf}}}\\ 	 \subfigure[$0=c(\alpha,\beta)<\widehat{c}_{\beta}<1$]{\scalebox{.4}{\input{Fig1b.pgf}}} 	 \subfigure[$0<c(\alpha,\beta)<\widehat{c}_{\beta}<1$]{\scalebox{.4}{\input{Fig1c.pgf}}} 	 
\caption{The three possible scenarios for the function $f_{\alpha,\beta,r}(c)$} \label{fig:f_alpha_beta_r}
\end{figure}
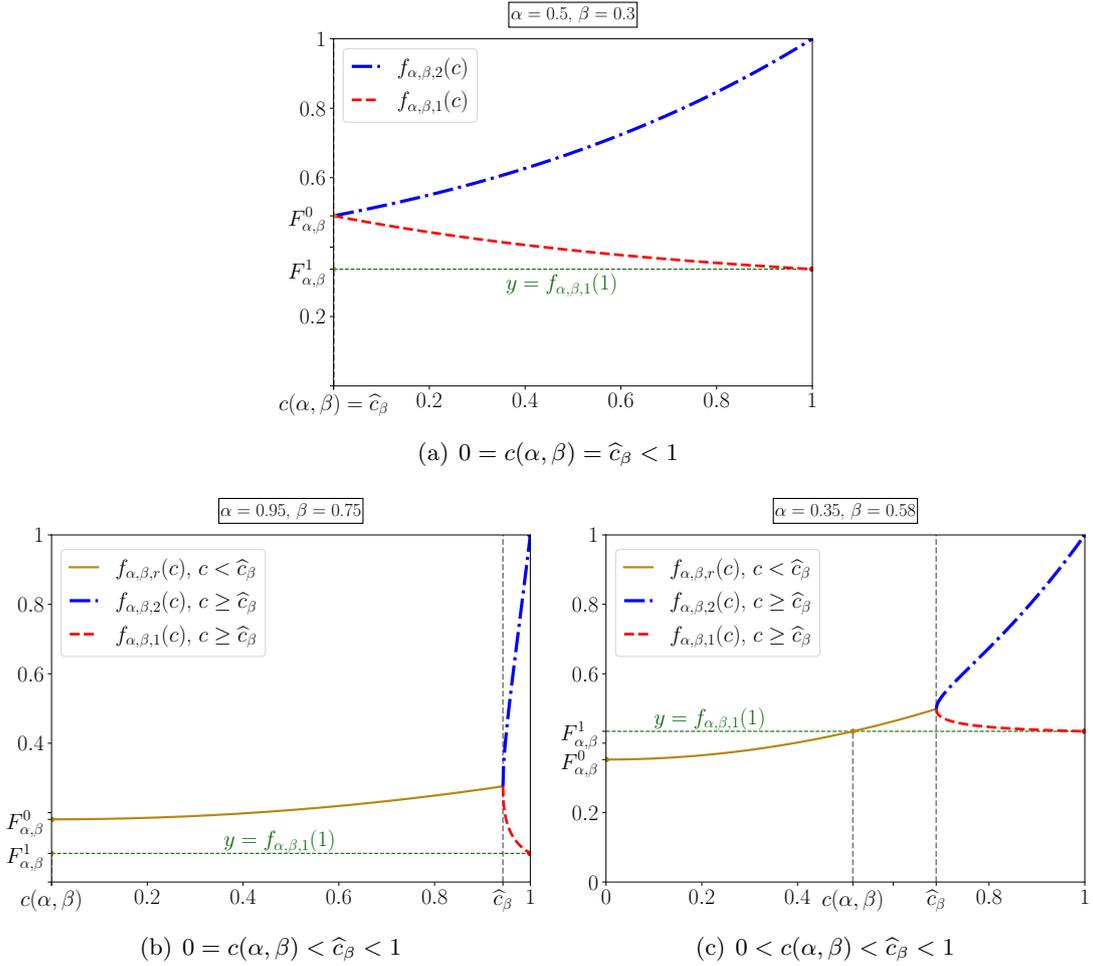

To conclude the proof, let us see that the subdominant eigenvalues are semisimple if and only if $\beta^2c_F^2-2\beta+1\neq 0$ or $c_F=0$.
The candidate eigenvalues to be subdominant are $\lambda_{1,1}$ and $\lambda_{s+1,2}$, possibly simultaneously.

Consider first the case where $\lambda_{1,1}=1-4\alpha\beta(1-\beta)$ is subdominant, and compute
\begin{equation*}
M_{\alpha,\beta}-\lambda_{1,1} I=2\alpha\beta\left(\begin{array}{cc}
-(2\beta-1)S^2 & -CS\\
(2\beta-1)CS & C^2-(2\beta-1)I_p
\end{array}\right).
\end{equation*}
Using the decomposition of $C$ and $S$ given in \eqref{decompCS}, we get
\begin{equation*}
M_{\alpha,\beta}-\lambda_{1,1} I=2\alpha\beta\left( \begin {array}{cccc}
0 & 0  & 0 & 0\\
0 & -(2\beta-1)\widetilde{S}^2  & 0 & -\widetilde{C}\widetilde{S}\\
0 & 0  & 2(1-\beta)I_s & 0\\
0 & (2\beta-1)\widetilde{C}\widetilde{S}  & 0 & \widetilde{C}^2-(2\beta-1)I_{p-s}\\
\end {array} \right),
\end{equation*}
and one has that
$\ker\left(M_{\alpha,\beta}-\lambda_{1,1} I\right)=\ker\left(\left(M_{\alpha,\beta}-\lambda_{1,1} I\right)^2\right)$ if and only if $\ker\left(M_1\right)=\ker\left(M_1^2\right)$, where $$M_1:=\left(\begin{array}{cc}
-(2\beta-1)\widetilde{S}^2 & -\widetilde{C}\widetilde{S}\\
(2\beta-1)\widetilde{C}\widetilde{S} & \widetilde{C}^2-(2\beta-1)I_{p-s}
\end{array}\right).$$
Since we are assuming that $\lambda_{1,1}$ is subdominant, it necessarily holds that $\frac{1}{2}<\beta<1$. Hence,
\begin{equation*}
\det\left(M_1\right)=\det\left((2\beta-1)^2\widetilde{S}^2\right)\neq 0,
\end{equation*}
and one trivially has that $\ker\left(M_1\right)=\ker\left(M_1^2\right)$, which proves that $\lambda_{1,1}$ is semisimple.

Consider now the case where $\lambda_{s+1,2}=1+2\alpha\beta\left(\beta c_F^2-1+c_F\sqrt{\beta^2c_F^2-2\beta+1}\right)$ is a subdominant eigenvalue. Denote by $\Delta_F:=\sqrt{\beta^2c_F^2-2\beta+1}$ and compute
\begin{equation*}
M_{\alpha,\beta}-\lambda_{s+1,2}I=2\alpha\beta\left(\begin{array}{cc}
(2\beta-1)C^2-c_F(\beta c_F+\Delta_F)I_p & -CS\\
(2\beta-1)CS&  C^2-c_F(\beta c_F+\Delta_F)I_p
\end{array}\right).
\end{equation*}
Let $k\in\{1,\ldots,p-s\}$ be such that $c_F=c_{s+1}=c_{s+2}=\cdots=c_{s+k}>c_{s+k+1}$. Then
\begin{equation*}
C=\left( \begin{array}{ccc} I_s & 0 & 0 \\ 0 & c_FI_k & 0\\ 0 & 0 & \widetilde{C}\end{array}\right)\quad\text{and}\quad S=\left( \begin{array}{ccc} 0_s & 0 & 0 \\ 0 & s_FI_k & 0\\ 0 & 0 & \widetilde{S}\end{array}\right),
\end{equation*}
where both $\widetilde{C}$ and $\widetilde{S}$ are diagonal matrices and $\widetilde{C}$ has entries strictly smaller than $c_F$. Hence, one has
\begin{equation*}
M_F:=M_{\alpha,\beta}-\lambda_{s+1,2}I=2\alpha\beta\left(\begin{array}{cccccc}
m_1 & 0 & 0 & 0 & 0 & 0 \\
0 & m_2 & 0 & 0 & m_{25} & 0\\
0 & 0 &  m_3 & 0 & 0 & m_{36}\\
0 & 0 & 0 & m_4 & 0 & 0 \\
0 & m_{52} & 0 & 0 & m_5 & 0 \\
0 & 0 & m_{63} & 0 & 0 & m_6
\end{array}\right),
\end{equation*}
where $m_{1}:=(2\beta-1-c_F(\beta c_F+\Delta_F))I_s$, $m_2:=-c_F(\Delta_F+(1-\beta)c_F)I_k$, $m_3:=(2\beta-1)\widetilde{C}^2-c_F(\beta c_F+\Delta_F)I_{p-k-s}$, $m_4:=(1-c_F(\beta c_F+\Delta_F))I_s$, $m_5:=-c_F(\Delta_F-(1-\beta)c_F)I_{k}$, $m_6:=\widetilde{C}^2-c_F(\beta c_F+\Delta_F)I_{p-k-s}$,
$m_{25}:=-c_Fs_FI_k$, $m_{36}:=-\widetilde{C}\widetilde{S}$, $m_{52}:=(2\beta-1)c_Fs_FI_k$ and $m_{63}:=(2\beta-1)\widetilde{C}\widetilde{S}$. Thus, if we denote by $M_{\{2,5\}}:=\begin{pmatrix}
m_2 & m_{25} \\
m_{52}& m_5 \\
\end{pmatrix}$ and by $M_{\{3,6\}}:=\begin{pmatrix}
m_3 & m_{36} \\
m_{63}& m_6 \\
\end{pmatrix}$, we get that
$$\ker\left(M_F\right)=\ker\left(M_F^2\right) \iff\ker\left(M_{\{2,5\}}\right) =\ker\left(M_{\{2,5\}}^2\right) \text{ and } \ker\left(M_{\{3,6\}}\right)=\ker\left(M_{\{3,6\}}^2\right).$$

On the one hand, by the block determinant formula we have that
\begin{align*}
\det\left( M_{\{3,6\}} \right)&=\det(m_3m_6-m_{63}m_{36})\\
&=\det\left((2\beta-1)\widetilde{C}^4+c_F^2(\beta c_F+\Delta_F)^2I_{p-k-s}-2\beta c_F(\beta c_F+\Delta_F)\widetilde{C}^2+(2\beta-1)\widetilde{C}^2\widetilde{S}^2\right)\\
&=\det\left((2\beta-1-2\beta c_F(\beta c_F+\Delta_F))\widetilde{C}^2+c_F^2(\beta c_F+\Delta_F)^2I_{p-k-s}\right)\\
&=\det\left((-(\beta^2c_F^2-2\beta+1)-\beta^2 c_F^2-2\beta c_F\Delta_F)\widetilde{C}^2+c_F^2(\beta c_F+\Delta_F)^2I_{p-k-s}\right)\\
&=\det\left(-\left(\beta c_F+\Delta_F\right)^2\left(\widetilde{C}^2-c_F^2I_{p-s-k}\right) \right).
\end{align*}
Observe that $\beta c_F+\Delta_F=0$ if and only if $\beta=\frac{1}{2}$ and $c_F=0$, in which case $M_{\{3,6\}}=0_{2p-k-s}$.
If $\beta c_F+\Delta_F\neq 0$, then $\det\left( M_{\{3,6\}} \right)\neq 0$. Thus, in either case, we get $\ker\left(M_{\{3,6\}}\right)=\ker\left(M_{\{3,6\}}^2\right)$.

On the other hand, we can rewrite
\begin{equation*}
M_{\{2,5\}}=-c_F\left(\begin{array}{cc}
(\Delta_F+(1-\beta)c_F)I_k & s_FI_k \\
-(2\beta-1)s_FI_k & (\Delta_F-(1-\beta)c_F)I_{k} \\
\end{array}\right),
\end{equation*}
and one has
\begin{equation*}
M_{\{2,5\}}^2=c_F^2\left(\begin{array}{cc}
\left( (\Delta_F+(1-\beta)c_F)^2-(2\beta-1)s_F^2 \right)I_k & 2\Delta_Fs_F I_k \\
-2\Delta_F(2\beta-1) s_F I_k & \left( (\Delta_F-(1-\beta)c_F)^2-(2\beta-1)s_F^2 \right)I_k\\
\end{array}\right).
\end{equation*}
Observing that
\begin{align*}
(\Delta_F-(1-\beta)c_F)^2-(2\beta-1)s_F^2 &= 2\Delta_F\left(\Delta_F-(1-\beta)c_F \right),\\
(\Delta_F+(1-\beta)c_F)^2-(2\beta-1)s_F^2 &= 2\Delta_F\left(\Delta_F+(1-\beta)c_F\right),
\end{align*}
we deduce that $M_{\{2,5\}}^2=-2\Delta_Fc_FM_{\{2,5\}}$. If $c_F=0$, then $M_{\{2,5\}}=0_{2k}$. Therefore,
$\ker\left(M_{\{2,5\}}\right) =\ker\left(M_{\{2,5\}}^2\right)$ if and only if $\Delta_F\neq 0$ or $c_F=0$.

Summarizing the discussion above, we have shown that
$$\ker\left(M_F\right)=\ker\left(M_F^2\right) \iff \Delta_F\neq 0\text{ or }c_F=0.$$
Finally, observe that $\Delta_F=0$ if and only if $c_F=\widehat{c}_\beta$, in which case $\lambda_{s+1,1}$ is a subdominant eigenvalue, and this proves the last assertion in the statement.

Case 2: $p+q\geq n$. We can take some $k\geq 1$ such that $n':=n+k>p+q$, and consider $U':=U\times\{0_{k\times 1}\}\subset\R^{n'}, V':=V\times\{0_{k\times 1}\}\subset\R^{n'}$, and $T'_{\alpha,\beta}:=T_{U',V',\alpha,\beta}=(1-\alpha)I+\alpha (2\beta P_{V'}-I)(2\beta P_{U'}-I)$. Since $P_{U'}=\left(\begin{array}{cc} P_U & 0 \\ 0 & 0_k\end{array}\right)$ and $P_{V'}=\left(\begin{array}{cc} P_V & 0 \\ 0 & 0_k\end{array}\right)$, it holds that
$$
T'_{\alpha,\beta}=\left(\begin{array}{cc} T_{\alpha,\beta} & 0 \\ 0 & I_k\end{array}\right).
$$

Therefore, $\sigma(T_{\alpha,\beta})\cup\{1\}=\sigma(T'_{\alpha,\beta})$ and $\gamma(T_{\alpha,\beta})=\gamma(T'_{\alpha,\beta})$. Note that the principal angles between $U'$ and $V'$ are the same that the ones between $U$ and $V$. Hence, the result follows from applying Case~1 to $T'_{\alpha,\beta}$.
\end{proof}

\begin{remark}
For simplicity, we have assumed that $\dim U=p\leq q=\dim V$. If this is not the case and $q<p$, observe that one has to exchange the matrix decomposition of $P_U$ and $P_V$ given in \eqref{eq:projections_principal_angles}. In this case, one can check that the matrix $T_{\alpha,\beta}$ obtained corresponds to the transpose of the one given in \eqref{eq:T_matrix}. Hence, the spectrum of $T_{\alpha,\beta}$ remains the same and thus all the results in \Cref{th:rate} also hold.
\end{remark}

\begin{remark}
The expression in~\eqref{eq:rate} corroborates what it was numerically observed in~\cite{AAMR}: there are values of $\alpha$ and $\beta$ for which the rate of convergence of AAMR does not depend on the value of the Friedrichs angle for all angles larger than $\arccos c(\alpha,\beta)$.
\end{remark}

We now look for the values of the parameters $\alpha$ and $\beta$ and $V$, in order to
that minimize the rate of convergence of the AAMR method obtained in \Cref{th:rate}.

\begin{theorem}\label{th:optim_rate}
The infimum of the linear convergence rates of the AAMR operator $T_{\alpha,\beta}$ attains its smallest value at $\alpha^\star=1$ and $\beta^\star=\frac{1}{1+\sin\theta_F}$, where $\theta_F$ is the Friedrichs angle between $U$ and $V$; i.e., it holds
\begin{equation*}
\frac{1-\sin\theta_F}{1+\sin\theta_F}=\gamma\left( T_{1,\beta^\star}\right)\leq \gamma\left(T_{\alpha,\beta}\right) \quad \text{ for all } (\alpha,\beta)\in\,]0,1]\times\,]0,1[\,.
\end{equation*}
Furthermore, $\gamma\left( T_{1,\beta^\star}\right)$ is an optimal linear convergence rate if and only if $\theta_F=\frac{\pi}{2}$.
\end{theorem}
\begin{proof}
Let us look for the values of parameters $\alpha$ and $\beta$ that minimize the rate $\gamma(T_{\alpha,\beta})$ given by~\eqref{eq:rate}. Define the sets $D:=\,]0,1]\times\,]0,1[\,$,
\begin{align*}
D_1:=&\left\{ (\alpha,\beta)\in D: \beta<\frac{1}{1+s_F}\right\},\\
D_2:=&\left\{ (\alpha,\beta)\in D: \frac{1}{1+s_F}\leq\beta\leq\frac{1}{1+s_F^2}\text{ or } \alpha\geq\frac{1-\beta(1+s_F^2)}{\beta\left( 4(1-\beta)^2-s_F^2 \right)}, \beta>\frac{1}{1+s_F^2}  \right\},\\
D_3:=&\left\{ (\alpha,\beta)\in D: \alpha<\frac{1-\beta(1+s_F^2)}{\beta\left( 4(1-\beta)^2-s_F^2 \right)} , \beta>\frac{1}{1+s_F^2}\right\},
\end{align*}
and the functions
\begin{align*}
\Gamma_1(\alpha,\beta)&:=1+2\alpha\beta\left(\beta c_F^2-1+c_F\sqrt{\beta^2 c_F^2-2\beta+1}\right),\quad\text{for }(\alpha,\beta)\in D_1,\\
\Gamma_2(\alpha,\beta)&:=\sqrt{4(1-\alpha)\alpha\beta^2 c_F^2+(1-2\alpha\beta)^2},\quad\text{for }(\alpha,\beta)\in D,\\
\Gamma_3(\alpha,\beta)&:=1-4\alpha\beta(1-\beta),\quad\text{for }(\alpha,\beta)\in D,
\end{align*}
having $D=D_1\cup D_2\cup D_3$.
Hence, we can define the convergence rate in terms of the parameters $\alpha$ and $\beta$ through the function
$$\Gamma(\alpha,\beta):=\gamma(T_{\alpha,\beta})=\left\{\begin{array}{ll}
\Gamma_1(\alpha,\beta), & \text{ if } (\alpha,\beta)\in D_1,\\
\Gamma_2(\alpha,\beta), & \text{ if } (\alpha,\beta)\in D_2,\\
\Gamma_3(\alpha,\beta), & \text{ if } (\alpha,\beta)\in D_3,
\end{array}\right.$$
see \Cref{fig:num_exp_alpha}.

\begin{figure}[ht!]
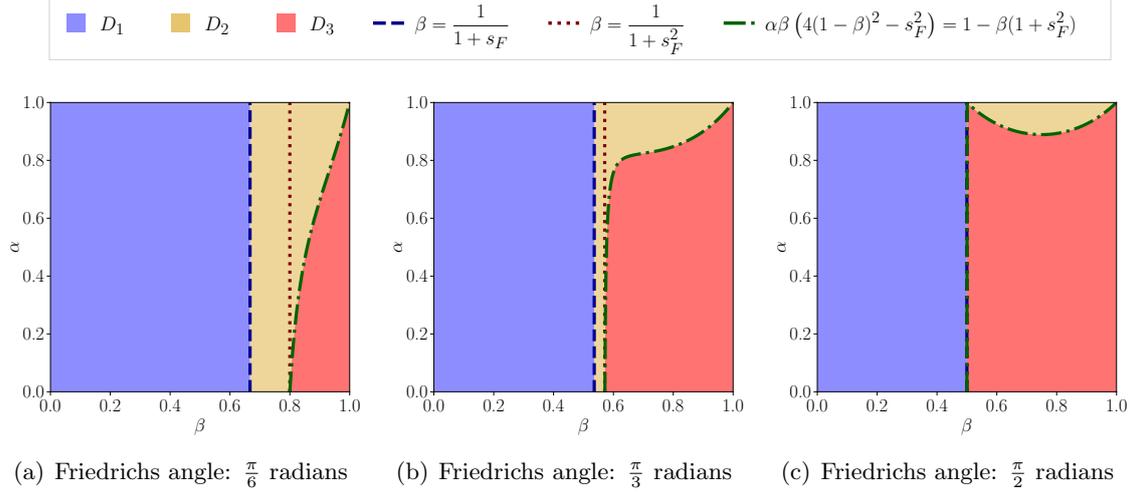
\addtocounter{subfigure}{-1}
	\centering
	\subfigure{\scalebox{.43}{\input{Fig2d.pgf}}} 	
    \subfigure[Friedrichs angle: $\frac{\pi}{6}$ radians]{\scalebox{0.4}{\input{Fig2a.pgf}}} 	 \subfigure[Friedrichs angle: $\frac{\pi}{3}$ radians]{\scalebox{0.4}{\input{Fig2b.pgf}}} 	 \subfigure[Friedrichs angle: $\frac{\pi}{2}$ radians]{\scalebox{0.4}{\input{Fig2c.pgf}}} 	 \caption{Piecewise domain of the function $\Gamma(\alpha,\beta)$ for three different values of the Friedrichs angle} \label{fig:num_exp_alpha}
\end{figure}

The function $\Gamma$ is piecewise defined, continuous and differentiable on the interior of each of the three regions $D_1$, $D_2$ and $D_3$, but is not differentiable on the boundaries. Let us analyze the three problems of minimizing the function $\Gamma$ over the closure of each of the three pieces. The gradient of the functions $\Gamma_1$, $\Gamma_2$ and $\Gamma_3$ are given by
\begin{gather*}
\nabla \Gamma_1(\alpha,\beta)=\begin{pmatrix}
2\beta\left(\beta c_F^2-1+c_F\sqrt{\beta^2c_F^2-2\beta+1}\right)\\
2\alpha\left(\beta c_F^2-1+c_F\sqrt{\beta^2c_F^2-2\beta+1}\right)\left(\frac{\beta c_F+\sqrt{\beta^2c^2_F-2\beta+1}}{\sqrt{\beta^2c^2_F-2\beta+1}}\right)\\
\end{pmatrix},\\
\nabla \Gamma_2(\alpha,\beta)=\frac{1}{\sqrt{4(1-\alpha)\alpha\beta^2 c_F^2+(1-2\alpha\beta)^2}}\begin{pmatrix}
2\beta\left(\beta c_F^2-1+2\alpha\beta(1-c_F^2)\right)\\
2\alpha\left(2\beta c_F^2-1+2\alpha\beta(1-c_F^2)\right)\\
\end{pmatrix},\\
\nabla \Gamma_3(\alpha,\beta)=\begin{pmatrix}
-4\beta\left(1-\beta\right)\\
-4\alpha\left(1-2\beta\right)\\
\end{pmatrix}.
\end{gather*}

To minimize $\Gamma$ over $\overline{D_1}$, we assert that
\begin{equation}\label{eq:gradD1}
\frac{\partial\Gamma_1}{\partial\alpha}(\alpha,\beta)< 0, \frac{\partial\Gamma_1}{\partial\beta}(\alpha,\beta)< 0, \quad \text{for all } (\alpha,\beta)\in D_1.
\end{equation}
Indeed, on the one hand, since $c_F^2<1$, then
$$c_F^2\left( \beta^2c_F^2-2\beta+1 \right)=\beta^2c_F^4-2\beta c_F^2+c_F^2<\beta^2c_F^4-2\beta c_F^2+1=\left(1-\beta c_F^2\right)^2.$$
Thus, taking square roots and reordering, we get that $\beta c_F^2-1+c_F\sqrt{\beta^2c_F^2-2\beta+1}<0$. On the other hand, one has $\beta c_F+\sqrt{\beta^2c^2_F-2\beta+1}\geq 0$, with equality if and only if $c_F=0$ and $\beta=\frac{1}{2}$. In this case, the point has the form $\left(\alpha,\frac{1}{2}\right)\notin D_1$ for $\alpha\in\,]0,1]$. We have therefore shown that~\eqref{eq:gradD1} holds, and thus the unique minimum of $\Gamma$ over $\overline{D_1}$ is attained at $\left(1,\frac{1}{1+s_F}\right)$.

Let us consider now the problem of minimizing $\Gamma$ over $\overline{D_2}$.  To address this problem, we consider two cases. Suppose first that $c_F=0$ and observe that
$$\Gamma_2(\alpha,\beta)=\sqrt{(1-2\alpha\beta)^2}\geq 0, \quad\text{for all } (\alpha,\beta)\in\overline{D_2},$$
having $\Gamma_2(\alpha,\beta)=0$ if and only if $2\alpha\beta=1$. Since $\left(1,\frac{1}{2}\right)$ is the only point in $\overline{D_2}$ satisfying this equation,  we deduce that it is the unique minimum.

Suppose now that $c_F>0$. In this case, we claim that $\Gamma_2$ attains its minimimum over the region $\overline{D_2\cup D_3}=[0,1]\times \left[\frac{1}{1+s_F},1\right]$ at the point $\left(1,\frac{1}{1+s_F}\right)\in {D_2}$, and so does $\Gamma$ over $\overline{D_2}$.  Indeed, observe that $\Gamma_2$ is smooth on the interior of the set ${D_2\cup D_3}$. Moreover, $\nabla \Gamma_2$ only vanishes at $(0,0)$. Therefore, the minimum has to be attained at some point in the boundary. Note that, for all $\beta\in\left[\frac{1}{1+s_F},1\right]$, the following holds:

\begin{itemize}[noitemsep]
\item[(i)]$\Gamma_2(0,\beta)= 1.$
\item[(ii)]$\Gamma_2(1,\beta)= 2\beta-1$, which attains its minimum at $\beta=\frac{1}{1+s_F}.$
\item[(iii)] The function $\alpha\mapsto\Gamma_2(\alpha,\beta)$ is the square root of a positive non-degenerated convex parabola,
\begin{equation*}
\alpha\mapsto\Gamma_2(\alpha,\beta)=\sqrt{4\beta^2s_F^2\alpha^2-4\beta(1-\beta c_F^2)\alpha+1};
\end{equation*}
which attains its minimum at $\alpha^\star(\beta):=\frac{1-\beta c_F^2}{2\beta s_F^2}$. Since $\alpha^\star\left(\frac{1}{1+s_F}\right)=\frac{1+s_F}{2s_F}\geq 1$, we have
\begin{equation*}
\Gamma_2\left(1,\frac{1}{1+s_F}\right)<\Gamma_2\left(\alpha,\frac{1}{1+s_F}\right),\quad \text{ for all } \alpha\in\left[0,1\right].
\end{equation*}
On the other hand, $\alpha^\star(1)=\frac{1}{2}$, which implies
\begin{equation*}
\Gamma_2\left(\frac{1}{2},1\right)<\Gamma_2(\alpha,1), \quad\text{ for all } \alpha\in\left[0,1\right].
\end{equation*}
\end{itemize}
Then, noting that
\begin{equation*}
\Gamma_2\left(1,\frac{1}{1+s_F}\right)=\frac{1-s_F}{1+s_F}<c_F=\Gamma_2\left(\frac{1}{2},1\right),
\end{equation*}
we have shown by (i)--(iii) that $\Gamma_2$ attains its minimum over $\overline{D_2\cup D_3}$ at $\left(1,\frac{1}{1+s_F}\right)\in {D_2}$, as claimed.

Finally, observe that if $(\alpha,\beta)\in D_3$, it holds that $\frac{1}{2}<\beta<1$. Then
\begin{equation*}
\frac{\partial\Gamma_3}{\partial\alpha}(\alpha,\beta)< 0, \frac{\partial\Gamma_3}{\partial\beta}(\alpha,\beta)>0, \quad \text{for all } (\alpha,\beta)\in D_3.
\end{equation*}
Thus, there exists some point $(\alpha^\star_3,\beta^\star_3)\in \overline{D_3}$ with $\alpha^\star_3\beta^\star_3\left( 4(1-\beta^\star_3)^2-s_F^2  \right)=1-\beta^\star_3(1+s_F^2)$ such that
\begin{equation*}
\Gamma_3(\alpha^\star_3,\beta^\star_3) < \Gamma_3(\alpha,\beta), \quad\text{ for all } (\alpha,\beta)\in{D_3}.
\end{equation*}
Note that $(\alpha^\star_3,\beta^\star_3)$ lies on the boundary curve between $D_2$ and $D_3$. Since $\Gamma$ is continuous on $D$, it holds that $\Gamma_2(\alpha^\star_3,\beta^\star_3)=\Gamma_3(\alpha^\star_3,\beta^\star_3)$ and hence,
\begin{equation*}
\Gamma\left(1,\frac{1}{1+s_F}\right)\leq\Gamma(\alpha^\star_3,\beta^\star_3) < \Gamma(\alpha,\beta),  \text{ for all } (\alpha,\beta)\in{D_3}.
\end{equation*}

Hence, all the reasoning above proves that
\begin{equation*}
\argmin_{(\alpha,\beta)\in D} \Gamma(\alpha,\beta)=\left(1,\frac{1}{1+s_F}\right),
\end{equation*}
with $\Gamma\left(1,\frac{1}{1+s_F}\right)=\frac{1-s_F}{1+s_F}$. Finally, by the last assertion in \Cref{th:rate}, $\gamma\left(T_{1,\frac{1}{1+s_F}}\right)$ is an optimal linear convergence rate if and only if $c_F=0$, as claimed.
\end{proof}

\begin{corollary}\label{cor:convergent}
Let $\alpha\in\,]0,1]$ and $\beta\in\,]0,1[\,$. Given $z\in\R^n$, choose any $x_0\in\R^n$ and consider the sequence generated, for $k=0,1,2,\ldots$, by
\begin{equation*}\label{eq:seq_xk}
x_{k+1}=T_{U-z,V-z,\alpha,\beta} (x_k)=(1-\alpha)x_k+\alpha(2\beta P_{V-z}-I)(2\beta P_{U-z} - I)(x_k).
\end{equation*}
Let $\gamma(T_{\alpha,\beta})$ be given by~\eqref{eq:rate}. Then, for every $\mu\in\,]\gamma(T_{\alpha,\beta}),1[\,$, the sequence $\left(x_k\right)_{k\geq 0}$ is R-linearly convergent to $P_{\Fix T_{U-z,V-z,\alpha,\beta}}(x_0)$ and the shadow sequence $\left( P_U(z+x_k)\right)_{k\geq 0}$ is
R-linearly convergent to $P_{U\cap V}(z)$, both with rate $\mu$, in the sense that
there exists a positive integer $k_0$ such that
\begin{equation}\label{eq:k_0}
\left\|P_U(z+x_k)-P_{U\cap V}(z)\right\|\leq\left\|x_k-P_{\Fix T_{U-z,V-z,\alpha,\beta}}(x_0)\right\|\leq \mu^k,\quad\text{for all }k\geq k_0.
\end{equation}
\end{corollary}
\begin{proof}
According to \Cref{fact:translation_aamr} we have that $\Fix T_{U-z,V-z,\alpha,\beta}\neq \emptyset$ and
$$x_{k+1}=T_{U-z,V-z,\alpha,\beta} (x_k)=T_{U,V,\alpha,\beta}(x_k-x^*)+x^*,$$
for $x^*:=P_{\Fix T_{U-z,V-z,\alpha,\beta}}(x_0)$.
Hence, one has
\begin{align*}
\|x_{k}-x^*\|&=\|T_{U,V,\alpha,\beta}(x_{k-1}-x^*)\|=\cdots=\|T_{U,V,\alpha,\beta}^k(x_0-x^*)\|.
\end{align*}
Again by \Cref{fact:translation_aamr}, one has
$$
\Fix T_{U-z,V-z,\alpha,\beta}=x^*+U^\perp\cap V^\perp,
$$
and by the translation formula for projections (c.f.~\cite[2.7(ii)]{D01}),
$$x^*=P_{\Fix T_{U-z,V-z,\alpha,\beta}}(x_0)=P_{x^*+U^\perp\cap V^\perp}(x_0)=P_{U^\perp\cap V^\perp}(x_0-x^*)+x^*,$$
which implies $P_{U^\perp\cap V^\perp}(x_0-x^*)=0$, and therefore,
$$\|x_{k}-x^*\|=\left\|T_{U,V,\alpha,\beta}^k(x_0-x^*)\right\|=\left\|\left(T_{U,V,\alpha,\beta}^k-P_{U^\perp\cap V^\perp}\right)(x_0-x^*)\right\|.$$
Let $\nu\in\,]\gamma(T_{\alpha,\beta}),\mu[$. Since $\nu>\gamma(T_{\alpha,\beta})$, by \Cref{th:rate}, there exists a positive integer $k_1$ and some $M>0$ such that
$$\left\|T_{U,V,\alpha,\beta}^k-P_{U^\perp\cap V^\perp}\right\|\leq M\nu^k,\quad \text{for all } k\geq k_1.$$
Let $k_0\geq k_1$ be a positive integer such that
$$\left(\frac{\mu}{\nu}\right)^k\geq M\|x_0-x^*\|, \quad \text{for all } k\geq k_0.$$
Then, we deduce that
\begin{align*}
\|x_{k}-x^*\|\leq \left\|T_{U,V,\alpha,\beta}^k-P_{U^\perp\cap V^\perp}\right\|\|x_0-x^*\|\leq M\nu^k\|x_0-x^*\|\leq \mu^k,
\end{align*}
for all $k\geq k_0$, which proves the second inequality in~\eqref{eq:k_0}.

By~\cite[Proposition~3.4]{AAMR} and the translation formula for projections (c.f.~\cite[2.7(ii)]{D01}), we can deduce that
\begin{equation*}
P_U\left(z+P_{\Fix T_{U-z,V-z,\alpha,\beta}}(x_0)\right)=P_{U\cap V}(z).
\end{equation*}
Thus, the first inequality in~\eqref{eq:k_0} is a consequence of this, and the linearity and nonexpansiveness of $P_U$. Indeed
\begin{align*}
\|P_U(z+x_k)-P_{U\cap V}(z)\|&=\|P_U(z+x_k)-P_U(z+P_{\Fix T_{U-z,V-z,\alpha,\beta}}(x_0))\|\\
&=\|P_U(x_k-P_{\Fix T_{U-z,V-z,\alpha,\beta}}(x_0))\|\\
&\leq \|x_k-P_{\Fix T_{U-z,V-z,\alpha,\beta}}(x_0)\|,
\end{align*}
which completes the proof.
\end{proof}

\section{Comparison with other projection methods}\label{sec:comparison}

In this section, we compare the rate of AAMR with optimal parameters obtained in \Cref{sec:rate} with the rates of various projection methods analyzed in~\cite{BCNPW15,FG17}.
We summarize the key features of these schemes in \Cref{tbl:Comparison_rates}, where we recall the operator defining the iteration of each method, as well as
the optimal parameters and rates of convergence when these schemes are applied to linear subspaces. Note that all these rates only depend on the Friedrichs angle $\theta_F$
between the subspaces.

\begin{table}[ht!]
\centering
\bgroup
\def\arraystretch{1.4}
\begin{tabular}{|c|c|c|}
\hline
Method & Optimal parameter(s) & Rate \tabularnewline
\hhline{===}
Alternating Projections &\multirow{2}{*}{--} & \multirow{2}{*}{$\cos^2\theta_F$} \tabularnewline
$AP=P_VP_U$ & &  \tabularnewline
\hline
Relaxed Alternating Projections &\multirow{2}{*}{$\alpha^\star=\frac{2}{1+\sin^2\theta_F}$} & \multirow{2}{*}{$\frac{1-\sin^2\theta_F}{1+\sin^2\theta_F}$} \tabularnewline
$RAP=(1-\alpha)I+\alpha P_VP_U$ & &  \tabularnewline
\hline
Generalized Alternating Projections &\multirow{1}{*}{$\alpha^\star=1$} & \multirow{2}{*}{$\frac{1-\sin\theta_F}{1+\sin\theta_F}$} \bigstrut[t] \tabularnewline
$GAP=(1-\alpha)I+\alpha\left(\alpha_1 P_V+(1-\alpha_1)I\right)\left(\alpha_2 P_U+(1-\alpha_2)I\right)$ & $\alpha_1^\star=\alpha_2^\star=\frac{2}{1+\sin\theta_F}$ & \bigstrut[b] \tabularnewline
\hline
Douglas--Rachford &\multirow{2}{*}{$\alpha^\star=\frac{1}{2}$} & \multirow{2}{*}{$\cos\theta_F$} \tabularnewline
$DR=(1-\alpha)I+\alpha(2P_V-I)(2P_U-I)$ & &  \tabularnewline
\hline
Averaged Alternating Modified Reflections &\multirow{2}{*}{$\alpha^\star=1$, $\beta^\star=\frac{1}{1+\sin\theta_F}$} & \multirow{2}{*}{$\frac{1-\sin\theta_F}{1+\sin\theta_F}$} \tabularnewline
$AAMR=(1-\alpha)I+\alpha(2\beta P_V-I)(2\beta P_U-I)$ & &  \tabularnewline
\hline
\end{tabular}
\egroup
\caption{Rates of convergence with optimal parameters of AP, RAP, GAP, DR and AAMR when they are applied to two subspaces}\label{tbl:Comparison_rates}
\end{table}

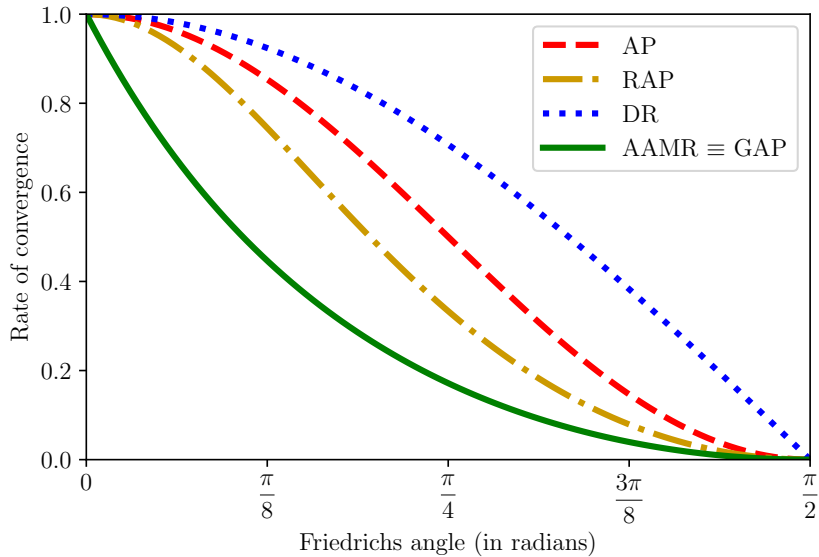
\begin{figure}[ht!]
    \centering
    \scalebox{0.77}{\input{Fig3.pgf}}
	\caption{Comparison of the rates of linear convergence with optimal parameters of AP, RAP, DR, AAMR and GAP}\label{fig:comparison_rates}
\end{figure}

On the one hand, we observe that the rates for AAMR and GAP coincide. Moreover, their optimal parameters are closely related, in the sense that
\begin{equation*}\label{eq:optGAP_AAMR}
\alpha^\star_{AAMR}=\alpha^\star_{GAP} \quad \text{and} \quad \alpha^\star_{1,GAP}=\alpha^\star_{2,GAP}=2\beta^\star_{AAMR}.
\end{equation*}
We explain this behavior in \Cref{subsec:GAP_AAMR}, where under some conditions, we show that the shadow sequences of GAP and AAMR coincide for linear subspaces (\Cref{th:GAP_AAMR}). On the other hand, we note that the rate for AAMR/GAP is considerably smaller than the one of other methods, see \Cref{fig:comparison_rates}. We numerically demonstrate this with a computational experiment in \Cref{sec:experiments}.

\subsection{Relationship between AAMR and GAP for subspaces}\label{subsec:GAP_AAMR}

F\"alt and Giselsson have recently obtained in~\cite{FG17} the rate of convergence with optimal parameters for the generalized alternating projections (GAP) method for two subspaces. This iterative scheme is defined by
\begin{equation}\label{eq:GAP}
z_{k+1}=(1-\alpha)z_k+\alpha\left(\alpha_1 P_V+(1-\alpha_1)I\right)\left(\alpha_2 P_U+(1-\alpha_2)I\right)(z_k),
\end{equation}
where $\alpha\in\,]0,1]$ and $\alpha_1,\alpha_2\in\,]0,2]$.

\begin{figure}[ht!]
    \centering
    \scalebox{.55}{\input{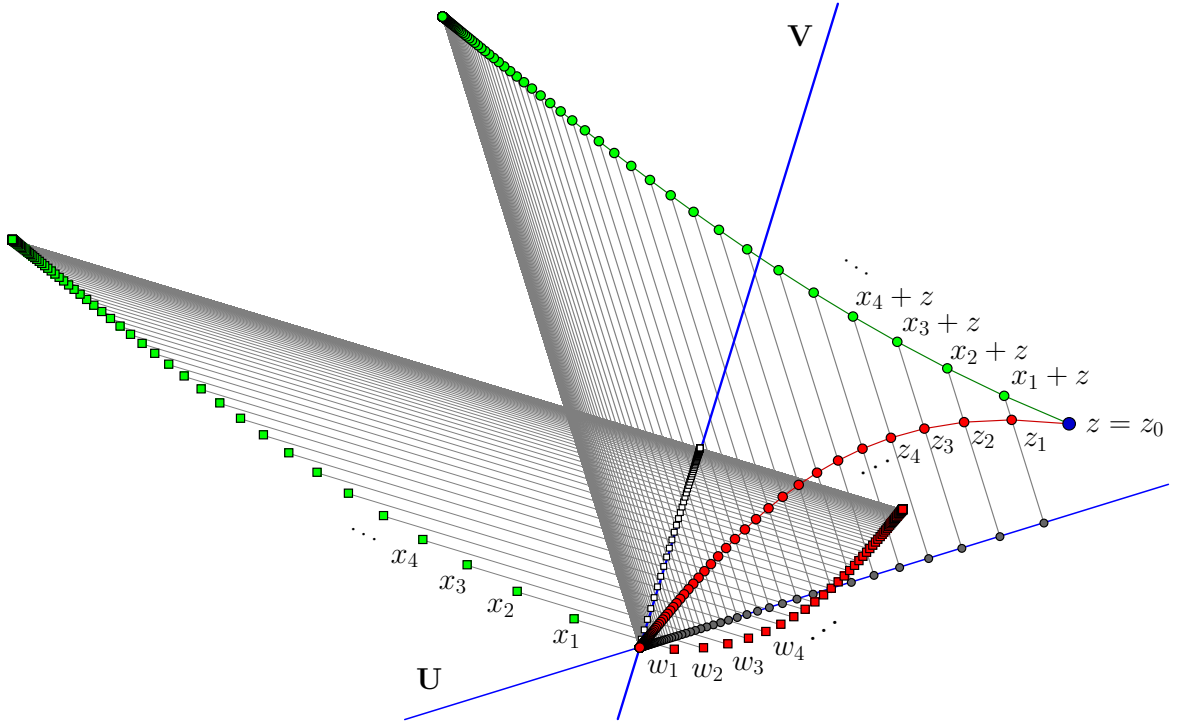}}
	\caption{Graphical representation of \Cref{th:GAP_AAMR} for two lines in $\R^2$. The sequence $\{x_k+z\}_{k=0}^\infty$ is generated by AAMR with $x_0=0$, while the sequence $\{z_k\}_{k=0}^\infty$ is generated by GAP. We also represent $w_k:=(2\beta-1)(z_k-z)$}\label{fig:th_GAP_AAMR}
\end{figure}

\begin{figure}[ht!]
    \centering
    \scalebox{0.55}{\input{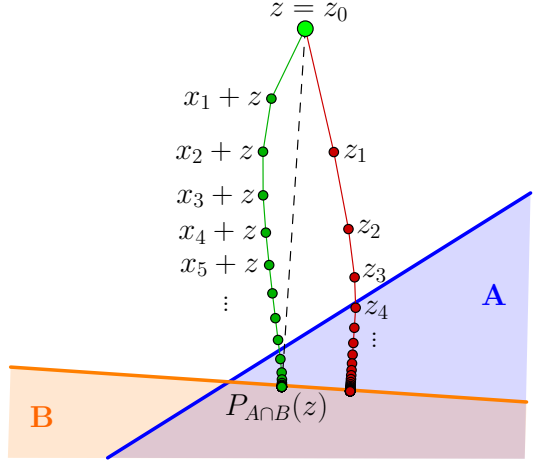}}
	\caption{The sequence $\{x_k+z\}_{k=0}^\infty$, with $x_0=0$, generated by AAMR converges to $P_{A\cap B}(z)$ and thus solves the best approximation problem, while the sequence $\{z_k\}_{k=0}^\infty$, with $z_0=z$, generated by GAP only converges to some point in $A\cap B$}\label{fig:GAP_vs_AAMR}
\end{figure}

The next result shows that, for subspaces, the shadow sequences of GAP and AAMR coincide when $\alpha_1=\alpha_2=2\beta$ and the starting point of AAMR is chosen as $x_0=0$; see \Cref{fig:th_GAP_AAMR} for a simple example in $\mathbb{R}^2$. This is not the case for general convex sets, as shown in \Cref{fig:GAP_vs_AAMR}: GAP gives a point in the intersection of the sets, while AAMR solves the best approximation problem~\eqref{eq:bestaproxprob}. ~\Cref{fig:th_GAP_AAMR,fig:GAP_vs_AAMR} were created with \texttt{Cinderella}~\cite{Cinderella}.


\begin{theorem}\label{th:GAP_AAMR}
Given $z\in\R^n$, set $x_0=0$ and consider the AAMR sequence~$\{x_k\}_{k=0}^\infty$ generated by~\eqref{eq:aamr_intro} with parameters $\alpha\in\,]0,1]$ and $\beta\in\,]0,1[$.
Let~$\{z_k\}_{k=0}^\infty$ be the sequence generated by GAP~\eqref{eq:GAP} with parameters $\alpha$ and $\alpha_1=\alpha_2=2\beta$ and starting point $z_0=z$. Then, one has
\begin{equation}\label{eq:induction}
P_U(x_k+z)=P_U(z_k)\quad\text{and}\quad P_V(x_{k})=(2\beta-1)P_V(z_k-z),
\end{equation}
for all $k=0, 1, 2,\ldots$.
\end{theorem}

\begin{proof}
To simplify the notation, let $\eta:=2\beta$. We shall prove~\eqref{eq:induction} by induction. Since both equalities clearly hold for $k=0$, we can assume that they are valid for some $k\geq 0$.
By~\eqref{eq:AAMR_subspaces}, the sequence generated by the AAMR scheme satisfies
\begin{align}
P_U(x_{k+1})&=P_U(x_k)+\alpha\eta\big( \eta P_UP_{V}P_{U}(x_k+z) \nonumber\\
&\quad+(1-\eta)P_UP_V(z) -P_UP_V(x_k)-P_{U}(x_k+z)\big)\nonumber\\
&=\left(\alpha\eta^2P_UP_{V}P_{U}-\alpha\eta P_UP_V+(1-\alpha\eta)P_U\right)(x_k)\nonumber\\
&\quad+\alpha\left(\eta^2P_UP_{V}P_{U}+\eta(1-\eta)P_UP_V-\eta P_U\right)(z)\label{eq:02},
\end{align}
and,
\begin{align}
P_V(x_{k+1})&=P_V(x_k)+\alpha\eta\left((\eta-1) P_VP_U(x_k+z)+P_V((1-\eta)z-x_k)\right)\nonumber\\
&=\alpha\eta(\eta-1)P_VP_U(x_k+z)+P_V\left(\alpha\eta(1-\eta)z+(1-\alpha\eta)x_k\right).\label{eq:03}
\end{align}
Thanks to the linearity of the projectors onto subspaces and using $\alpha_1=\alpha_2=\eta$, the GAP iteration~\eqref{eq:GAP} takes the form
\begin{equation*}\label{eq:GAP_subspaces}
z_{k+1}=\left(1-\alpha\eta(2-\eta)\right)z_k+\alpha\left(\eta^2P_VP_Uz_k+\eta(1-\eta)P_Vz_k+\eta(1-\eta)P_Uz_k\right);
\end{equation*}
and thus this scheme verifies
\begin{equation}\label{eq:00}
P_U(z_{k+1})=\left(\alpha\eta^2P_UP_VP_U+\alpha\eta (1-\eta)P_UP_V+(1-\alpha\eta)P_U\right)(z_k),
\end{equation}
and
\begin{equation}\label{eq:01}
P_V(z_{k+1})=(1-\alpha\eta)P_V(z_k)+\alpha\eta P_VP_U(z_k).
\end{equation}
Then, by~\eqref{eq:03}, the induction hypothesis~\eqref{eq:induction} and~\eqref{eq:01}, we obtain
\begin{align*}
P_V(x_{k+1})&=\alpha\eta(\eta-1)P_VP_U(z_k)+(1-\alpha\eta)P_V(x_k)+P_V(\alpha\eta(1-\eta)z)\\
&=\alpha\eta(\eta-1)P_VP_U(z_k)+(1-\alpha\eta)(\eta-1)P_V(z_k-z)+P_V(\alpha\eta(1-\eta)z)\\
&=(\eta-1)\left(\alpha\eta P_VP_U(z_k)+(1-\alpha\eta)P_V(z_k)\right)+(1-\eta)P_V(z)\\
&=(\eta-1)P_V(z_{k+1}-z),
\end{align*}
which proves the second equation in~\eqref{eq:induction} for $k+1$. Finally,  by~\eqref{eq:00},~\eqref{eq:induction} and~\eqref{eq:02}, we have that
\begin{align*}
P_U(z_{k+1})&=\alpha\eta^2P_UP_VP_U(x_k+z)+\alpha\eta P_UP_V(-x_k+(1-\eta)z)+(1-\alpha\eta)P_U(x_k+z)\\
&=\left(\alpha\eta^2P_UP_VP_U-\alpha\eta P_UP_V+(1-\alpha\eta)P_U\right)(x_k)\\
&\quad+\alpha\left(\eta^2P_UP_VP_U+\eta(1-\eta)P_UP_V-\eta P_U\right)(z)+P_U(z)\\
&=P_U(x_{k+1}+z),
\end{align*}
which proves the first equation in~\eqref{eq:induction} for $k+1$ and completes the proof.
\end{proof}

\section{Computational experiments}\label{sec:experiments}

In this section we demonstrate the theoretical results obtained in the previous sections with two different numerical experiments. In both experiments we consider randomly generated subspaces $U$ and $V$ in $\R^{50}$ with $U\cap V\neq\{0\}$. We have implemented all the algorithms in \emph{Python 2.7} and the figures were drawn with \emph{Matplotlib}~\cite{matplotlib}.

\begin{figure}[ht!]
    \centering
   \scalebox{0.81}{\input{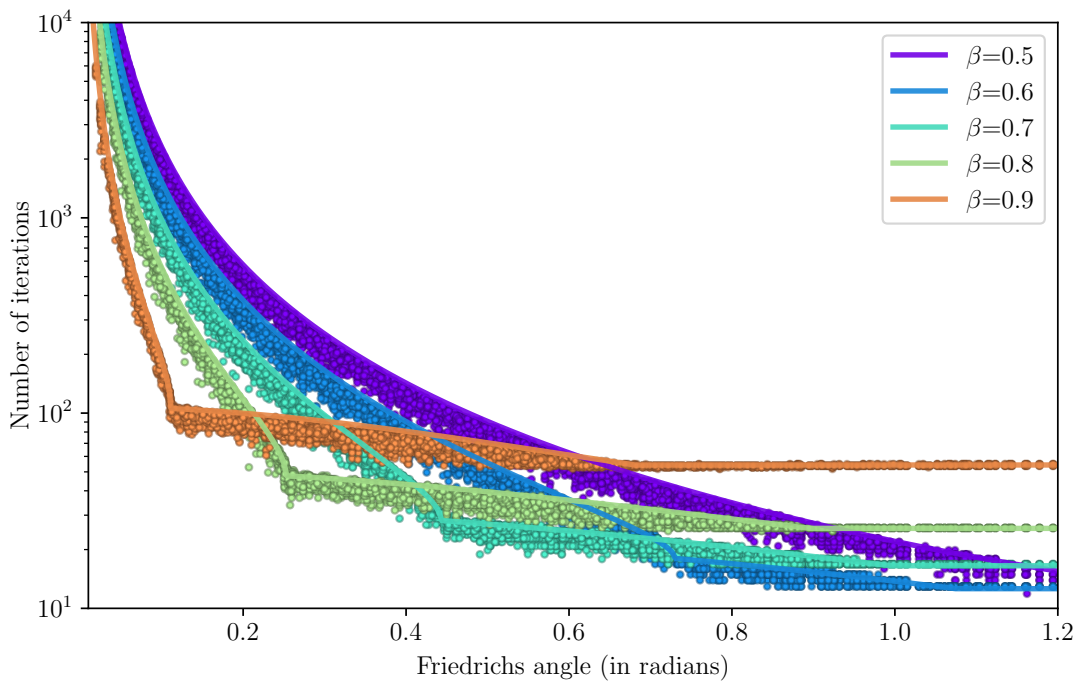}}
	\caption{Number of iterations required to converge for the AAMR algorithm with $\alpha=0.8$ and five different values of the parameter~$\beta$, with respect to the Friedrichs angle. The lines correspond to the approximate upper bounds given by~\eqref{eq:maxiter} and the theoretical rates~\eqref{eq:rate}}	
\label{fig:theo_rate}
\end{figure}

The purpose of our first computational experiment is to exhibit the piecewise expression of the convergence rate $\gamma(T_{\alpha,\beta})$ given in~\Cref{th:rate}. To this aim, we generated $500$ pairs of random subspaces. For each pair of subspaces, we chose $10$ random starting points with $\|x_0\|=1$.
Then, for each of these instances, we ran the AAMR method with $\alpha=0.8$ and $\beta\in\{0.5, 0.6, 0.7, 0.8, 0.9\}$. The algorithm was stopped when the shadow sequence satisfies
\begin{equation}\label{eq:stopping_criterion}
\|P_U(x_n+x_0)-P_{U\cap V}(x_0)\|<\epsilon:=10^{-8},
\end{equation}
where $(x_n)_{n=0}^{+\infty}$ is the sequence iteratively defined by \eqref{eq:aamr_intro} with $z=x_0$. According to~\Cref{cor:convergent}, for any $\mu\in\,]\gamma(T_{\alpha,\beta}),1[$, the left-hand side of~\eqref{eq:stopping_criterion} is bounded by $\mu^n$ for $n$ big enough. Therefore, an estimate of the maximum number of iterations is given by
\begin{equation}\label{eq:maxiter}
\frac{\log \epsilon}{\log \gamma(T_{\alpha,\beta})}.
\end{equation}
The results are shown in~\Cref{fig:theo_rate}, where the points represent the number of iterations required by AAMR to satisfy~\eqref{eq:stopping_criterion}, and the lines correspond to the estimated upper bounds given by~\eqref{eq:maxiter}. We clearly observe that the algorithm behaves in accordance with the theoretical rates. We emphasize the fact that~\eqref{eq:maxiter} is expected to be a good upper bound on the number of iterations only when this number is \emph{sufficiently} large. We can indeed find a few instances in the plot, especially those which require a small number of iterations, exceeding its estimated upper bound.

\begin{figure}[ht!]
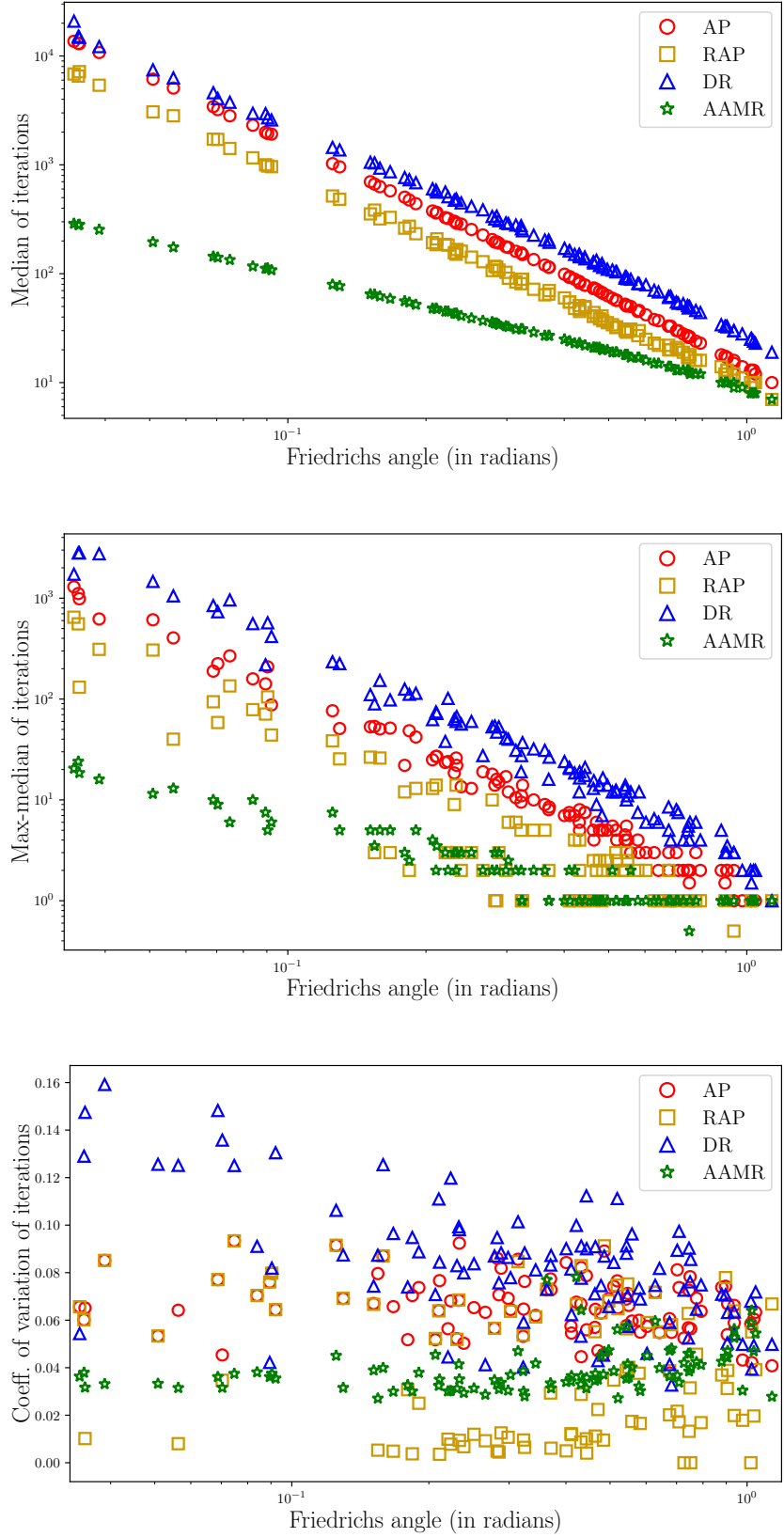

    \centering
    \subfigure{\scalebox{0.56}{\input{Fig7a.pgf}}} 
    \subfigure{\scalebox{0.56}{\input{Fig7b.pgf}}} 
    \subfigure{\scalebox{0.56}{\input{Fig7c.pgf}}} 	 \caption{Median, difference between the maximum and the median, and coefficient of variation of the required number of iterations with respect to the Friedrichs angle of AP, RAP, DR and AAMR for their respective optimal parameters}\label{fig:median_std_optimal}
\end{figure}

In our second experiment we compare the performance of AP, RAP, DR and AAMR, when their parameters are selected to be optimal (see \Cref{tbl:Comparison_rates}).
For $100$~pairs of subspaces, we generated $50$ random starting points with $\|x_0\|=1$. For a fair comparison, we monitored the shadow sequence for all the algorithms.
We also used the stopping criterion~\eqref{eq:stopping_criterion}, with $\epsilon=10^{-8}$. The results of this experiment are summarized in \Cref{fig:median_std_optimal}, where we show in three different graphics the median, the difference between the maximum and the median, and the coefficient of variation of the number of iterations needed to converge for each pair of subspaces. As expected, since the rate of convergence of AAMR is the smallest amongst all the compared methods (see \Cref{tbl:Comparison_rates,fig:comparison_rates}), this algorithm is clearly the fastest, particularly for small angles. Moreover, we can observe that AAMR is one of the most robust methods (together with RAP), which makes the median to be a good representative of the rate of convergence.

\section{Conclusions}\label{sec:conclusions}

We have computed the rate of linear convergence of the \emph{averaged alternating modified reflections (AAMR) method}, which was originally introduced in~\cite{AAMR}, for the case of two subspaces in a Euclidean space. We have additionally found the optimal selection of the parameters defining the scheme that minimizes this rate, in terms of the Friedrichs angle. The rate with optimal parameters coincides with the one of the generalized alternating projections (GAP) method, which is the best among all known rates of projection methods. This coincidence motivated us to study the relationship between AAMR and GAP when they are applied to subspaces. We have discovered that, under some conditions, their associated shadow sequences overlap, which explains the coincidence of the rates. This is not the case for general convex sets.

The developed theoretical results validate the conclusions drawn in the numerical analysis of the convergence rate developed in~\cite[Section~7]{AAMR}. The sharpness of these theoretical results were additionally demonstrated in this work with two computational experiments.

The analysis in this work was done for the case of linear subspaces in a finite-dimensional space. It would be interesting to investigate in future research whether the results can be extended to infinite-dimensional spaces; or even more, to study the rate of convergence of the method when it is applied to two arbitrary convex sets.

\paragraph{Acknowledgements} The authors would like to thank one of the anonymous referees of~\cite{AAMR} for suggesting us this interesting project.

This work was partially supported by  Ministerio de Econom\'ia, Industria y Competitividad (MINECO) of Spain  and  European Regional Development Fund (ERDF), grant MTM2014-59179-C2-1-P. F.J. Arag\'on Artacho was supported by the Ram\'on y Cajal program by MINECO  and  ERDF (RYC-2013-13327) and R. Campoy was supported by MINECO and European Social Fund (BES-2015-073360) under the program ``Ayudas para contratos predoctorales para la formaci\'on de doctores 2015''.

\end{document}

%% file: Fig1a.pgf
\begingroup%
\makeatletter%
\begin{pgfpicture}%
\pgfpathrectangle{\pgfpointorigin}{\pgfqpoint{7.178831in}{5.545448in}}%
\pgfusepath{use as bounding box, clip}%
\begin{pgfscope}%
\pgfsetbuttcap%
\pgfsetmiterjoin%
\definecolor{currentfill}{rgb}{1.000000,1.000000,1.000000}%
\pgfsetfillcolor{currentfill}%
\pgfsetlinewidth{0.000000pt}%
\definecolor{currentstroke}{rgb}{1.000000,1.000000,1.000000}%
\pgfsetstrokecolor{currentstroke}%
\pgfsetdash{}{0pt}%
\pgfpathmoveto{\pgfqpoint{0.000000in}{0.000000in}}%
\pgfpathlineto{\pgfqpoint{7.178831in}{0.000000in}}%
\pgfpathlineto{\pgfqpoint{7.178831in}{5.545448in}}%
\pgfpathlineto{\pgfqpoint{0.000000in}{5.545448in}}%
\pgfpathclose%
\pgfusepath{fill}%
\end{pgfscope}%
\begin{pgfscope}%
\pgfsetbuttcap%
\pgfsetmiterjoin%
\definecolor{currentfill}{rgb}{1.000000,1.000000,1.000000}%
\pgfsetfillcolor{currentfill}%
\pgfsetlinewidth{0.000000pt}%
\definecolor{currentstroke}{rgb}{0.000000,0.000000,0.000000}%
\pgfsetstrokecolor{currentstroke}%
\pgfsetstrokeopacity{0.000000}%
\pgfsetdash{}{0pt}%
\pgfpathmoveto{\pgfqpoint{0.812777in}{0.495102in}}%
\pgfpathlineto{\pgfqpoint{7.012777in}{0.495102in}}%
\pgfpathlineto{\pgfqpoint{7.012777in}{5.025102in}}%
\pgfpathlineto{\pgfqpoint{0.812777in}{5.025102in}}%
\pgfpathclose%
\pgfusepath{fill}%
\end{pgfscope}%
\begin{pgfscope}%
\pgfsetbuttcap%
\pgfsetroundjoin%
\definecolor{currentfill}{rgb}{0.000000,0.000000,0.000000}%
\pgfsetfillcolor{currentfill}%
\pgfsetlinewidth{0.803000pt}%
\definecolor{currentstroke}{rgb}{0.000000,0.000000,0.000000}%
\pgfsetstrokecolor{currentstroke}%
\pgfsetdash{}{0pt}%
\pgfsys@defobject{currentmarker}{\pgfqpoint{0.000000in}{-0.048611in}}{\pgfqpoint{0.000000in}{0.000000in}}{%
\pgfpathmoveto{\pgfqpoint{0.000000in}{0.000000in}}%
\pgfpathlineto{\pgfqpoint{0.000000in}{-0.048611in}}%
\pgfusepath{stroke,fill}%
}%
\begin{pgfscope}%
\pgfsys@transformshift{0.812777in}{0.495102in}%
\pgfsys@useobject{currentmarker}{}%
\end{pgfscope}%
\end{pgfscope}%
\begin{pgfscope}%
\pgftext[x=0.812777in,y=0.397879in,,top]{\rmfamily\fontsize{20.000000}{24.000000}\selectfont  }%
\end{pgfscope}%
\begin{pgfscope}%
\pgfsetbuttcap%
\pgfsetroundjoin%
\definecolor{currentfill}{rgb}{0.000000,0.000000,0.000000}%
\pgfsetfillcolor{currentfill}%
\pgfsetlinewidth{0.803000pt}%
\definecolor{currentstroke}{rgb}{0.000000,0.000000,0.000000}%
\pgfsetstrokecolor{currentstroke}%
\pgfsetdash{}{0pt}%
\pgfsys@defobject{currentmarker}{\pgfqpoint{0.000000in}{-0.048611in}}{\pgfqpoint{0.000000in}{0.000000in}}{%
\pgfpathmoveto{\pgfqpoint{0.000000in}{0.000000in}}%
\pgfpathlineto{\pgfqpoint{0.000000in}{-0.048611in}}%
\pgfusepath{stroke,fill}%
}%
\begin{pgfscope}%
\pgfsys@transformshift{2.052777in}{0.495102in}%
\pgfsys@useobject{currentmarker}{}%
\end{pgfscope}%
\end{pgfscope}%
\begin{pgfscope}%
\pgftext[x=2.052777in,y=0.397879in,,top]{\rmfamily\fontsize{20.000000}{24.000000}\selectfont \(\displaystyle 0.2\)}%
\end{pgfscope}%
\begin{pgfscope}%
\pgfsetbuttcap%
\pgfsetroundjoin%
\definecolor{currentfill}{rgb}{0.000000,0.000000,0.000000}%
\pgfsetfillcolor{currentfill}%
\pgfsetlinewidth{0.803000pt}%
\definecolor{currentstroke}{rgb}{0.000000,0.000000,0.000000}%
\pgfsetstrokecolor{currentstroke}%
\pgfsetdash{}{0pt}%
\pgfsys@defobject{currentmarker}{\pgfqpoint{0.000000in}{-0.048611in}}{\pgfqpoint{0.000000in}{0.000000in}}{%
\pgfpathmoveto{\pgfqpoint{0.000000in}{0.000000in}}%
\pgfpathlineto{\pgfqpoint{0.000000in}{-0.048611in}}%
\pgfusepath{stroke,fill}%
}%
\begin{pgfscope}%
\pgfsys@transformshift{3.292777in}{0.495102in}%
\pgfsys@useobject{currentmarker}{}%
\end{pgfscope}%
\end{pgfscope}%
\begin{pgfscope}%
\pgftext[x=3.292777in,y=0.397879in,,top]{\rmfamily\fontsize{20.000000}{24.000000}\selectfont \(\displaystyle 0.4\)}%
\end{pgfscope}%
\begin{pgfscope}%
\pgfsetbuttcap%
\pgfsetroundjoin%
\definecolor{currentfill}{rgb}{0.000000,0.000000,0.000000}%
\pgfsetfillcolor{currentfill}%
\pgfsetlinewidth{0.803000pt}%
\definecolor{currentstroke}{rgb}{0.000000,0.000000,0.000000}%
\pgfsetstrokecolor{currentstroke}%
\pgfsetdash{}{0pt}%
\pgfsys@defobject{currentmarker}{\pgfqpoint{0.000000in}{-0.048611in}}{\pgfqpoint{0.000000in}{0.000000in}}{%
\pgfpathmoveto{\pgfqpoint{0.000000in}{0.000000in}}%
\pgfpathlineto{\pgfqpoint{0.000000in}{-0.048611in}}%
\pgfusepath{stroke,fill}%
}%
\begin{pgfscope}%
\pgfsys@transformshift{4.532777in}{0.495102in}%
\pgfsys@useobject{currentmarker}{}%
\end{pgfscope}%
\end{pgfscope}%
\begin{pgfscope}%
\pgftext[x=4.532777in,y=0.397879in,,top]{\rmfamily\fontsize{20.000000}{24.000000}\selectfont \(\displaystyle 0.6\)}%
\end{pgfscope}%
\begin{pgfscope}%
\pgfsetbuttcap%
\pgfsetroundjoin%
\definecolor{currentfill}{rgb}{0.000000,0.000000,0.000000}%
\pgfsetfillcolor{currentfill}%
\pgfsetlinewidth{0.803000pt}%
\definecolor{currentstroke}{rgb}{0.000000,0.000000,0.000000}%
\pgfsetstrokecolor{currentstroke}%
\pgfsetdash{}{0pt}%
\pgfsys@defobject{currentmarker}{\pgfqpoint{0.000000in}{-0.048611in}}{\pgfqpoint{0.000000in}{0.000000in}}{%
\pgfpathmoveto{\pgfqpoint{0.000000in}{0.000000in}}%
\pgfpathlineto{\pgfqpoint{0.000000in}{-0.048611in}}%
\pgfusepath{stroke,fill}%
}%
\begin{pgfscope}%
\pgfsys@transformshift{5.772777in}{0.495102in}%
\pgfsys@useobject{currentmarker}{}%
\end{pgfscope}%
\end{pgfscope}%
\begin{pgfscope}%
\pgftext[x=5.772777in,y=0.397879in,,top]{\rmfamily\fontsize{20.000000}{24.000000}\selectfont \(\displaystyle 0.8\)}%
\end{pgfscope}%
\begin{pgfscope}%
\pgfsetbuttcap%
\pgfsetroundjoin%
\definecolor{currentfill}{rgb}{0.000000,0.000000,0.000000}%
\pgfsetfillcolor{currentfill}%
\pgfsetlinewidth{0.803000pt}%
\definecolor{currentstroke}{rgb}{0.000000,0.000000,0.000000}%
\pgfsetstrokecolor{currentstroke}%
\pgfsetdash{}{0pt}%
\pgfsys@defobject{currentmarker}{\pgfqpoint{0.000000in}{-0.048611in}}{\pgfqpoint{0.000000in}{0.000000in}}{%
\pgfpathmoveto{\pgfqpoint{0.000000in}{0.000000in}}%
\pgfpathlineto{\pgfqpoint{0.000000in}{-0.048611in}}%
\pgfusepath{stroke,fill}%
}%
\begin{pgfscope}%
\pgfsys@transformshift{7.012777in}{0.495102in}%
\pgfsys@useobject{currentmarker}{}%
\end{pgfscope}%
\end{pgfscope}%
\begin{pgfscope}%
\pgftext[x=7.012777in,y=0.397879in,,top]{\rmfamily\fontsize{20.000000}{24.000000}\selectfont \(\displaystyle 1\)}%
\end{pgfscope}%
\begin{pgfscope}%
\pgfsetbuttcap%
\pgfsetroundjoin%
\definecolor{currentfill}{rgb}{0.000000,0.000000,0.000000}%
\pgfsetfillcolor{currentfill}%
\pgfsetlinewidth{0.803000pt}%
\definecolor{currentstroke}{rgb}{0.000000,0.000000,0.000000}%
\pgfsetstrokecolor{currentstroke}%
\pgfsetdash{}{0pt}%
\pgfsys@defobject{currentmarker}{\pgfqpoint{0.000000in}{-0.048611in}}{\pgfqpoint{0.000000in}{0.000000in}}{%
\pgfpathmoveto{\pgfqpoint{0.000000in}{0.000000in}}%
\pgfpathlineto{\pgfqpoint{0.000000in}{-0.048611in}}%
\pgfusepath{stroke,fill}%
}%
\begin{pgfscope}%
\pgfsys@transformshift{0.812777in}{0.495102in}%
\pgfsys@useobject{currentmarker}{}%
\end{pgfscope}%
\end{pgfscope}%
\begin{pgfscope}%
\pgftext[x=0.812777in,y=0.397879in,,top]{\rmfamily\fontsize{20.000000}{24.000000}\selectfont \(\displaystyle c(\alpha,\beta)=\widehat{c}_{\beta}\)}%
\end{pgfscope}%
\begin{pgfscope}%
\pgfsetbuttcap%
\pgfsetroundjoin%
\definecolor{currentfill}{rgb}{0.000000,0.000000,0.000000}%
\pgfsetfillcolor{currentfill}%
\pgfsetlinewidth{0.803000pt}%
\definecolor{currentstroke}{rgb}{0.000000,0.000000,0.000000}%
\pgfsetstrokecolor{currentstroke}%
\pgfsetdash{}{0pt}%
\pgfsys@defobject{currentmarker}{\pgfqpoint{-0.048611in}{0.000000in}}{\pgfqpoint{0.000000in}{0.000000in}}{%
\pgfpathmoveto{\pgfqpoint{0.000000in}{0.000000in}}%
\pgfpathlineto{\pgfqpoint{-0.048611in}{0.000000in}}%
\pgfusepath{stroke,fill}%
}%
\begin{pgfscope}%
\pgfsys@transformshift{0.812777in}{0.495102in}%
\pgfsys@useobject{currentmarker}{}%
\end{pgfscope}%
\end{pgfscope}%
\begin{pgfscope}%
\pgftext[x=0.715555in,y=0.398713in,left,base]{\rmfamily\fontsize{20.000000}{24.000000}\selectfont  }%
\end{pgfscope}%
\begin{pgfscope}%
\pgfsetbuttcap%
\pgfsetroundjoin%
\definecolor{currentfill}{rgb}{0.000000,0.000000,0.000000}%
\pgfsetfillcolor{currentfill}%
\pgfsetlinewidth{0.803000pt}%
\definecolor{currentstroke}{rgb}{0.000000,0.000000,0.000000}%
\pgfsetstrokecolor{currentstroke}%
\pgfsetdash{}{0pt}%
\pgfsys@defobject{currentmarker}{\pgfqpoint{-0.048611in}{0.000000in}}{\pgfqpoint{0.000000in}{0.000000in}}{%
\pgfpathmoveto{\pgfqpoint{0.000000in}{0.000000in}}%
\pgfpathlineto{\pgfqpoint{-0.048611in}{0.000000in}}%
\pgfusepath{stroke,fill}%
}%
\begin{pgfscope}%
\pgfsys@transformshift{0.812777in}{1.401102in}%
\pgfsys@useobject{currentmarker}{}%
\end{pgfscope}%
\end{pgfscope}%
\begin{pgfscope}%
\pgftext[x=0.372993in,y=1.304713in,left,base]{\rmfamily\fontsize{20.000000}{24.000000}\selectfont \(\displaystyle 0.2\)}%
\end{pgfscope}%
\begin{pgfscope}%
\pgfsetbuttcap%
\pgfsetroundjoin%
\definecolor{currentfill}{rgb}{0.000000,0.000000,0.000000}%
\pgfsetfillcolor{currentfill}%
\pgfsetlinewidth{0.803000pt}%
\definecolor{currentstroke}{rgb}{0.000000,0.000000,0.000000}%
\pgfsetstrokecolor{currentstroke}%
\pgfsetdash{}{0pt}%
\pgfsys@defobject{currentmarker}{\pgfqpoint{-0.048611in}{0.000000in}}{\pgfqpoint{0.000000in}{0.000000in}}{%
\pgfpathmoveto{\pgfqpoint{0.000000in}{0.000000in}}%
\pgfpathlineto{\pgfqpoint{-0.048611in}{0.000000in}}%
\pgfusepath{stroke,fill}%
}%
\begin{pgfscope}%
\pgfsys@transformshift{0.812777in}{2.307102in}%
\pgfsys@useobject{currentmarker}{}%
\end{pgfscope}%
\end{pgfscope}%
\begin{pgfscope}%
\pgftext[x=0.715555in,y=2.210713in,left,base]{\rmfamily\fontsize{20.000000}{24.000000}\selectfont  }%
\end{pgfscope}%
\begin{pgfscope}%
\pgfsetbuttcap%
\pgfsetroundjoin%
\definecolor{currentfill}{rgb}{0.000000,0.000000,0.000000}%
\pgfsetfillcolor{currentfill}%
\pgfsetlinewidth{0.803000pt}%
\definecolor{currentstroke}{rgb}{0.000000,0.000000,0.000000}%
\pgfsetstrokecolor{currentstroke}%
\pgfsetdash{}{0pt}%
\pgfsys@defobject{currentmarker}{\pgfqpoint{-0.048611in}{0.000000in}}{\pgfqpoint{0.000000in}{0.000000in}}{%
\pgfpathmoveto{\pgfqpoint{0.000000in}{0.000000in}}%
\pgfpathlineto{\pgfqpoint{-0.048611in}{0.000000in}}%
\pgfusepath{stroke,fill}%
}%
\begin{pgfscope}%
\pgfsys@transformshift{0.812777in}{3.213102in}%
\pgfsys@useobject{currentmarker}{}%
\end{pgfscope}%
\end{pgfscope}%
\begin{pgfscope}%
\pgftext[x=0.372993in,y=3.116713in,left,base]{\rmfamily\fontsize{20.000000}{24.000000}\selectfont \(\displaystyle 0.6\)}%
\end{pgfscope}%
\begin{pgfscope}%
\pgfsetbuttcap%
\pgfsetroundjoin%
\definecolor{currentfill}{rgb}{0.000000,0.000000,0.000000}%
\pgfsetfillcolor{currentfill}%
\pgfsetlinewidth{0.803000pt}%
\definecolor{currentstroke}{rgb}{0.000000,0.000000,0.000000}%
\pgfsetstrokecolor{currentstroke}%
\pgfsetdash{}{0pt}%
\pgfsys@defobject{currentmarker}{\pgfqpoint{-0.048611in}{0.000000in}}{\pgfqpoint{0.000000in}{0.000000in}}{%
\pgfpathmoveto{\pgfqpoint{0.000000in}{0.000000in}}%
\pgfpathlineto{\pgfqpoint{-0.048611in}{0.000000in}}%
\pgfusepath{stroke,fill}%
}%
\begin{pgfscope}%
\pgfsys@transformshift{0.812777in}{4.119102in}%
\pgfsys@useobject{currentmarker}{}%
\end{pgfscope}%
\end{pgfscope}%
\begin{pgfscope}%
\pgftext[x=0.372993in,y=4.022713in,left,base]{\rmfamily\fontsize{20.000000}{24.000000}\selectfont \(\displaystyle 0.8\)}%
\end{pgfscope}%
\begin{pgfscope}%
\pgfsetbuttcap%
\pgfsetroundjoin%
\definecolor{currentfill}{rgb}{0.000000,0.000000,0.000000}%
\pgfsetfillcolor{currentfill}%
\pgfsetlinewidth{0.803000pt}%
\definecolor{currentstroke}{rgb}{0.000000,0.000000,0.000000}%
\pgfsetstrokecolor{currentstroke}%
\pgfsetdash{}{0pt}%
\pgfsys@defobject{currentmarker}{\pgfqpoint{-0.048611in}{0.000000in}}{\pgfqpoint{0.000000in}{0.000000in}}{%
\pgfpathmoveto{\pgfqpoint{0.000000in}{0.000000in}}%
\pgfpathlineto{\pgfqpoint{-0.048611in}{0.000000in}}%
\pgfusepath{stroke,fill}%
}%
\begin{pgfscope}%
\pgfsys@transformshift{0.812777in}{5.025102in}%
\pgfsys@useobject{currentmarker}{}%
\end{pgfscope}%
\end{pgfscope}%
\begin{pgfscope}%
\pgftext[x=0.583448in,y=4.928713in,left,base]{\rmfamily\fontsize{20.000000}{24.000000}\selectfont \(\displaystyle 1\)}%
\end{pgfscope}%
\begin{pgfscope}%
\pgfsetbuttcap%
\pgfsetroundjoin%
\definecolor{currentfill}{rgb}{0.000000,0.000000,0.000000}%
\pgfsetfillcolor{currentfill}%
\pgfsetlinewidth{0.803000pt}%
\definecolor{currentstroke}{rgb}{0.000000,0.000000,0.000000}%
\pgfsetstrokecolor{currentstroke}%
\pgfsetdash{}{0pt}%
\pgfsys@defobject{currentmarker}{\pgfqpoint{-0.048611in}{0.000000in}}{\pgfqpoint{0.000000in}{0.000000in}}{%
\pgfpathmoveto{\pgfqpoint{0.000000in}{0.000000in}}%
\pgfpathlineto{\pgfqpoint{-0.048611in}{0.000000in}}%
\pgfusepath{stroke,fill}%
}%
\begin{pgfscope}%
\pgfsys@transformshift{0.812777in}{2.018994in}%
\pgfsys@useobject{currentmarker}{}%
\end{pgfscope}%
\end{pgfscope}%
\begin{pgfscope}%
\pgftext[x=0.224155in,y=1.895081in,left,base]{\rmfamily\fontsize{20.000000}{24.000000}\selectfont \(\displaystyle F_{\alpha,\beta}^1\)}%
\end{pgfscope}%
\begin{pgfscope}%
\pgfsetbuttcap%
\pgfsetroundjoin%
\definecolor{currentfill}{rgb}{0.000000,0.000000,0.000000}%
\pgfsetfillcolor{currentfill}%
\pgfsetlinewidth{0.803000pt}%
\definecolor{currentstroke}{rgb}{0.000000,0.000000,0.000000}%
\pgfsetstrokecolor{currentstroke}%
\pgfsetdash{}{0pt}%
\pgfsys@defobject{currentmarker}{\pgfqpoint{-0.048611in}{0.000000in}}{\pgfqpoint{0.000000in}{0.000000in}}{%
\pgfpathmoveto{\pgfqpoint{0.000000in}{0.000000in}}%
\pgfpathlineto{\pgfqpoint{-0.048611in}{0.000000in}}%
\pgfusepath{stroke,fill}%
}%
\begin{pgfscope}%
\pgfsys@transformshift{0.812777in}{2.714802in}%
\pgfsys@useobject{currentmarker}{}%
\end{pgfscope}%
\end{pgfscope}%
\begin{pgfscope}%
\pgftext[x=0.224155in,y=2.590889in,left,base]{\rmfamily\fontsize{20.000000}{24.000000}\selectfont \(\displaystyle F_{\alpha,\beta}^0\)}%
\end{pgfscope}%
\begin{pgfscope}%
\pgfpathrectangle{\pgfqpoint{0.812777in}{0.495102in}}{\pgfqpoint{6.200000in}{4.530000in}} %
\pgfusepath{clip}%
\pgfsetbuttcap%
\pgfsetroundjoin%
\pgfsetlinewidth{1.505625pt}%
\definecolor{currentstroke}{rgb}{0.501961,0.501961,0.501961}%
\pgfsetstrokecolor{currentstroke}%
\pgfsetdash{{5.550000pt}{2.400000pt}}{0.000000pt}%
\pgfpathmoveto{\pgfqpoint{0.812777in}{0.495102in}}%
\pgfpathlineto{\pgfqpoint{0.812777in}{5.038990in}}%
\pgfusepath{stroke}%
\end{pgfscope}%
\begin{pgfscope}%
\pgfpathrectangle{\pgfqpoint{0.812777in}{0.495102in}}{\pgfqpoint{6.200000in}{4.530000in}} %
\pgfusepath{clip}%
\pgfsetbuttcap%
\pgfsetroundjoin%
\pgfsetlinewidth{1.003750pt}%
\definecolor{currentstroke}{rgb}{0.000000,0.392157,0.000000}%
\pgfsetstrokecolor{currentstroke}%
\pgfsetdash{{3.700000pt}{1.600000pt}}{0.000000pt}%
\pgfpathmoveto{\pgfqpoint{0.812777in}{2.018994in}}%
\pgfpathlineto{\pgfqpoint{7.012777in}{2.018994in}}%
\pgfusepath{stroke}%
\end{pgfscope}%
\begin{pgfscope}%
\pgfpathrectangle{\pgfqpoint{0.812777in}{0.495102in}}{\pgfqpoint{6.200000in}{4.530000in}} %
\pgfusepath{clip}%
\pgfsetbuttcap%
\pgfsetroundjoin%
\pgfsetlinewidth{1.505625pt}%
\definecolor{currentstroke}{rgb}{0.501961,0.501961,0.501961}%
\pgfsetstrokecolor{currentstroke}%
\pgfsetdash{{5.550000pt}{2.400000pt}}{0.000000pt}%
\pgfpathmoveto{\pgfqpoint{0.812777in}{0.495102in}}%
\pgfpathlineto{\pgfqpoint{0.812777in}{2.018994in}}%
\pgfusepath{stroke}%
\end{pgfscope}%
\begin{pgfscope}%
\pgfpathrectangle{\pgfqpoint{0.812777in}{0.495102in}}{\pgfqpoint{6.200000in}{4.530000in}} %
\pgfusepath{clip}%
\pgfsetbuttcap%
\pgfsetroundjoin%
\pgfsetlinewidth{3.011250pt}%
\definecolor{currentstroke}{rgb}{0.000000,0.000000,1.000000}%
\pgfsetstrokecolor{currentstroke}%
\pgfsetdash{{19.200000pt}{4.800000pt}{3.000000pt}{4.800000pt}}{0.000000pt}%
\pgfpathmoveto{\pgfqpoint{0.812777in}{2.714802in}}%
\pgfpathlineto{\pgfqpoint{1.097977in}{2.771735in}}%
\pgfpathlineto{\pgfqpoint{1.376977in}{2.830604in}}%
\pgfpathlineto{\pgfqpoint{1.649777in}{2.891361in}}%
\pgfpathlineto{\pgfqpoint{1.916377in}{2.953954in}}%
\pgfpathlineto{\pgfqpoint{2.176777in}{3.018314in}}%
\pgfpathlineto{\pgfqpoint{2.430977in}{3.084367in}}%
\pgfpathlineto{\pgfqpoint{2.678977in}{3.152024in}}%
\pgfpathlineto{\pgfqpoint{2.920777in}{3.221188in}}%
\pgfpathlineto{\pgfqpoint{3.156377in}{3.291751in}}%
\pgfpathlineto{\pgfqpoint{3.391977in}{3.365580in}}%
\pgfpathlineto{\pgfqpoint{3.621377in}{3.440738in}}%
\pgfpathlineto{\pgfqpoint{3.844577in}{3.517090in}}%
\pgfpathlineto{\pgfqpoint{4.067777in}{3.596750in}}%
\pgfpathlineto{\pgfqpoint{4.284777in}{3.677491in}}%
\pgfpathlineto{\pgfqpoint{4.495577in}{3.759154in}}%
\pgfpathlineto{\pgfqpoint{4.706377in}{3.844115in}}%
\pgfpathlineto{\pgfqpoint{4.910977in}{3.929844in}}%
\pgfpathlineto{\pgfqpoint{5.115577in}{4.018904in}}%
\pgfpathlineto{\pgfqpoint{5.313977in}{4.108555in}}%
\pgfpathlineto{\pgfqpoint{5.512377in}{4.201555in}}%
\pgfpathlineto{\pgfqpoint{5.710777in}{4.298013in}}%
\pgfpathlineto{\pgfqpoint{5.902977in}{4.394861in}}%
\pgfpathlineto{\pgfqpoint{6.095177in}{4.495165in}}%
\pgfpathlineto{\pgfqpoint{6.287377in}{4.599033in}}%
\pgfpathlineto{\pgfqpoint{6.473377in}{4.703048in}}%
\pgfpathlineto{\pgfqpoint{6.659377in}{4.810606in}}%
\pgfpathlineto{\pgfqpoint{6.845377in}{4.921811in}}%
\pgfpathlineto{\pgfqpoint{7.012777in}{5.025102in}}%
\pgfpathlineto{\pgfqpoint{7.012777in}{5.025102in}}%
\pgfusepath{stroke}%
\end{pgfscope}%
\begin{pgfscope}%
\pgfpathrectangle{\pgfqpoint{0.812777in}{0.495102in}}{\pgfqpoint{6.200000in}{4.530000in}} %
\pgfusepath{clip}%
\pgfsetbuttcap%
\pgfsetroundjoin%
\pgfsetlinewidth{2.509375pt}%
\definecolor{currentstroke}{rgb}{1.000000,0.000000,0.000000}%
\pgfsetstrokecolor{currentstroke}%
\pgfsetdash{{9.250000pt}{4.000000pt}}{0.000000pt}%
\pgfpathmoveto{\pgfqpoint{0.812777in}{2.714802in}}%
\pgfpathlineto{\pgfqpoint{1.104177in}{2.659837in}}%
\pgfpathlineto{\pgfqpoint{1.401777in}{2.606867in}}%
\pgfpathlineto{\pgfqpoint{1.705577in}{2.555907in}}%
\pgfpathlineto{\pgfqpoint{2.021777in}{2.506014in}}%
\pgfpathlineto{\pgfqpoint{2.344177in}{2.458249in}}%
\pgfpathlineto{\pgfqpoint{2.678977in}{2.411761in}}%
\pgfpathlineto{\pgfqpoint{3.026177in}{2.366686in}}%
\pgfpathlineto{\pgfqpoint{3.385777in}{2.323140in}}%
\pgfpathlineto{\pgfqpoint{3.757777in}{2.281217in}}%
\pgfpathlineto{\pgfqpoint{4.142177in}{2.240993in}}%
\pgfpathlineto{\pgfqpoint{4.545177in}{2.201944in}}%
\pgfpathlineto{\pgfqpoint{4.966777in}{2.164241in}}%
\pgfpathlineto{\pgfqpoint{5.406977in}{2.128025in}}%
\pgfpathlineto{\pgfqpoint{5.865777in}{2.093405in}}%
\pgfpathlineto{\pgfqpoint{6.349377in}{2.060053in}}%
\pgfpathlineto{\pgfqpoint{6.857777in}{2.028134in}}%
\pgfpathlineto{\pgfqpoint{7.012777in}{2.018994in}}%
\pgfpathlineto{\pgfqpoint{7.012777in}{2.018994in}}%
\pgfusepath{stroke}%
\end{pgfscope}%
\begin{pgfscope}%
\pgfpathrectangle{\pgfqpoint{0.812777in}{0.495102in}}{\pgfqpoint{6.200000in}{4.530000in}} %
\pgfusepath{clip}%
\pgfsetbuttcap%
\pgfsetroundjoin%
\definecolor{currentfill}{rgb}{0.721569,0.525490,0.043137}%
\pgfsetfillcolor{currentfill}%
\pgfsetlinewidth{1.003750pt}%
\definecolor{currentstroke}{rgb}{0.721569,0.525490,0.043137}%
\pgfsetstrokecolor{currentstroke}%
\pgfsetdash{}{0pt}%
\pgfsys@defobject{currentmarker}{\pgfqpoint{-0.027778in}{-0.027778in}}{\pgfqpoint{0.027778in}{0.027778in}}{%
\pgfpathmoveto{\pgfqpoint{0.000000in}{-0.027778in}}%
\pgfpathcurveto{\pgfqpoint{0.007367in}{-0.027778in}}{\pgfqpoint{0.014433in}{-0.024851in}}{\pgfqpoint{0.019642in}{-0.019642in}}%
\pgfpathcurveto{\pgfqpoint{0.024851in}{-0.014433in}}{\pgfqpoint{0.027778in}{-0.007367in}}{\pgfqpoint{0.027778in}{0.000000in}}%
\pgfpathcurveto{\pgfqpoint{0.027778in}{0.007367in}}{\pgfqpoint{0.024851in}{0.014433in}}{\pgfqpoint{0.019642in}{0.019642in}}%
\pgfpathcurveto{\pgfqpoint{0.014433in}{0.024851in}}{\pgfqpoint{0.007367in}{0.027778in}}{\pgfqpoint{0.000000in}{0.027778in}}%
\pgfpathcurveto{\pgfqpoint{-0.007367in}{0.027778in}}{\pgfqpoint{-0.014433in}{0.024851in}}{\pgfqpoint{-0.019642in}{0.019642in}}%
\pgfpathcurveto{\pgfqpoint{-0.024851in}{0.014433in}}{\pgfqpoint{-0.027778in}{0.007367in}}{\pgfqpoint{-0.027778in}{0.000000in}}%
\pgfpathcurveto{\pgfqpoint{-0.027778in}{-0.007367in}}{\pgfqpoint{-0.024851in}{-0.014433in}}{\pgfqpoint{-0.019642in}{-0.019642in}}%
\pgfpathcurveto{\pgfqpoint{-0.014433in}{-0.024851in}}{\pgfqpoint{-0.007367in}{-0.027778in}}{\pgfqpoint{0.000000in}{-0.027778in}}%
\pgfpathclose%
\pgfusepath{stroke,fill}%
}%
\begin{pgfscope}%
\pgfsys@transformshift{0.818977in}{2.714802in}%
\pgfsys@useobject{currentmarker}{}%
\end{pgfscope}%
\end{pgfscope}%
\begin{pgfscope}%
\pgfpathrectangle{\pgfqpoint{0.812777in}{0.495102in}}{\pgfqpoint{6.200000in}{4.530000in}} %
\pgfusepath{clip}%
\pgfsetbuttcap%
\pgfsetroundjoin%
\definecolor{currentfill}{rgb}{0.000000,0.000000,1.000000}%
\pgfsetfillcolor{currentfill}%
\pgfsetlinewidth{1.003750pt}%
\definecolor{currentstroke}{rgb}{0.000000,0.000000,1.000000}%
\pgfsetstrokecolor{currentstroke}%
\pgfsetdash{}{0pt}%
\pgfsys@defobject{currentmarker}{\pgfqpoint{-0.027778in}{-0.027778in}}{\pgfqpoint{0.027778in}{0.027778in}}{%
\pgfpathmoveto{\pgfqpoint{0.000000in}{-0.027778in}}%
\pgfpathcurveto{\pgfqpoint{0.007367in}{-0.027778in}}{\pgfqpoint{0.014433in}{-0.024851in}}{\pgfqpoint{0.019642in}{-0.019642in}}%
\pgfpathcurveto{\pgfqpoint{0.024851in}{-0.014433in}}{\pgfqpoint{0.027778in}{-0.007367in}}{\pgfqpoint{0.027778in}{0.000000in}}%
\pgfpathcurveto{\pgfqpoint{0.027778in}{0.007367in}}{\pgfqpoint{0.024851in}{0.014433in}}{\pgfqpoint{0.019642in}{0.019642in}}%
\pgfpathcurveto{\pgfqpoint{0.014433in}{0.024851in}}{\pgfqpoint{0.007367in}{0.027778in}}{\pgfqpoint{0.000000in}{0.027778in}}%
\pgfpathcurveto{\pgfqpoint{-0.007367in}{0.027778in}}{\pgfqpoint{-0.014433in}{0.024851in}}{\pgfqpoint{-0.019642in}{0.019642in}}%
\pgfpathcurveto{\pgfqpoint{-0.024851in}{0.014433in}}{\pgfqpoint{-0.027778in}{0.007367in}}{\pgfqpoint{-0.027778in}{0.000000in}}%
\pgfpathcurveto{\pgfqpoint{-0.027778in}{-0.007367in}}{\pgfqpoint{-0.024851in}{-0.014433in}}{\pgfqpoint{-0.019642in}{-0.019642in}}%
\pgfpathcurveto{\pgfqpoint{-0.014433in}{-0.024851in}}{\pgfqpoint{-0.007367in}{-0.027778in}}{\pgfqpoint{0.000000in}{-0.027778in}}%
\pgfpathclose%
\pgfusepath{stroke,fill}%
}%
\begin{pgfscope}%
\pgfsys@transformshift{7.006577in}{5.025102in}%
\pgfsys@useobject{currentmarker}{}%
\end{pgfscope}%
\end{pgfscope}%
\begin{pgfscope}%
\pgfpathrectangle{\pgfqpoint{0.812777in}{0.495102in}}{\pgfqpoint{6.200000in}{4.530000in}} %
\pgfusepath{clip}%
\pgfsetbuttcap%
\pgfsetroundjoin%
\definecolor{currentfill}{rgb}{1.000000,0.000000,0.000000}%
\pgfsetfillcolor{currentfill}%
\pgfsetlinewidth{1.003750pt}%
\definecolor{currentstroke}{rgb}{1.000000,0.000000,0.000000}%
\pgfsetstrokecolor{currentstroke}%
\pgfsetdash{}{0pt}%
\pgfsys@defobject{currentmarker}{\pgfqpoint{-0.027778in}{-0.027778in}}{\pgfqpoint{0.027778in}{0.027778in}}{%
\pgfpathmoveto{\pgfqpoint{0.000000in}{-0.027778in}}%
\pgfpathcurveto{\pgfqpoint{0.007367in}{-0.027778in}}{\pgfqpoint{0.014433in}{-0.024851in}}{\pgfqpoint{0.019642in}{-0.019642in}}%
\pgfpathcurveto{\pgfqpoint{0.024851in}{-0.014433in}}{\pgfqpoint{0.027778in}{-0.007367in}}{\pgfqpoint{0.027778in}{0.000000in}}%
\pgfpathcurveto{\pgfqpoint{0.027778in}{0.007367in}}{\pgfqpoint{0.024851in}{0.014433in}}{\pgfqpoint{0.019642in}{0.019642in}}%
\pgfpathcurveto{\pgfqpoint{0.014433in}{0.024851in}}{\pgfqpoint{0.007367in}{0.027778in}}{\pgfqpoint{0.000000in}{0.027778in}}%
\pgfpathcurveto{\pgfqpoint{-0.007367in}{0.027778in}}{\pgfqpoint{-0.014433in}{0.024851in}}{\pgfqpoint{-0.019642in}{0.019642in}}%
\pgfpathcurveto{\pgfqpoint{-0.024851in}{0.014433in}}{\pgfqpoint{-0.027778in}{0.007367in}}{\pgfqpoint{-0.027778in}{0.000000in}}%
\pgfpathcurveto{\pgfqpoint{-0.027778in}{-0.007367in}}{\pgfqpoint{-0.024851in}{-0.014433in}}{\pgfqpoint{-0.019642in}{-0.019642in}}%
\pgfpathcurveto{\pgfqpoint{-0.014433in}{-0.024851in}}{\pgfqpoint{-0.007367in}{-0.027778in}}{\pgfqpoint{0.000000in}{-0.027778in}}%
\pgfpathclose%
\pgfusepath{stroke,fill}%
}%
\begin{pgfscope}%
\pgfsys@transformshift{7.006577in}{2.018994in}%
\pgfsys@useobject{currentmarker}{}%
\end{pgfscope}%
\end{pgfscope}%
\begin{pgfscope}%
\pgfpathrectangle{\pgfqpoint{0.812777in}{0.495102in}}{\pgfqpoint{6.200000in}{4.530000in}} %
\pgfusepath{clip}%
\pgfsetbuttcap%
\pgfsetroundjoin%
\definecolor{currentfill}{rgb}{0.721569,0.525490,0.043137}%
\pgfsetfillcolor{currentfill}%
\pgfsetlinewidth{1.003750pt}%
\definecolor{currentstroke}{rgb}{0.721569,0.525490,0.043137}%
\pgfsetstrokecolor{currentstroke}%
\pgfsetdash{}{0pt}%
\pgfsys@defobject{currentmarker}{\pgfqpoint{-0.027778in}{-0.027778in}}{\pgfqpoint{0.027778in}{0.027778in}}{%
\pgfpathmoveto{\pgfqpoint{0.000000in}{-0.027778in}}%
\pgfpathcurveto{\pgfqpoint{0.007367in}{-0.027778in}}{\pgfqpoint{0.014433in}{-0.024851in}}{\pgfqpoint{0.019642in}{-0.019642in}}%
\pgfpathcurveto{\pgfqpoint{0.024851in}{-0.014433in}}{\pgfqpoint{0.027778in}{-0.007367in}}{\pgfqpoint{0.027778in}{0.000000in}}%
\pgfpathcurveto{\pgfqpoint{0.027778in}{0.007367in}}{\pgfqpoint{0.024851in}{0.014433in}}{\pgfqpoint{0.019642in}{0.019642in}}%
\pgfpathcurveto{\pgfqpoint{0.014433in}{0.024851in}}{\pgfqpoint{0.007367in}{0.027778in}}{\pgfqpoint{0.000000in}{0.027778in}}%
\pgfpathcurveto{\pgfqpoint{-0.007367in}{0.027778in}}{\pgfqpoint{-0.014433in}{0.024851in}}{\pgfqpoint{-0.019642in}{0.019642in}}%
\pgfpathcurveto{\pgfqpoint{-0.024851in}{0.014433in}}{\pgfqpoint{-0.027778in}{0.007367in}}{\pgfqpoint{-0.027778in}{0.000000in}}%
\pgfpathcurveto{\pgfqpoint{-0.027778in}{-0.007367in}}{\pgfqpoint{-0.024851in}{-0.014433in}}{\pgfqpoint{-0.019642in}{-0.019642in}}%
\pgfpathcurveto{\pgfqpoint{-0.014433in}{-0.024851in}}{\pgfqpoint{-0.007367in}{-0.027778in}}{\pgfqpoint{0.000000in}{-0.027778in}}%
\pgfpathclose%
\pgfusepath{stroke,fill}%
}%
\begin{pgfscope}%
\pgfsys@transformshift{0.812777in}{2.018994in}%
\pgfsys@useobject{currentmarker}{}%
\end{pgfscope}%
\end{pgfscope}%
\begin{pgfscope}%
\pgfsetrectcap%
\pgfsetmiterjoin%
\pgfsetlinewidth{0.803000pt}%
\definecolor{currentstroke}{rgb}{0.000000,0.000000,0.000000}%
\pgfsetstrokecolor{currentstroke}%
\pgfsetdash{}{0pt}%
\pgfpathmoveto{\pgfqpoint{0.812777in}{0.495102in}}%
\pgfpathlineto{\pgfqpoint{0.812777in}{5.025102in}}%
\pgfusepath{stroke}%
\end{pgfscope}%
\begin{pgfscope}%
\pgfsetrectcap%
\pgfsetmiterjoin%
\pgfsetlinewidth{0.803000pt}%
\definecolor{currentstroke}{rgb}{0.000000,0.000000,0.000000}%
\pgfsetstrokecolor{currentstroke}%
\pgfsetdash{}{0pt}%
\pgfpathmoveto{\pgfqpoint{7.012777in}{0.495102in}}%
\pgfpathlineto{\pgfqpoint{7.012777in}{5.025102in}}%
\pgfusepath{stroke}%
\end{pgfscope}%
\begin{pgfscope}%
\pgfsetrectcap%
\pgfsetmiterjoin%
\pgfsetlinewidth{0.803000pt}%
\definecolor{currentstroke}{rgb}{0.000000,0.000000,0.000000}%
\pgfsetstrokecolor{currentstroke}%
\pgfsetdash{}{0pt}%
\pgfpathmoveto{\pgfqpoint{0.812777in}{0.495102in}}%
\pgfpathlineto{\pgfqpoint{7.012777in}{0.495102in}}%
\pgfusepath{stroke}%
\end{pgfscope}%
\begin{pgfscope}%
\pgfsetrectcap%
\pgfsetmiterjoin%
\pgfsetlinewidth{0.803000pt}%
\definecolor{currentstroke}{rgb}{0.000000,0.000000,0.000000}%
\pgfsetstrokecolor{currentstroke}%
\pgfsetdash{}{0pt}%
\pgfpathmoveto{\pgfqpoint{0.812777in}{5.025102in}}%
\pgfpathlineto{\pgfqpoint{7.012777in}{5.025102in}}%
\pgfusepath{stroke}%
\end{pgfscope}%
\begin{pgfscope}%
\definecolor{textcolor}{rgb}{0.000000,0.392157,0.000000}%
\pgfsetstrokecolor{textcolor}%
\pgfsetfillcolor{textcolor}%
\pgftext[x=3.044777in,y=1.747194in,left,base]{\color{textcolor}\rmfamily\fontsize{20.000000}{24.000000}\selectfont \(\displaystyle y=f_{\alpha,\beta,1}(1)\)}%
\end{pgfscope}%
\begin{pgfscope}%
\pgfsetbuttcap%
\pgfsetmiterjoin%
\definecolor{currentfill}{rgb}{1.000000,1.000000,1.000000}%
\pgfsetfillcolor{currentfill}%
\pgfsetlinewidth{1.003750pt}%
\definecolor{currentstroke}{rgb}{0.000000,0.000000,0.000000}%
\pgfsetstrokecolor{currentstroke}%
\pgfsetdash{}{0pt}%
\pgfpathmoveto{\pgfqpoint{3.081310in}{5.176559in}}%
\pgfpathlineto{\pgfqpoint{4.744244in}{5.176559in}}%
\pgfpathlineto{\pgfqpoint{4.744244in}{5.501003in}}%
\pgfpathlineto{\pgfqpoint{3.081310in}{5.501003in}}%
\pgfpathclose%
\pgfusepath{stroke,fill}%
\end{pgfscope}%
\begin{pgfscope}%
\pgftext[x=3.912777in,y=5.278782in,,base]{\rmfamily\fontsize{16.000000}{19.200000}\selectfont \(\displaystyle \alpha=\) 0.5, \(\displaystyle \beta=\) 0.3}%
\end{pgfscope}%
\begin{pgfscope}%
\pgfsetbuttcap%
\pgfsetmiterjoin%
\definecolor{currentfill}{rgb}{1.000000,1.000000,1.000000}%
\pgfsetfillcolor{currentfill}%
\pgfsetfillopacity{0.800000}%
\pgfsetlinewidth{1.003750pt}%
\definecolor{currentstroke}{rgb}{0.800000,0.800000,0.800000}%
\pgfsetstrokecolor{currentstroke}%
\pgfsetstrokeopacity{0.800000}%
\pgfsetdash{}{0pt}%
\pgfpathmoveto{\pgfqpoint{0.992333in}{3.937902in}}%
\pgfpathlineto{\pgfqpoint{2.617894in}{3.937902in}}%
\pgfpathquadraticcurveto{\pgfqpoint{2.673450in}{3.937902in}}{\pgfqpoint{2.673450in}{3.993457in}}%
\pgfpathlineto{\pgfqpoint{2.673450in}{4.839216in}}%
\pgfpathquadraticcurveto{\pgfqpoint{2.673450in}{4.894771in}}{\pgfqpoint{2.617894in}{4.894771in}}%
\pgfpathlineto{\pgfqpoint{0.992333in}{4.894771in}}%
\pgfpathquadraticcurveto{\pgfqpoint{0.936777in}{4.894771in}}{\pgfqpoint{0.936777in}{4.839216in}}%
\pgfpathlineto{\pgfqpoint{0.936777in}{3.993457in}}%
\pgfpathquadraticcurveto{\pgfqpoint{0.936777in}{3.937902in}}{\pgfqpoint{0.992333in}{3.937902in}}%
\pgfpathclose%
\pgfusepath{stroke,fill}%
\end{pgfscope}%
\begin{pgfscope}%
\pgfsetbuttcap%
\pgfsetroundjoin%
\pgfsetlinewidth{3.011250pt}%
\definecolor{currentstroke}{rgb}{0.000000,0.000000,1.000000}%
\pgfsetstrokecolor{currentstroke}%
\pgfsetdash{{19.200000pt}{4.800000pt}{3.000000pt}{4.800000pt}}{0.000000pt}%
\pgfpathmoveto{\pgfqpoint{1.047888in}{4.665100in}}%
\pgfpathlineto{\pgfqpoint{1.414555in}{4.665100in}}%
\pgfusepath{stroke}%
\end{pgfscope}%
\begin{pgfscope}%
\pgftext[x=1.636777in,y=4.567878in,left,base]{\rmfamily\fontsize{20.000000}{24.000000}\selectfont \(\displaystyle f_{\alpha,\beta,2}(c)\)}%
\end{pgfscope}%
\begin{pgfscope}%
\pgfsetbuttcap%
\pgfsetroundjoin%
\pgfsetlinewidth{2.509375pt}%
\definecolor{currentstroke}{rgb}{1.000000,0.000000,0.000000}%
\pgfsetstrokecolor{currentstroke}%
\pgfsetdash{{9.250000pt}{4.000000pt}}{0.000000pt}%
\pgfpathmoveto{\pgfqpoint{1.047888in}{4.228332in}}%
\pgfpathlineto{\pgfqpoint{1.414555in}{4.228332in}}%
\pgfusepath{stroke}%
\end{pgfscope}%
\begin{pgfscope}%
\pgftext[x=1.636777in,y=4.131110in,left,base]{\rmfamily\fontsize{20.000000}{24.000000}\selectfont \(\displaystyle f_{\alpha,\beta,1}(c)\)}%
\end{pgfscope}%
\end{pgfpicture}%
\makeatother%
\endgroup%

%% file: Fig1b.pgf
\begingroup%
\makeatletter%
\begin{pgfpicture}%
\pgfpathrectangle{\pgfpointorigin}{\pgfqpoint{7.054676in}{5.535124in}}%
\pgfusepath{use as bounding box, clip}%
\begin{pgfscope}%
\pgfsetbuttcap%
\pgfsetmiterjoin%
\definecolor{currentfill}{rgb}{1.000000,1.000000,1.000000}%
\pgfsetfillcolor{currentfill}%
\pgfsetlinewidth{0.000000pt}%
\definecolor{currentstroke}{rgb}{1.000000,1.000000,1.000000}%
\pgfsetstrokecolor{currentstroke}%
\pgfsetdash{}{0pt}%
\pgfpathmoveto{\pgfqpoint{0.000000in}{0.000000in}}%
\pgfpathlineto{\pgfqpoint{7.054676in}{0.000000in}}%
\pgfpathlineto{\pgfqpoint{7.054676in}{5.535124in}}%
\pgfpathlineto{\pgfqpoint{0.000000in}{5.535124in}}%
\pgfpathclose%
\pgfusepath{fill}%
\end{pgfscope}%
\begin{pgfscope}%
\pgfsetbuttcap%
\pgfsetmiterjoin%
\definecolor{currentfill}{rgb}{1.000000,1.000000,1.000000}%
\pgfsetfillcolor{currentfill}%
\pgfsetlinewidth{0.000000pt}%
\definecolor{currentstroke}{rgb}{0.000000,0.000000,0.000000}%
\pgfsetstrokecolor{currentstroke}%
\pgfsetstrokeopacity{0.000000}%
\pgfsetdash{}{0pt}%
\pgfpathmoveto{\pgfqpoint{0.688622in}{0.484778in}}%
\pgfpathlineto{\pgfqpoint{6.888622in}{0.484778in}}%
\pgfpathlineto{\pgfqpoint{6.888622in}{5.014778in}}%
\pgfpathlineto{\pgfqpoint{0.688622in}{5.014778in}}%
\pgfpathclose%
\pgfusepath{fill}%
\end{pgfscope}%
\begin{pgfscope}%
\pgfsetbuttcap%
\pgfsetroundjoin%
\definecolor{currentfill}{rgb}{0.000000,0.000000,0.000000}%
\pgfsetfillcolor{currentfill}%
\pgfsetlinewidth{0.803000pt}%
\definecolor{currentstroke}{rgb}{0.000000,0.000000,0.000000}%
\pgfsetstrokecolor{currentstroke}%
\pgfsetdash{}{0pt}%
\pgfsys@defobject{currentmarker}{\pgfqpoint{0.000000in}{-0.048611in}}{\pgfqpoint{0.000000in}{0.000000in}}{%
\pgfpathmoveto{\pgfqpoint{0.000000in}{0.000000in}}%
\pgfpathlineto{\pgfqpoint{0.000000in}{-0.048611in}}%
\pgfusepath{stroke,fill}%
}%
\begin{pgfscope}%
\pgfsys@transformshift{0.688622in}{0.484778in}%
\pgfsys@useobject{currentmarker}{}%
\end{pgfscope}%
\end{pgfscope}%
\begin{pgfscope}%
\pgftext[x=0.688622in,y=0.387555in,,top]{\rmfamily\fontsize{20.000000}{24.000000}\selectfont  }%
\end{pgfscope}%
\begin{pgfscope}%
\pgfsetbuttcap%
\pgfsetroundjoin%
\definecolor{currentfill}{rgb}{0.000000,0.000000,0.000000}%
\pgfsetfillcolor{currentfill}%
\pgfsetlinewidth{0.803000pt}%
\definecolor{currentstroke}{rgb}{0.000000,0.000000,0.000000}%
\pgfsetstrokecolor{currentstroke}%
\pgfsetdash{}{0pt}%
\pgfsys@defobject{currentmarker}{\pgfqpoint{0.000000in}{-0.048611in}}{\pgfqpoint{0.000000in}{0.000000in}}{%
\pgfpathmoveto{\pgfqpoint{0.000000in}{0.000000in}}%
\pgfpathlineto{\pgfqpoint{0.000000in}{-0.048611in}}%
\pgfusepath{stroke,fill}%
}%
\begin{pgfscope}%
\pgfsys@transformshift{1.928622in}{0.484778in}%
\pgfsys@useobject{currentmarker}{}%
\end{pgfscope}%
\end{pgfscope}%
\begin{pgfscope}%
\pgftext[x=1.928622in,y=0.387555in,,top]{\rmfamily\fontsize{20.000000}{24.000000}\selectfont \(\displaystyle 0.2\)}%
\end{pgfscope}%
\begin{pgfscope}%
\pgfsetbuttcap%
\pgfsetroundjoin%
\definecolor{currentfill}{rgb}{0.000000,0.000000,0.000000}%
\pgfsetfillcolor{currentfill}%
\pgfsetlinewidth{0.803000pt}%
\definecolor{currentstroke}{rgb}{0.000000,0.000000,0.000000}%
\pgfsetstrokecolor{currentstroke}%
\pgfsetdash{}{0pt}%
\pgfsys@defobject{currentmarker}{\pgfqpoint{0.000000in}{-0.048611in}}{\pgfqpoint{0.000000in}{0.000000in}}{%
\pgfpathmoveto{\pgfqpoint{0.000000in}{0.000000in}}%
\pgfpathlineto{\pgfqpoint{0.000000in}{-0.048611in}}%
\pgfusepath{stroke,fill}%
}%
\begin{pgfscope}%
\pgfsys@transformshift{3.168622in}{0.484778in}%
\pgfsys@useobject{currentmarker}{}%
\end{pgfscope}%
\end{pgfscope}%
\begin{pgfscope}%
\pgftext[x=3.168622in,y=0.387555in,,top]{\rmfamily\fontsize{20.000000}{24.000000}\selectfont \(\displaystyle 0.4\)}%
\end{pgfscope}%
\begin{pgfscope}%
\pgfsetbuttcap%
\pgfsetroundjoin%
\definecolor{currentfill}{rgb}{0.000000,0.000000,0.000000}%
\pgfsetfillcolor{currentfill}%
\pgfsetlinewidth{0.803000pt}%
\definecolor{currentstroke}{rgb}{0.000000,0.000000,0.000000}%
\pgfsetstrokecolor{currentstroke}%
\pgfsetdash{}{0pt}%
\pgfsys@defobject{currentmarker}{\pgfqpoint{0.000000in}{-0.048611in}}{\pgfqpoint{0.000000in}{0.000000in}}{%
\pgfpathmoveto{\pgfqpoint{0.000000in}{0.000000in}}%
\pgfpathlineto{\pgfqpoint{0.000000in}{-0.048611in}}%
\pgfusepath{stroke,fill}%
}%
\begin{pgfscope}%
\pgfsys@transformshift{4.408622in}{0.484778in}%
\pgfsys@useobject{currentmarker}{}%
\end{pgfscope}%
\end{pgfscope}%
\begin{pgfscope}%
\pgftext[x=4.408622in,y=0.387555in,,top]{\rmfamily\fontsize{20.000000}{24.000000}\selectfont \(\displaystyle 0.6\)}%
\end{pgfscope}%
\begin{pgfscope}%
\pgfsetbuttcap%
\pgfsetroundjoin%
\definecolor{currentfill}{rgb}{0.000000,0.000000,0.000000}%
\pgfsetfillcolor{currentfill}%
\pgfsetlinewidth{0.803000pt}%
\definecolor{currentstroke}{rgb}{0.000000,0.000000,0.000000}%
\pgfsetstrokecolor{currentstroke}%
\pgfsetdash{}{0pt}%
\pgfsys@defobject{currentmarker}{\pgfqpoint{0.000000in}{-0.048611in}}{\pgfqpoint{0.000000in}{0.000000in}}{%
\pgfpathmoveto{\pgfqpoint{0.000000in}{0.000000in}}%
\pgfpathlineto{\pgfqpoint{0.000000in}{-0.048611in}}%
\pgfusepath{stroke,fill}%
}%
\begin{pgfscope}%
\pgfsys@transformshift{5.648622in}{0.484778in}%
\pgfsys@useobject{currentmarker}{}%
\end{pgfscope}%
\end{pgfscope}%
\begin{pgfscope}%
\pgftext[x=5.648622in,y=0.387555in,,top]{\rmfamily\fontsize{20.000000}{24.000000}\selectfont \(\displaystyle 0.8\)}%
\end{pgfscope}%
\begin{pgfscope}%
\pgfsetbuttcap%
\pgfsetroundjoin%
\definecolor{currentfill}{rgb}{0.000000,0.000000,0.000000}%
\pgfsetfillcolor{currentfill}%
\pgfsetlinewidth{0.803000pt}%
\definecolor{currentstroke}{rgb}{0.000000,0.000000,0.000000}%
\pgfsetstrokecolor{currentstroke}%
\pgfsetdash{}{0pt}%
\pgfsys@defobject{currentmarker}{\pgfqpoint{0.000000in}{-0.048611in}}{\pgfqpoint{0.000000in}{0.000000in}}{%
\pgfpathmoveto{\pgfqpoint{0.000000in}{0.000000in}}%
\pgfpathlineto{\pgfqpoint{0.000000in}{-0.048611in}}%
\pgfusepath{stroke,fill}%
}%
\begin{pgfscope}%
\pgfsys@transformshift{6.888622in}{0.484778in}%
\pgfsys@useobject{currentmarker}{}%
\end{pgfscope}%
\end{pgfscope}%
\begin{pgfscope}%
\pgftext[x=6.888622in,y=0.387555in,,top]{\rmfamily\fontsize{20.000000}{24.000000}\selectfont \(\displaystyle 1\)}%
\end{pgfscope}%
\begin{pgfscope}%
\pgfsetbuttcap%
\pgfsetroundjoin%
\definecolor{currentfill}{rgb}{0.000000,0.000000,0.000000}%
\pgfsetfillcolor{currentfill}%
\pgfsetlinewidth{0.803000pt}%
\definecolor{currentstroke}{rgb}{0.000000,0.000000,0.000000}%
\pgfsetstrokecolor{currentstroke}%
\pgfsetdash{}{0pt}%
\pgfsys@defobject{currentmarker}{\pgfqpoint{0.000000in}{-0.048611in}}{\pgfqpoint{0.000000in}{0.000000in}}{%
\pgfpathmoveto{\pgfqpoint{0.000000in}{0.000000in}}%
\pgfpathlineto{\pgfqpoint{0.000000in}{-0.048611in}}%
\pgfusepath{stroke,fill}%
}%
\begin{pgfscope}%
\pgfsys@transformshift{6.534038in}{0.484778in}%
\pgfsys@useobject{currentmarker}{}%
\end{pgfscope}%
\end{pgfscope}%
\begin{pgfscope}%
\pgftext[x=6.534038in,y=0.387555in,,top]{\rmfamily\fontsize{20.000000}{24.000000}\selectfont \(\displaystyle \widehat{c}_{\beta}\)}%
\end{pgfscope}%
\begin{pgfscope}%
\pgfsetbuttcap%
\pgfsetroundjoin%
\definecolor{currentfill}{rgb}{0.000000,0.000000,0.000000}%
\pgfsetfillcolor{currentfill}%
\pgfsetlinewidth{0.803000pt}%
\definecolor{currentstroke}{rgb}{0.000000,0.000000,0.000000}%
\pgfsetstrokecolor{currentstroke}%
\pgfsetdash{}{0pt}%
\pgfsys@defobject{currentmarker}{\pgfqpoint{0.000000in}{-0.048611in}}{\pgfqpoint{0.000000in}{0.000000in}}{%
\pgfpathmoveto{\pgfqpoint{0.000000in}{0.000000in}}%
\pgfpathlineto{\pgfqpoint{0.000000in}{-0.048611in}}%
\pgfusepath{stroke,fill}%
}%
\begin{pgfscope}%
\pgfsys@transformshift{0.688622in}{0.484778in}%
\pgfsys@useobject{currentmarker}{}%
\end{pgfscope}%
\end{pgfscope}%
\begin{pgfscope}%
\pgftext[x=0.688622in,y=0.387555in,,top]{\rmfamily\fontsize{20.000000}{24.000000}\selectfont \(\displaystyle c(\alpha,\beta)\)}%
\end{pgfscope}%
\begin{pgfscope}%
\pgfsetbuttcap%
\pgfsetroundjoin%
\definecolor{currentfill}{rgb}{0.000000,0.000000,0.000000}%
\pgfsetfillcolor{currentfill}%
\pgfsetlinewidth{0.803000pt}%
\definecolor{currentstroke}{rgb}{0.000000,0.000000,0.000000}%
\pgfsetstrokecolor{currentstroke}%
\pgfsetdash{}{0pt}%
\pgfsys@defobject{currentmarker}{\pgfqpoint{-0.048611in}{0.000000in}}{\pgfqpoint{0.000000in}{0.000000in}}{%
\pgfpathmoveto{\pgfqpoint{0.000000in}{0.000000in}}%
\pgfpathlineto{\pgfqpoint{-0.048611in}{0.000000in}}%
\pgfusepath{stroke,fill}%
}%
\begin{pgfscope}%
\pgfsys@transformshift{0.688622in}{0.484778in}%
\pgfsys@useobject{currentmarker}{}%
\end{pgfscope}%
\end{pgfscope}%
\begin{pgfscope}%
\pgftext[x=0.591400in,y=0.388389in,left,base]{\rmfamily\fontsize{20.000000}{24.000000}\selectfont  }%
\end{pgfscope}%
\begin{pgfscope}%
\pgfsetbuttcap%
\pgfsetroundjoin%
\definecolor{currentfill}{rgb}{0.000000,0.000000,0.000000}%
\pgfsetfillcolor{currentfill}%
\pgfsetlinewidth{0.803000pt}%
\definecolor{currentstroke}{rgb}{0.000000,0.000000,0.000000}%
\pgfsetstrokecolor{currentstroke}%
\pgfsetdash{}{0pt}%
\pgfsys@defobject{currentmarker}{\pgfqpoint{-0.048611in}{0.000000in}}{\pgfqpoint{0.000000in}{0.000000in}}{%
\pgfpathmoveto{\pgfqpoint{0.000000in}{0.000000in}}%
\pgfpathlineto{\pgfqpoint{-0.048611in}{0.000000in}}%
\pgfusepath{stroke,fill}%
}%
\begin{pgfscope}%
\pgfsys@transformshift{0.688622in}{1.390778in}%
\pgfsys@useobject{currentmarker}{}%
\end{pgfscope}%
\end{pgfscope}%
\begin{pgfscope}%
\pgftext[x=0.591400in,y=1.294389in,left,base]{\rmfamily\fontsize{20.000000}{24.000000}\selectfont  }%
\end{pgfscope}%
\begin{pgfscope}%
\pgfsetbuttcap%
\pgfsetroundjoin%
\definecolor{currentfill}{rgb}{0.000000,0.000000,0.000000}%
\pgfsetfillcolor{currentfill}%
\pgfsetlinewidth{0.803000pt}%
\definecolor{currentstroke}{rgb}{0.000000,0.000000,0.000000}%
\pgfsetstrokecolor{currentstroke}%
\pgfsetdash{}{0pt}%
\pgfsys@defobject{currentmarker}{\pgfqpoint{-0.048611in}{0.000000in}}{\pgfqpoint{0.000000in}{0.000000in}}{%
\pgfpathmoveto{\pgfqpoint{0.000000in}{0.000000in}}%
\pgfpathlineto{\pgfqpoint{-0.048611in}{0.000000in}}%
\pgfusepath{stroke,fill}%
}%
\begin{pgfscope}%
\pgfsys@transformshift{0.688622in}{2.296778in}%
\pgfsys@useobject{currentmarker}{}%
\end{pgfscope}%
\end{pgfscope}%
\begin{pgfscope}%
\pgftext[x=0.248838in,y=2.200389in,left,base]{\rmfamily\fontsize{20.000000}{24.000000}\selectfont \(\displaystyle 0.4\)}%
\end{pgfscope}%
\begin{pgfscope}%
\pgfsetbuttcap%
\pgfsetroundjoin%
\definecolor{currentfill}{rgb}{0.000000,0.000000,0.000000}%
\pgfsetfillcolor{currentfill}%
\pgfsetlinewidth{0.803000pt}%
\definecolor{currentstroke}{rgb}{0.000000,0.000000,0.000000}%
\pgfsetstrokecolor{currentstroke}%
\pgfsetdash{}{0pt}%
\pgfsys@defobject{currentmarker}{\pgfqpoint{-0.048611in}{0.000000in}}{\pgfqpoint{0.000000in}{0.000000in}}{%
\pgfpathmoveto{\pgfqpoint{0.000000in}{0.000000in}}%
\pgfpathlineto{\pgfqpoint{-0.048611in}{0.000000in}}%
\pgfusepath{stroke,fill}%
}%
\begin{pgfscope}%
\pgfsys@transformshift{0.688622in}{3.202777in}%
\pgfsys@useobject{currentmarker}{}%
\end{pgfscope}%
\end{pgfscope}%
\begin{pgfscope}%
\pgftext[x=0.248838in,y=3.106389in,left,base]{\rmfamily\fontsize{20.000000}{24.000000}\selectfont \(\displaystyle 0.6\)}%
\end{pgfscope}%
\begin{pgfscope}%
\pgfsetbuttcap%
\pgfsetroundjoin%
\definecolor{currentfill}{rgb}{0.000000,0.000000,0.000000}%
\pgfsetfillcolor{currentfill}%
\pgfsetlinewidth{0.803000pt}%
\definecolor{currentstroke}{rgb}{0.000000,0.000000,0.000000}%
\pgfsetstrokecolor{currentstroke}%
\pgfsetdash{}{0pt}%
\pgfsys@defobject{currentmarker}{\pgfqpoint{-0.048611in}{0.000000in}}{\pgfqpoint{0.000000in}{0.000000in}}{%
\pgfpathmoveto{\pgfqpoint{0.000000in}{0.000000in}}%
\pgfpathlineto{\pgfqpoint{-0.048611in}{0.000000in}}%
\pgfusepath{stroke,fill}%
}%
\begin{pgfscope}%
\pgfsys@transformshift{0.688622in}{4.108778in}%
\pgfsys@useobject{currentmarker}{}%
\end{pgfscope}%
\end{pgfscope}%
\begin{pgfscope}%
\pgftext[x=0.248838in,y=4.012389in,left,base]{\rmfamily\fontsize{20.000000}{24.000000}\selectfont \(\displaystyle 0.8\)}%
\end{pgfscope}%
\begin{pgfscope}%
\pgfsetbuttcap%
\pgfsetroundjoin%
\definecolor{currentfill}{rgb}{0.000000,0.000000,0.000000}%
\pgfsetfillcolor{currentfill}%
\pgfsetlinewidth{0.803000pt}%
\definecolor{currentstroke}{rgb}{0.000000,0.000000,0.000000}%
\pgfsetstrokecolor{currentstroke}%
\pgfsetdash{}{0pt}%
\pgfsys@defobject{currentmarker}{\pgfqpoint{-0.048611in}{0.000000in}}{\pgfqpoint{0.000000in}{0.000000in}}{%
\pgfpathmoveto{\pgfqpoint{0.000000in}{0.000000in}}%
\pgfpathlineto{\pgfqpoint{-0.048611in}{0.000000in}}%
\pgfusepath{stroke,fill}%
}%
\begin{pgfscope}%
\pgfsys@transformshift{0.688622in}{5.014778in}%
\pgfsys@useobject{currentmarker}{}%
\end{pgfscope}%
\end{pgfscope}%
\begin{pgfscope}%
\pgftext[x=0.459293in,y=4.918389in,left,base]{\rmfamily\fontsize{20.000000}{24.000000}\selectfont \(\displaystyle 1\)}%
\end{pgfscope}%
\begin{pgfscope}%
\pgfsetbuttcap%
\pgfsetroundjoin%
\definecolor{currentfill}{rgb}{0.000000,0.000000,0.000000}%
\pgfsetfillcolor{currentfill}%
\pgfsetlinewidth{0.803000pt}%
\definecolor{currentstroke}{rgb}{0.000000,0.000000,0.000000}%
\pgfsetstrokecolor{currentstroke}%
\pgfsetdash{}{0pt}%
\pgfsys@defobject{currentmarker}{\pgfqpoint{-0.048611in}{0.000000in}}{\pgfqpoint{0.000000in}{0.000000in}}{%
\pgfpathmoveto{\pgfqpoint{0.000000in}{0.000000in}}%
\pgfpathlineto{\pgfqpoint{-0.048611in}{0.000000in}}%
\pgfusepath{stroke,fill}%
}%
\begin{pgfscope}%
\pgfsys@transformshift{0.688622in}{0.859210in}%
\pgfsys@useobject{currentmarker}{}%
\end{pgfscope}%
\end{pgfscope}%
\begin{pgfscope}%
\pgftext[x=0.100000in,y=0.735298in,left,base]{\rmfamily\fontsize{20.000000}{24.000000}\selectfont \(\displaystyle F_{\alpha,\beta}^1\)}%
\end{pgfscope}%
\begin{pgfscope}%
\pgfsetbuttcap%
\pgfsetroundjoin%
\definecolor{currentfill}{rgb}{0.000000,0.000000,0.000000}%
\pgfsetfillcolor{currentfill}%
\pgfsetlinewidth{0.803000pt}%
\definecolor{currentstroke}{rgb}{0.000000,0.000000,0.000000}%
\pgfsetstrokecolor{currentstroke}%
\pgfsetdash{}{0pt}%
\pgfsys@defobject{currentmarker}{\pgfqpoint{-0.048611in}{0.000000in}}{\pgfqpoint{0.000000in}{0.000000in}}{%
\pgfpathmoveto{\pgfqpoint{0.000000in}{0.000000in}}%
\pgfpathlineto{\pgfqpoint{-0.048611in}{0.000000in}}%
\pgfusepath{stroke,fill}%
}%
\begin{pgfscope}%
\pgfsys@transformshift{0.688622in}{1.303009in}%
\pgfsys@useobject{currentmarker}{}%
\end{pgfscope}%
\end{pgfscope}%
\begin{pgfscope}%
\pgftext[x=0.100000in,y=1.179096in,left,base]{\rmfamily\fontsize{20.000000}{24.000000}\selectfont \(\displaystyle F_{\alpha,\beta}^0\)}%
\end{pgfscope}%
\begin{pgfscope}%
\pgfpathrectangle{\pgfqpoint{0.688622in}{0.484778in}}{\pgfqpoint{6.200000in}{4.530000in}} %
\pgfusepath{clip}%
\pgfsetbuttcap%
\pgfsetroundjoin%
\pgfsetlinewidth{1.505625pt}%
\definecolor{currentstroke}{rgb}{0.501961,0.501961,0.501961}%
\pgfsetstrokecolor{currentstroke}%
\pgfsetdash{{5.550000pt}{2.400000pt}}{0.000000pt}%
\pgfpathmoveto{\pgfqpoint{6.534038in}{0.484778in}}%
\pgfpathlineto{\pgfqpoint{6.534038in}{5.028666in}}%
\pgfusepath{stroke}%
\end{pgfscope}%
\begin{pgfscope}%
\pgfpathrectangle{\pgfqpoint{0.688622in}{0.484778in}}{\pgfqpoint{6.200000in}{4.530000in}} %
\pgfusepath{clip}%
\pgfsetbuttcap%
\pgfsetroundjoin%
\pgfsetlinewidth{1.003750pt}%
\definecolor{currentstroke}{rgb}{0.000000,0.392157,0.000000}%
\pgfsetstrokecolor{currentstroke}%
\pgfsetdash{{3.700000pt}{1.600000pt}}{0.000000pt}%
\pgfpathmoveto{\pgfqpoint{0.688622in}{0.859210in}}%
\pgfpathlineto{\pgfqpoint{6.888622in}{0.859210in}}%
\pgfusepath{stroke}%
\end{pgfscope}%
\begin{pgfscope}%
\pgfpathrectangle{\pgfqpoint{0.688622in}{0.484778in}}{\pgfqpoint{6.200000in}{4.530000in}} %
\pgfusepath{clip}%
\pgfsetbuttcap%
\pgfsetroundjoin%
\pgfsetlinewidth{1.505625pt}%
\definecolor{currentstroke}{rgb}{0.501961,0.501961,0.501961}%
\pgfsetstrokecolor{currentstroke}%
\pgfsetdash{{5.550000pt}{2.400000pt}}{0.000000pt}%
\pgfpathmoveto{\pgfqpoint{0.688622in}{0.484778in}}%
\pgfpathlineto{\pgfqpoint{0.688622in}{0.859210in}}%
\pgfusepath{stroke}%
\end{pgfscope}%
\begin{pgfscope}%
\pgfpathrectangle{\pgfqpoint{0.688622in}{0.484778in}}{\pgfqpoint{6.200000in}{4.530000in}} %
\pgfusepath{clip}%
\pgfsetrectcap%
\pgfsetroundjoin%
\pgfsetlinewidth{2.007500pt}%
\definecolor{currentstroke}{rgb}{0.721569,0.525490,0.043137}%
\pgfsetstrokecolor{currentstroke}%
\pgfsetdash{}{0pt}%
\pgfpathmoveto{\pgfqpoint{0.688622in}{1.303009in}}%
\pgfpathlineto{\pgfqpoint{1.035822in}{1.304527in}}%
\pgfpathlineto{\pgfqpoint{1.383022in}{1.309082in}}%
\pgfpathlineto{\pgfqpoint{1.730222in}{1.316673in}}%
\pgfpathlineto{\pgfqpoint{2.077422in}{1.327301in}}%
\pgfpathlineto{\pgfqpoint{2.424622in}{1.340966in}}%
\pgfpathlineto{\pgfqpoint{2.771822in}{1.357667in}}%
\pgfpathlineto{\pgfqpoint{3.119022in}{1.377404in}}%
\pgfpathlineto{\pgfqpoint{3.472422in}{1.400613in}}%
\pgfpathlineto{\pgfqpoint{3.825822in}{1.426967in}}%
\pgfpathlineto{\pgfqpoint{4.179222in}{1.456467in}}%
\pgfpathlineto{\pgfqpoint{4.532622in}{1.489114in}}%
\pgfpathlineto{\pgfqpoint{4.886022in}{1.524906in}}%
\pgfpathlineto{\pgfqpoint{5.239422in}{1.563844in}}%
\pgfpathlineto{\pgfqpoint{5.592822in}{1.605928in}}%
\pgfpathlineto{\pgfqpoint{5.946222in}{1.651158in}}%
\pgfpathlineto{\pgfqpoint{6.299622in}{1.699535in}}%
\pgfpathlineto{\pgfqpoint{6.535222in}{1.733533in}}%
\pgfpathlineto{\pgfqpoint{6.534038in}{1.733359in}}%
\pgfusepath{stroke}%
\end{pgfscope}%
\begin{pgfscope}%
\pgfpathrectangle{\pgfqpoint{0.688622in}{0.484778in}}{\pgfqpoint{6.200000in}{4.530000in}} %
\pgfusepath{clip}%
\pgfsetbuttcap%
\pgfsetroundjoin%
\pgfsetlinewidth{3.011250pt}%
\definecolor{currentstroke}{rgb}{0.000000,0.000000,1.000000}%
\pgfsetstrokecolor{currentstroke}%
\pgfsetdash{{19.200000pt}{4.800000pt}{3.000000pt}{4.800000pt}}{0.000000pt}%
\pgfpathmoveto{\pgfqpoint{6.534038in}{1.733359in}}%
\pgfpathlineto{\pgfqpoint{6.541422in}{1.983693in}}%
\pgfpathlineto{\pgfqpoint{6.547622in}{2.085139in}}%
\pgfpathlineto{\pgfqpoint{6.553822in}{2.170019in}}%
\pgfpathlineto{\pgfqpoint{6.560022in}{2.246034in}}%
\pgfpathlineto{\pgfqpoint{6.566222in}{2.316375in}}%
\pgfpathlineto{\pgfqpoint{6.572422in}{2.382726in}}%
\pgfpathlineto{\pgfqpoint{6.578622in}{2.446104in}}%
\pgfpathlineto{\pgfqpoint{6.584822in}{2.507180in}}%
\pgfpathlineto{\pgfqpoint{6.591022in}{2.566418in}}%
\pgfpathlineto{\pgfqpoint{6.597222in}{2.624162in}}%
\pgfpathlineto{\pgfqpoint{6.603422in}{2.680668in}}%
\pgfpathlineto{\pgfqpoint{6.609622in}{2.736136in}}%
\pgfpathlineto{\pgfqpoint{6.615822in}{2.790727in}}%
\pgfpathlineto{\pgfqpoint{6.622022in}{2.844569in}}%
\pgfpathlineto{\pgfqpoint{6.628222in}{2.897769in}}%
\pgfpathlineto{\pgfqpoint{6.634422in}{2.950415in}}%
\pgfpathlineto{\pgfqpoint{6.640622in}{3.002583in}}%
\pgfpathlineto{\pgfqpoint{6.646822in}{3.054335in}}%
\pgfpathlineto{\pgfqpoint{6.653022in}{3.105727in}}%
\pgfpathlineto{\pgfqpoint{6.659222in}{3.156808in}}%
\pgfpathlineto{\pgfqpoint{6.665422in}{3.207618in}}%
\pgfpathlineto{\pgfqpoint{6.671622in}{3.258194in}}%
\pgfpathlineto{\pgfqpoint{6.677822in}{3.308570in}}%
\pgfpathlineto{\pgfqpoint{6.684022in}{3.358773in}}%
\pgfpathlineto{\pgfqpoint{6.690222in}{3.408831in}}%
\pgfpathlineto{\pgfqpoint{6.696422in}{3.458767in}}%
\pgfpathlineto{\pgfqpoint{6.702622in}{3.508601in}}%
\pgfpathlineto{\pgfqpoint{6.708822in}{3.558353in}}%
\pgfpathlineto{\pgfqpoint{6.715022in}{3.608040in}}%
\pgfpathlineto{\pgfqpoint{6.721222in}{3.657678in}}%
\pgfpathlineto{\pgfqpoint{6.727422in}{3.707283in}}%
\pgfpathlineto{\pgfqpoint{6.733622in}{3.756866in}}%
\pgfpathlineto{\pgfqpoint{6.739822in}{3.806441in}}%
\pgfpathlineto{\pgfqpoint{6.746022in}{3.856019in}}%
\pgfpathlineto{\pgfqpoint{6.752222in}{3.905611in}}%
\pgfpathlineto{\pgfqpoint{6.758422in}{3.955226in}}%
\pgfpathlineto{\pgfqpoint{6.764622in}{4.004874in}}%
\pgfpathlineto{\pgfqpoint{6.770822in}{4.054563in}}%
\pgfpathlineto{\pgfqpoint{6.777022in}{4.104301in}}%
\pgfpathlineto{\pgfqpoint{6.783222in}{4.154096in}}%
\pgfpathlineto{\pgfqpoint{6.789422in}{4.203954in}}%
\pgfpathlineto{\pgfqpoint{6.795622in}{4.253883in}}%
\pgfpathlineto{\pgfqpoint{6.801822in}{4.303888in}}%
\pgfpathlineto{\pgfqpoint{6.808022in}{4.353976in}}%
\pgfpathlineto{\pgfqpoint{6.814222in}{4.404151in}}%
\pgfpathlineto{\pgfqpoint{6.820422in}{4.454419in}}%
\pgfpathlineto{\pgfqpoint{6.826622in}{4.504786in}}%
\pgfpathlineto{\pgfqpoint{6.832822in}{4.555255in}}%
\pgfpathlineto{\pgfqpoint{6.839022in}{4.605831in}}%
\pgfpathlineto{\pgfqpoint{6.845222in}{4.656519in}}%
\pgfpathlineto{\pgfqpoint{6.851422in}{4.707323in}}%
\pgfpathlineto{\pgfqpoint{6.857622in}{4.758246in}}%
\pgfpathlineto{\pgfqpoint{6.863822in}{4.809292in}}%
\pgfpathlineto{\pgfqpoint{6.870022in}{4.860465in}}%
\pgfpathlineto{\pgfqpoint{6.876222in}{4.911768in}}%
\pgfpathlineto{\pgfqpoint{6.882422in}{4.963205in}}%
\pgfpathlineto{\pgfqpoint{6.888622in}{5.014778in}}%
\pgfusepath{stroke}%
\end{pgfscope}%
\begin{pgfscope}%
\pgfpathrectangle{\pgfqpoint{0.688622in}{0.484778in}}{\pgfqpoint{6.200000in}{4.530000in}} %
\pgfusepath{clip}%
\pgfsetbuttcap%
\pgfsetroundjoin%
\pgfsetlinewidth{2.509375pt}%
\definecolor{currentstroke}{rgb}{1.000000,0.000000,0.000000}%
\pgfsetstrokecolor{currentstroke}%
\pgfsetdash{{9.250000pt}{4.000000pt}}{0.000000pt}%
\pgfpathmoveto{\pgfqpoint{6.534038in}{1.733359in}}%
\pgfpathlineto{\pgfqpoint{6.541422in}{1.526646in}}%
\pgfpathlineto{\pgfqpoint{6.547622in}{1.462032in}}%
\pgfpathlineto{\pgfqpoint{6.553822in}{1.414170in}}%
\pgfpathlineto{\pgfqpoint{6.560022in}{1.375360in}}%
\pgfpathlineto{\pgfqpoint{6.566222in}{1.342414in}}%
\pgfpathlineto{\pgfqpoint{6.572422in}{1.313645in}}%
\pgfpathlineto{\pgfqpoint{6.578622in}{1.288037in}}%
\pgfpathlineto{\pgfqpoint{6.584822in}{1.264922in}}%
\pgfpathlineto{\pgfqpoint{6.591022in}{1.243833in}}%
\pgfpathlineto{\pgfqpoint{6.597222in}{1.224430in}}%
\pgfpathlineto{\pgfqpoint{6.603422in}{1.206455in}}%
\pgfpathlineto{\pgfqpoint{6.609622in}{1.189708in}}%
\pgfpathlineto{\pgfqpoint{6.615822in}{1.174031in}}%
\pgfpathlineto{\pgfqpoint{6.622022in}{1.159295in}}%
\pgfpathlineto{\pgfqpoint{6.628222in}{1.145393in}}%
\pgfpathlineto{\pgfqpoint{6.634422in}{1.132239in}}%
\pgfpathlineto{\pgfqpoint{6.640622in}{1.119756in}}%
\pgfpathlineto{\pgfqpoint{6.646822in}{1.107883in}}%
\pgfpathlineto{\pgfqpoint{6.653022in}{1.096564in}}%
\pgfpathlineto{\pgfqpoint{6.659222in}{1.085753in}}%
\pgfpathlineto{\pgfqpoint{6.665422in}{1.075407in}}%
\pgfpathlineto{\pgfqpoint{6.671622in}{1.065490in}}%
\pgfpathlineto{\pgfqpoint{6.677822in}{1.055971in}}%
\pgfpathlineto{\pgfqpoint{6.684022in}{1.046821in}}%
\pgfpathlineto{\pgfqpoint{6.690222in}{1.038014in}}%
\pgfpathlineto{\pgfqpoint{6.696422in}{1.029527in}}%
\pgfpathlineto{\pgfqpoint{6.702622in}{1.021340in}}%
\pgfpathlineto{\pgfqpoint{6.708822in}{1.013434in}}%
\pgfpathlineto{\pgfqpoint{6.715022in}{1.005792in}}%
\pgfpathlineto{\pgfqpoint{6.721222in}{0.998398in}}%
\pgfpathlineto{\pgfqpoint{6.727422in}{0.991239in}}%
\pgfpathlineto{\pgfqpoint{6.733622in}{0.984301in}}%
\pgfpathlineto{\pgfqpoint{6.739822in}{0.977574in}}%
\pgfpathlineto{\pgfqpoint{6.746022in}{0.971044in}}%
\pgfpathlineto{\pgfqpoint{6.752222in}{0.964704in}}%
\pgfpathlineto{\pgfqpoint{6.758422in}{0.958542in}}%
\pgfpathlineto{\pgfqpoint{6.764622in}{0.952552in}}%
\pgfpathlineto{\pgfqpoint{6.770822in}{0.946723in}}%
\pgfpathlineto{\pgfqpoint{6.777022in}{0.941050in}}%
\pgfpathlineto{\pgfqpoint{6.783222in}{0.935525in}}%
\pgfpathlineto{\pgfqpoint{6.789422in}{0.930141in}}%
\pgfpathlineto{\pgfqpoint{6.795622in}{0.924893in}}%
\pgfpathlineto{\pgfqpoint{6.801822in}{0.919774in}}%
\pgfpathlineto{\pgfqpoint{6.808022in}{0.914779in}}%
\pgfpathlineto{\pgfqpoint{6.814222in}{0.909904in}}%
\pgfpathlineto{\pgfqpoint{6.820422in}{0.905143in}}%
\pgfpathlineto{\pgfqpoint{6.826622in}{0.900493in}}%
\pgfpathlineto{\pgfqpoint{6.832822in}{0.895948in}}%
\pgfpathlineto{\pgfqpoint{6.839022in}{0.891504in}}%
\pgfpathlineto{\pgfqpoint{6.845222in}{0.887159in}}%
\pgfpathlineto{\pgfqpoint{6.851422in}{0.882909in}}%
\pgfpathlineto{\pgfqpoint{6.857622in}{0.878749in}}%
\pgfpathlineto{\pgfqpoint{6.863822in}{0.874677in}}%
\pgfpathlineto{\pgfqpoint{6.870022in}{0.870690in}}%
\pgfpathlineto{\pgfqpoint{6.876222in}{0.866785in}}%
\pgfpathlineto{\pgfqpoint{6.882422in}{0.862960in}}%
\pgfpathlineto{\pgfqpoint{6.888622in}{0.859210in}}%
\pgfusepath{stroke}%
\end{pgfscope}%
\begin{pgfscope}%
\pgfpathrectangle{\pgfqpoint{0.688622in}{0.484778in}}{\pgfqpoint{6.200000in}{4.530000in}} %
\pgfusepath{clip}%
\pgfsetbuttcap%
\pgfsetroundjoin%
\definecolor{currentfill}{rgb}{0.721569,0.525490,0.043137}%
\pgfsetfillcolor{currentfill}%
\pgfsetlinewidth{1.003750pt}%
\definecolor{currentstroke}{rgb}{0.721569,0.525490,0.043137}%
\pgfsetstrokecolor{currentstroke}%
\pgfsetdash{}{0pt}%
\pgfsys@defobject{currentmarker}{\pgfqpoint{-0.027778in}{-0.027778in}}{\pgfqpoint{0.027778in}{0.027778in}}{%
\pgfpathmoveto{\pgfqpoint{0.000000in}{-0.027778in}}%
\pgfpathcurveto{\pgfqpoint{0.007367in}{-0.027778in}}{\pgfqpoint{0.014433in}{-0.024851in}}{\pgfqpoint{0.019642in}{-0.019642in}}%
\pgfpathcurveto{\pgfqpoint{0.024851in}{-0.014433in}}{\pgfqpoint{0.027778in}{-0.007367in}}{\pgfqpoint{0.027778in}{0.000000in}}%
\pgfpathcurveto{\pgfqpoint{0.027778in}{0.007367in}}{\pgfqpoint{0.024851in}{0.014433in}}{\pgfqpoint{0.019642in}{0.019642in}}%
\pgfpathcurveto{\pgfqpoint{0.014433in}{0.024851in}}{\pgfqpoint{0.007367in}{0.027778in}}{\pgfqpoint{0.000000in}{0.027778in}}%
\pgfpathcurveto{\pgfqpoint{-0.007367in}{0.027778in}}{\pgfqpoint{-0.014433in}{0.024851in}}{\pgfqpoint{-0.019642in}{0.019642in}}%
\pgfpathcurveto{\pgfqpoint{-0.024851in}{0.014433in}}{\pgfqpoint{-0.027778in}{0.007367in}}{\pgfqpoint{-0.027778in}{0.000000in}}%
\pgfpathcurveto{\pgfqpoint{-0.027778in}{-0.007367in}}{\pgfqpoint{-0.024851in}{-0.014433in}}{\pgfqpoint{-0.019642in}{-0.019642in}}%
\pgfpathcurveto{\pgfqpoint{-0.014433in}{-0.024851in}}{\pgfqpoint{-0.007367in}{-0.027778in}}{\pgfqpoint{0.000000in}{-0.027778in}}%
\pgfpathclose%
\pgfusepath{stroke,fill}%
}%
\begin{pgfscope}%
\pgfsys@transformshift{0.694822in}{1.303009in}%
\pgfsys@useobject{currentmarker}{}%
\end{pgfscope}%
\end{pgfscope}%
\begin{pgfscope}%
\pgfpathrectangle{\pgfqpoint{0.688622in}{0.484778in}}{\pgfqpoint{6.200000in}{4.530000in}} %
\pgfusepath{clip}%
\pgfsetbuttcap%
\pgfsetroundjoin%
\definecolor{currentfill}{rgb}{0.000000,0.000000,1.000000}%
\pgfsetfillcolor{currentfill}%
\pgfsetlinewidth{1.003750pt}%
\definecolor{currentstroke}{rgb}{0.000000,0.000000,1.000000}%
\pgfsetstrokecolor{currentstroke}%
\pgfsetdash{}{0pt}%
\pgfsys@defobject{currentmarker}{\pgfqpoint{-0.027778in}{-0.027778in}}{\pgfqpoint{0.027778in}{0.027778in}}{%
\pgfpathmoveto{\pgfqpoint{0.000000in}{-0.027778in}}%
\pgfpathcurveto{\pgfqpoint{0.007367in}{-0.027778in}}{\pgfqpoint{0.014433in}{-0.024851in}}{\pgfqpoint{0.019642in}{-0.019642in}}%
\pgfpathcurveto{\pgfqpoint{0.024851in}{-0.014433in}}{\pgfqpoint{0.027778in}{-0.007367in}}{\pgfqpoint{0.027778in}{0.000000in}}%
\pgfpathcurveto{\pgfqpoint{0.027778in}{0.007367in}}{\pgfqpoint{0.024851in}{0.014433in}}{\pgfqpoint{0.019642in}{0.019642in}}%
\pgfpathcurveto{\pgfqpoint{0.014433in}{0.024851in}}{\pgfqpoint{0.007367in}{0.027778in}}{\pgfqpoint{0.000000in}{0.027778in}}%
\pgfpathcurveto{\pgfqpoint{-0.007367in}{0.027778in}}{\pgfqpoint{-0.014433in}{0.024851in}}{\pgfqpoint{-0.019642in}{0.019642in}}%
\pgfpathcurveto{\pgfqpoint{-0.024851in}{0.014433in}}{\pgfqpoint{-0.027778in}{0.007367in}}{\pgfqpoint{-0.027778in}{0.000000in}}%
\pgfpathcurveto{\pgfqpoint{-0.027778in}{-0.007367in}}{\pgfqpoint{-0.024851in}{-0.014433in}}{\pgfqpoint{-0.019642in}{-0.019642in}}%
\pgfpathcurveto{\pgfqpoint{-0.014433in}{-0.024851in}}{\pgfqpoint{-0.007367in}{-0.027778in}}{\pgfqpoint{0.000000in}{-0.027778in}}%
\pgfpathclose%
\pgfusepath{stroke,fill}%
}%
\begin{pgfscope}%
\pgfsys@transformshift{6.882422in}{5.014778in}%
\pgfsys@useobject{currentmarker}{}%
\end{pgfscope}%
\end{pgfscope}%
\begin{pgfscope}%
\pgfpathrectangle{\pgfqpoint{0.688622in}{0.484778in}}{\pgfqpoint{6.200000in}{4.530000in}} %
\pgfusepath{clip}%
\pgfsetbuttcap%
\pgfsetroundjoin%
\definecolor{currentfill}{rgb}{1.000000,0.000000,0.000000}%
\pgfsetfillcolor{currentfill}%
\pgfsetlinewidth{1.003750pt}%
\definecolor{currentstroke}{rgb}{1.000000,0.000000,0.000000}%
\pgfsetstrokecolor{currentstroke}%
\pgfsetdash{}{0pt}%
\pgfsys@defobject{currentmarker}{\pgfqpoint{-0.027778in}{-0.027778in}}{\pgfqpoint{0.027778in}{0.027778in}}{%
\pgfpathmoveto{\pgfqpoint{0.000000in}{-0.027778in}}%
\pgfpathcurveto{\pgfqpoint{0.007367in}{-0.027778in}}{\pgfqpoint{0.014433in}{-0.024851in}}{\pgfqpoint{0.019642in}{-0.019642in}}%
\pgfpathcurveto{\pgfqpoint{0.024851in}{-0.014433in}}{\pgfqpoint{0.027778in}{-0.007367in}}{\pgfqpoint{0.027778in}{0.000000in}}%
\pgfpathcurveto{\pgfqpoint{0.027778in}{0.007367in}}{\pgfqpoint{0.024851in}{0.014433in}}{\pgfqpoint{0.019642in}{0.019642in}}%
\pgfpathcurveto{\pgfqpoint{0.014433in}{0.024851in}}{\pgfqpoint{0.007367in}{0.027778in}}{\pgfqpoint{0.000000in}{0.027778in}}%
\pgfpathcurveto{\pgfqpoint{-0.007367in}{0.027778in}}{\pgfqpoint{-0.014433in}{0.024851in}}{\pgfqpoint{-0.019642in}{0.019642in}}%
\pgfpathcurveto{\pgfqpoint{-0.024851in}{0.014433in}}{\pgfqpoint{-0.027778in}{0.007367in}}{\pgfqpoint{-0.027778in}{0.000000in}}%
\pgfpathcurveto{\pgfqpoint{-0.027778in}{-0.007367in}}{\pgfqpoint{-0.024851in}{-0.014433in}}{\pgfqpoint{-0.019642in}{-0.019642in}}%
\pgfpathcurveto{\pgfqpoint{-0.014433in}{-0.024851in}}{\pgfqpoint{-0.007367in}{-0.027778in}}{\pgfqpoint{0.000000in}{-0.027778in}}%
\pgfpathclose%
\pgfusepath{stroke,fill}%
}%
\begin{pgfscope}%
\pgfsys@transformshift{6.882422in}{0.859210in}%
\pgfsys@useobject{currentmarker}{}%
\end{pgfscope}%
\end{pgfscope}%
\begin{pgfscope}%
\pgfpathrectangle{\pgfqpoint{0.688622in}{0.484778in}}{\pgfqpoint{6.200000in}{4.530000in}} %
\pgfusepath{clip}%
\pgfsetbuttcap%
\pgfsetroundjoin%
\definecolor{currentfill}{rgb}{0.721569,0.525490,0.043137}%
\pgfsetfillcolor{currentfill}%
\pgfsetlinewidth{1.003750pt}%
\definecolor{currentstroke}{rgb}{0.721569,0.525490,0.043137}%
\pgfsetstrokecolor{currentstroke}%
\pgfsetdash{}{0pt}%
\pgfsys@defobject{currentmarker}{\pgfqpoint{-0.027778in}{-0.027778in}}{\pgfqpoint{0.027778in}{0.027778in}}{%
\pgfpathmoveto{\pgfqpoint{0.000000in}{-0.027778in}}%
\pgfpathcurveto{\pgfqpoint{0.007367in}{-0.027778in}}{\pgfqpoint{0.014433in}{-0.024851in}}{\pgfqpoint{0.019642in}{-0.019642in}}%
\pgfpathcurveto{\pgfqpoint{0.024851in}{-0.014433in}}{\pgfqpoint{0.027778in}{-0.007367in}}{\pgfqpoint{0.027778in}{0.000000in}}%
\pgfpathcurveto{\pgfqpoint{0.027778in}{0.007367in}}{\pgfqpoint{0.024851in}{0.014433in}}{\pgfqpoint{0.019642in}{0.019642in}}%
\pgfpathcurveto{\pgfqpoint{0.014433in}{0.024851in}}{\pgfqpoint{0.007367in}{0.027778in}}{\pgfqpoint{0.000000in}{0.027778in}}%
\pgfpathcurveto{\pgfqpoint{-0.007367in}{0.027778in}}{\pgfqpoint{-0.014433in}{0.024851in}}{\pgfqpoint{-0.019642in}{0.019642in}}%
\pgfpathcurveto{\pgfqpoint{-0.024851in}{0.014433in}}{\pgfqpoint{-0.027778in}{0.007367in}}{\pgfqpoint{-0.027778in}{0.000000in}}%
\pgfpathcurveto{\pgfqpoint{-0.027778in}{-0.007367in}}{\pgfqpoint{-0.024851in}{-0.014433in}}{\pgfqpoint{-0.019642in}{-0.019642in}}%
\pgfpathcurveto{\pgfqpoint{-0.014433in}{-0.024851in}}{\pgfqpoint{-0.007367in}{-0.027778in}}{\pgfqpoint{0.000000in}{-0.027778in}}%
\pgfpathclose%
\pgfusepath{stroke,fill}%
}%
\begin{pgfscope}%
\pgfsys@transformshift{0.688622in}{0.859210in}%
\pgfsys@useobject{currentmarker}{}%
\end{pgfscope}%
\end{pgfscope}%
\begin{pgfscope}%
\pgfsetrectcap%
\pgfsetmiterjoin%
\pgfsetlinewidth{0.803000pt}%
\definecolor{currentstroke}{rgb}{0.000000,0.000000,0.000000}%
\pgfsetstrokecolor{currentstroke}%
\pgfsetdash{}{0pt}%
\pgfpathmoveto{\pgfqpoint{0.688622in}{0.484778in}}%
\pgfpathlineto{\pgfqpoint{0.688622in}{5.014778in}}%
\pgfusepath{stroke}%
\end{pgfscope}%
\begin{pgfscope}%
\pgfsetrectcap%
\pgfsetmiterjoin%
\pgfsetlinewidth{0.803000pt}%
\definecolor{currentstroke}{rgb}{0.000000,0.000000,0.000000}%
\pgfsetstrokecolor{currentstroke}%
\pgfsetdash{}{0pt}%
\pgfpathmoveto{\pgfqpoint{6.888622in}{0.484778in}}%
\pgfpathlineto{\pgfqpoint{6.888622in}{5.014778in}}%
\pgfusepath{stroke}%
\end{pgfscope}%
\begin{pgfscope}%
\pgfsetrectcap%
\pgfsetmiterjoin%
\pgfsetlinewidth{0.803000pt}%
\definecolor{currentstroke}{rgb}{0.000000,0.000000,0.000000}%
\pgfsetstrokecolor{currentstroke}%
\pgfsetdash{}{0pt}%
\pgfpathmoveto{\pgfqpoint{0.688622in}{0.484778in}}%
\pgfpathlineto{\pgfqpoint{6.888622in}{0.484778in}}%
\pgfusepath{stroke}%
\end{pgfscope}%
\begin{pgfscope}%
\pgfsetrectcap%
\pgfsetmiterjoin%
\pgfsetlinewidth{0.803000pt}%
\definecolor{currentstroke}{rgb}{0.000000,0.000000,0.000000}%
\pgfsetstrokecolor{currentstroke}%
\pgfsetdash{}{0pt}%
\pgfpathmoveto{\pgfqpoint{0.688622in}{5.014778in}}%
\pgfpathlineto{\pgfqpoint{6.888622in}{5.014778in}}%
\pgfusepath{stroke}%
\end{pgfscope}%
\begin{pgfscope}%
\definecolor{textcolor}{rgb}{0.000000,0.392157,0.000000}%
\pgfsetstrokecolor{textcolor}%
\pgfsetfillcolor{textcolor}%
\pgftext[x=2.920622in,y=0.949810in,left,base]{\color{textcolor}\rmfamily\fontsize{20.000000}{24.000000}\selectfont \(\displaystyle y=f_{\alpha,\beta,1}(1)\)}%
\end{pgfscope}%
\begin{pgfscope}%
\pgfsetbuttcap%
\pgfsetmiterjoin%
\definecolor{currentfill}{rgb}{1.000000,1.000000,1.000000}%
\pgfsetfillcolor{currentfill}%
\pgfsetlinewidth{1.003750pt}%
\definecolor{currentstroke}{rgb}{0.000000,0.000000,0.000000}%
\pgfsetstrokecolor{currentstroke}%
\pgfsetdash{}{0pt}%
\pgfpathmoveto{\pgfqpoint{2.855155in}{5.166235in}}%
\pgfpathlineto{\pgfqpoint{4.722089in}{5.166235in}}%
\pgfpathlineto{\pgfqpoint{4.722089in}{5.490679in}}%
\pgfpathlineto{\pgfqpoint{2.855155in}{5.490679in}}%
\pgfpathclose%
\pgfusepath{stroke,fill}%
\end{pgfscope}%
\begin{pgfscope}%
\pgftext[x=3.788622in,y=5.268458in,,base]{\rmfamily\fontsize{16.000000}{19.200000}\selectfont \(\displaystyle \alpha=\) 0.95, \(\displaystyle \beta=\) 0.75}%
\end{pgfscope}%
\begin{pgfscope}%
\pgfsetbuttcap%
\pgfsetmiterjoin%
\definecolor{currentfill}{rgb}{1.000000,1.000000,1.000000}%
\pgfsetfillcolor{currentfill}%
\pgfsetfillopacity{0.800000}%
\pgfsetlinewidth{1.003750pt}%
\definecolor{currentstroke}{rgb}{0.800000,0.800000,0.800000}%
\pgfsetstrokecolor{currentstroke}%
\pgfsetstrokeopacity{0.800000}%
\pgfsetdash{}{0pt}%
\pgfpathmoveto{\pgfqpoint{0.868178in}{3.429277in}}%
\pgfpathlineto{\pgfqpoint{3.402511in}{3.429277in}}%
\pgfpathquadraticcurveto{\pgfqpoint{3.458067in}{3.429277in}}{\pgfqpoint{3.458067in}{3.484833in}}%
\pgfpathlineto{\pgfqpoint{3.458067in}{4.767360in}}%
\pgfpathquadraticcurveto{\pgfqpoint{3.458067in}{4.822915in}}{\pgfqpoint{3.402511in}{4.822915in}}%
\pgfpathlineto{\pgfqpoint{0.868178in}{4.822915in}}%
\pgfpathquadraticcurveto{\pgfqpoint{0.812622in}{4.822915in}}{\pgfqpoint{0.812622in}{4.767360in}}%
\pgfpathlineto{\pgfqpoint{0.812622in}{3.484833in}}%
\pgfpathquadraticcurveto{\pgfqpoint{0.812622in}{3.429277in}}{\pgfqpoint{0.868178in}{3.429277in}}%
\pgfpathclose%
\pgfusepath{stroke,fill}%
\end{pgfscope}%
\begin{pgfscope}%
\pgfsetrectcap%
\pgfsetroundjoin%
\pgfsetlinewidth{2.007500pt}%
\definecolor{currentstroke}{rgb}{0.721569,0.525490,0.043137}%
\pgfsetstrokecolor{currentstroke}%
\pgfsetdash{}{0pt}%
\pgfpathmoveto{\pgfqpoint{0.923733in}{4.593244in}}%
\pgfpathlineto{\pgfqpoint{1.290400in}{4.593244in}}%
\pgfusepath{stroke}%
\end{pgfscope}%
\begin{pgfscope}%
\pgftext[x=1.512622in,y=4.496022in,left,base]{\rmfamily\fontsize{20.000000}{24.000000}\selectfont \(\displaystyle f_{\alpha,\beta,r}(c)\), \(\displaystyle c< \widehat{c}_{\beta}\)}%
\end{pgfscope}%
\begin{pgfscope}%
\pgfsetbuttcap%
\pgfsetroundjoin%
\pgfsetlinewidth{3.011250pt}%
\definecolor{currentstroke}{rgb}{0.000000,0.000000,1.000000}%
\pgfsetstrokecolor{currentstroke}%
\pgfsetdash{{19.200000pt}{4.800000pt}{3.000000pt}{4.800000pt}}{0.000000pt}%
\pgfpathmoveto{\pgfqpoint{0.923733in}{4.156476in}}%
\pgfpathlineto{\pgfqpoint{1.290400in}{4.156476in}}%
\pgfusepath{stroke}%
\end{pgfscope}%
\begin{pgfscope}%
\pgftext[x=1.512622in,y=4.059254in,left,base]{\rmfamily\fontsize{20.000000}{24.000000}\selectfont \(\displaystyle f_{\alpha,\beta,2}(c)\), \(\displaystyle c\geq \widehat{c}_{\beta}\)}%
\end{pgfscope}%
\begin{pgfscope}%
\pgfsetbuttcap%
\pgfsetroundjoin%
\pgfsetlinewidth{2.509375pt}%
\definecolor{currentstroke}{rgb}{1.000000,0.000000,0.000000}%
\pgfsetstrokecolor{currentstroke}%
\pgfsetdash{{9.250000pt}{4.000000pt}}{0.000000pt}%
\pgfpathmoveto{\pgfqpoint{0.923733in}{3.719708in}}%
\pgfpathlineto{\pgfqpoint{1.290400in}{3.719708in}}%
\pgfusepath{stroke}%
\end{pgfscope}%
\begin{pgfscope}%
\pgftext[x=1.512622in,y=3.622486in,left,base]{\rmfamily\fontsize{20.000000}{24.000000}\selectfont \(\displaystyle f_{\alpha,\beta,1}(c)\), \(\displaystyle c\geq \widehat{c}_{\beta}\)}%
\end{pgfscope}%
\end{pgfpicture}%
\makeatother%
\endgroup%

%% file: Fig1c.pgf
\begingroup%
\makeatletter%
\begin{pgfpicture}%
\pgfpathrectangle{\pgfpointorigin}{\pgfqpoint{7.054676in}{5.535124in}}%
\pgfusepath{use as bounding box, clip}%
\begin{pgfscope}%
\pgfsetbuttcap%
\pgfsetmiterjoin%
\definecolor{currentfill}{rgb}{1.000000,1.000000,1.000000}%
\pgfsetfillcolor{currentfill}%
\pgfsetlinewidth{0.000000pt}%
\definecolor{currentstroke}{rgb}{1.000000,1.000000,1.000000}%
\pgfsetstrokecolor{currentstroke}%
\pgfsetdash{}{0pt}%
\pgfpathmoveto{\pgfqpoint{0.000000in}{0.000000in}}%
\pgfpathlineto{\pgfqpoint{7.054676in}{0.000000in}}%
\pgfpathlineto{\pgfqpoint{7.054676in}{5.535124in}}%
\pgfpathlineto{\pgfqpoint{0.000000in}{5.535124in}}%
\pgfpathclose%
\pgfusepath{fill}%
\end{pgfscope}%
\begin{pgfscope}%
\pgfsetbuttcap%
\pgfsetmiterjoin%
\definecolor{currentfill}{rgb}{1.000000,1.000000,1.000000}%
\pgfsetfillcolor{currentfill}%
\pgfsetlinewidth{0.000000pt}%
\definecolor{currentstroke}{rgb}{0.000000,0.000000,0.000000}%
\pgfsetstrokecolor{currentstroke}%
\pgfsetstrokeopacity{0.000000}%
\pgfsetdash{}{0pt}%
\pgfpathmoveto{\pgfqpoint{0.688622in}{0.484778in}}%
\pgfpathlineto{\pgfqpoint{6.888622in}{0.484778in}}%
\pgfpathlineto{\pgfqpoint{6.888622in}{5.014778in}}%
\pgfpathlineto{\pgfqpoint{0.688622in}{5.014778in}}%
\pgfpathclose%
\pgfusepath{fill}%
\end{pgfscope}%
\begin{pgfscope}%
\pgfsetbuttcap%
\pgfsetroundjoin%
\definecolor{currentfill}{rgb}{0.000000,0.000000,0.000000}%
\pgfsetfillcolor{currentfill}%
\pgfsetlinewidth{0.803000pt}%
\definecolor{currentstroke}{rgb}{0.000000,0.000000,0.000000}%
\pgfsetstrokecolor{currentstroke}%
\pgfsetdash{}{0pt}%
\pgfsys@defobject{currentmarker}{\pgfqpoint{0.000000in}{-0.048611in}}{\pgfqpoint{0.000000in}{0.000000in}}{%
\pgfpathmoveto{\pgfqpoint{0.000000in}{0.000000in}}%
\pgfpathlineto{\pgfqpoint{0.000000in}{-0.048611in}}%
\pgfusepath{stroke,fill}%
}%
\begin{pgfscope}%
\pgfsys@transformshift{0.688622in}{0.484778in}%
\pgfsys@useobject{currentmarker}{}%
\end{pgfscope}%
\end{pgfscope}%
\begin{pgfscope}%
\pgftext[x=0.688622in,y=0.387555in,,top]{\rmfamily\fontsize{20.000000}{24.000000}\selectfont \(\displaystyle 0\)}%
\end{pgfscope}%
\begin{pgfscope}%
\pgfsetbuttcap%
\pgfsetroundjoin%
\definecolor{currentfill}{rgb}{0.000000,0.000000,0.000000}%
\pgfsetfillcolor{currentfill}%
\pgfsetlinewidth{0.803000pt}%
\definecolor{currentstroke}{rgb}{0.000000,0.000000,0.000000}%
\pgfsetstrokecolor{currentstroke}%
\pgfsetdash{}{0pt}%
\pgfsys@defobject{currentmarker}{\pgfqpoint{0.000000in}{-0.048611in}}{\pgfqpoint{0.000000in}{0.000000in}}{%
\pgfpathmoveto{\pgfqpoint{0.000000in}{0.000000in}}%
\pgfpathlineto{\pgfqpoint{0.000000in}{-0.048611in}}%
\pgfusepath{stroke,fill}%
}%
\begin{pgfscope}%
\pgfsys@transformshift{1.928622in}{0.484778in}%
\pgfsys@useobject{currentmarker}{}%
\end{pgfscope}%
\end{pgfscope}%
\begin{pgfscope}%
\pgftext[x=1.928622in,y=0.387555in,,top]{\rmfamily\fontsize{20.000000}{24.000000}\selectfont \(\displaystyle 0.2\)}%
\end{pgfscope}%
\begin{pgfscope}%
\pgfsetbuttcap%
\pgfsetroundjoin%
\definecolor{currentfill}{rgb}{0.000000,0.000000,0.000000}%
\pgfsetfillcolor{currentfill}%
\pgfsetlinewidth{0.803000pt}%
\definecolor{currentstroke}{rgb}{0.000000,0.000000,0.000000}%
\pgfsetstrokecolor{currentstroke}%
\pgfsetdash{}{0pt}%
\pgfsys@defobject{currentmarker}{\pgfqpoint{0.000000in}{-0.048611in}}{\pgfqpoint{0.000000in}{0.000000in}}{%
\pgfpathmoveto{\pgfqpoint{0.000000in}{0.000000in}}%
\pgfpathlineto{\pgfqpoint{0.000000in}{-0.048611in}}%
\pgfusepath{stroke,fill}%
}%
\begin{pgfscope}%
\pgfsys@transformshift{3.168622in}{0.484778in}%
\pgfsys@useobject{currentmarker}{}%
\end{pgfscope}%
\end{pgfscope}%
\begin{pgfscope}%
\pgftext[x=3.168622in,y=0.387555in,,top]{\rmfamily\fontsize{20.000000}{24.000000}\selectfont \(\displaystyle 0.4\)}%
\end{pgfscope}%
\begin{pgfscope}%
\pgfsetbuttcap%
\pgfsetroundjoin%
\definecolor{currentfill}{rgb}{0.000000,0.000000,0.000000}%
\pgfsetfillcolor{currentfill}%
\pgfsetlinewidth{0.803000pt}%
\definecolor{currentstroke}{rgb}{0.000000,0.000000,0.000000}%
\pgfsetstrokecolor{currentstroke}%
\pgfsetdash{}{0pt}%
\pgfsys@defobject{currentmarker}{\pgfqpoint{0.000000in}{-0.048611in}}{\pgfqpoint{0.000000in}{0.000000in}}{%
\pgfpathmoveto{\pgfqpoint{0.000000in}{0.000000in}}%
\pgfpathlineto{\pgfqpoint{0.000000in}{-0.048611in}}%
\pgfusepath{stroke,fill}%
}%
\begin{pgfscope}%
\pgfsys@transformshift{4.408622in}{0.484778in}%
\pgfsys@useobject{currentmarker}{}%
\end{pgfscope}%
\end{pgfscope}%
\begin{pgfscope}%
\pgftext[x=4.408622in,y=0.387555in,,top]{\rmfamily\fontsize{20.000000}{24.000000}\selectfont  }%
\end{pgfscope}%
\begin{pgfscope}%
\pgfsetbuttcap%
\pgfsetroundjoin%
\definecolor{currentfill}{rgb}{0.000000,0.000000,0.000000}%
\pgfsetfillcolor{currentfill}%
\pgfsetlinewidth{0.803000pt}%
\definecolor{currentstroke}{rgb}{0.000000,0.000000,0.000000}%
\pgfsetstrokecolor{currentstroke}%
\pgfsetdash{}{0pt}%
\pgfsys@defobject{currentmarker}{\pgfqpoint{0.000000in}{-0.048611in}}{\pgfqpoint{0.000000in}{0.000000in}}{%
\pgfpathmoveto{\pgfqpoint{0.000000in}{0.000000in}}%
\pgfpathlineto{\pgfqpoint{0.000000in}{-0.048611in}}%
\pgfusepath{stroke,fill}%
}%
\begin{pgfscope}%
\pgfsys@transformshift{5.648622in}{0.484778in}%
\pgfsys@useobject{currentmarker}{}%
\end{pgfscope}%
\end{pgfscope}%
\begin{pgfscope}%
\pgftext[x=5.648622in,y=0.387555in,,top]{\rmfamily\fontsize{20.000000}{24.000000}\selectfont \(\displaystyle 0.8\)}%
\end{pgfscope}%
\begin{pgfscope}%
\pgfsetbuttcap%
\pgfsetroundjoin%
\definecolor{currentfill}{rgb}{0.000000,0.000000,0.000000}%
\pgfsetfillcolor{currentfill}%
\pgfsetlinewidth{0.803000pt}%
\definecolor{currentstroke}{rgb}{0.000000,0.000000,0.000000}%
\pgfsetstrokecolor{currentstroke}%
\pgfsetdash{}{0pt}%
\pgfsys@defobject{currentmarker}{\pgfqpoint{0.000000in}{-0.048611in}}{\pgfqpoint{0.000000in}{0.000000in}}{%
\pgfpathmoveto{\pgfqpoint{0.000000in}{0.000000in}}%
\pgfpathlineto{\pgfqpoint{0.000000in}{-0.048611in}}%
\pgfusepath{stroke,fill}%
}%
\begin{pgfscope}%
\pgfsys@transformshift{6.888622in}{0.484778in}%
\pgfsys@useobject{currentmarker}{}%
\end{pgfscope}%
\end{pgfscope}%
\begin{pgfscope}%
\pgftext[x=6.888622in,y=0.387555in,,top]{\rmfamily\fontsize{20.000000}{24.000000}\selectfont \(\displaystyle 1\)}%
\end{pgfscope}%
\begin{pgfscope}%
\pgfsetbuttcap%
\pgfsetroundjoin%
\definecolor{currentfill}{rgb}{0.000000,0.000000,0.000000}%
\pgfsetfillcolor{currentfill}%
\pgfsetlinewidth{0.803000pt}%
\definecolor{currentstroke}{rgb}{0.000000,0.000000,0.000000}%
\pgfsetstrokecolor{currentstroke}%
\pgfsetdash{}{0pt}%
\pgfsys@defobject{currentmarker}{\pgfqpoint{0.000000in}{-0.048611in}}{\pgfqpoint{0.000000in}{0.000000in}}{%
\pgfpathmoveto{\pgfqpoint{0.000000in}{0.000000in}}%
\pgfpathlineto{\pgfqpoint{0.000000in}{-0.048611in}}%
\pgfusepath{stroke,fill}%
}%
\begin{pgfscope}%
\pgfsys@transformshift{4.964484in}{0.484778in}%
\pgfsys@useobject{currentmarker}{}%
\end{pgfscope}%
\end{pgfscope}%
\begin{pgfscope}%
\pgftext[x=4.964484in,y=0.387555in,,top]{\rmfamily\fontsize{20.000000}{24.000000}\selectfont \(\displaystyle \widehat{c}_{\beta}\)}%
\end{pgfscope}%
\begin{pgfscope}%
\pgfsetbuttcap%
\pgfsetroundjoin%
\definecolor{currentfill}{rgb}{0.000000,0.000000,0.000000}%
\pgfsetfillcolor{currentfill}%
\pgfsetlinewidth{0.803000pt}%
\definecolor{currentstroke}{rgb}{0.000000,0.000000,0.000000}%
\pgfsetstrokecolor{currentstroke}%
\pgfsetdash{}{0pt}%
\pgfsys@defobject{currentmarker}{\pgfqpoint{0.000000in}{-0.048611in}}{\pgfqpoint{0.000000in}{0.000000in}}{%
\pgfpathmoveto{\pgfqpoint{0.000000in}{0.000000in}}%
\pgfpathlineto{\pgfqpoint{0.000000in}{-0.048611in}}%
\pgfusepath{stroke,fill}%
}%
\begin{pgfscope}%
\pgfsys@transformshift{3.885563in}{0.484778in}%
\pgfsys@useobject{currentmarker}{}%
\end{pgfscope}%
\end{pgfscope}%
\begin{pgfscope}%
\pgftext[x=3.885563in,y=0.387555in,,top]{\rmfamily\fontsize{20.000000}{24.000000}\selectfont \(\displaystyle c(\alpha,\beta)\)}%
\end{pgfscope}%
\begin{pgfscope}%
\pgfsetbuttcap%
\pgfsetroundjoin%
\definecolor{currentfill}{rgb}{0.000000,0.000000,0.000000}%
\pgfsetfillcolor{currentfill}%
\pgfsetlinewidth{0.803000pt}%
\definecolor{currentstroke}{rgb}{0.000000,0.000000,0.000000}%
\pgfsetstrokecolor{currentstroke}%
\pgfsetdash{}{0pt}%
\pgfsys@defobject{currentmarker}{\pgfqpoint{-0.048611in}{0.000000in}}{\pgfqpoint{0.000000in}{0.000000in}}{%
\pgfpathmoveto{\pgfqpoint{0.000000in}{0.000000in}}%
\pgfpathlineto{\pgfqpoint{-0.048611in}{0.000000in}}%
\pgfusepath{stroke,fill}%
}%
\begin{pgfscope}%
\pgfsys@transformshift{0.688622in}{0.484778in}%
\pgfsys@useobject{currentmarker}{}%
\end{pgfscope}%
\end{pgfscope}%
\begin{pgfscope}%
\pgftext[x=0.459293in,y=0.388389in,left,base]{\rmfamily\fontsize{20.000000}{24.000000}\selectfont \(\displaystyle 0\)}%
\end{pgfscope}%
\begin{pgfscope}%
\pgfsetbuttcap%
\pgfsetroundjoin%
\definecolor{currentfill}{rgb}{0.000000,0.000000,0.000000}%
\pgfsetfillcolor{currentfill}%
\pgfsetlinewidth{0.803000pt}%
\definecolor{currentstroke}{rgb}{0.000000,0.000000,0.000000}%
\pgfsetstrokecolor{currentstroke}%
\pgfsetdash{}{0pt}%
\pgfsys@defobject{currentmarker}{\pgfqpoint{-0.048611in}{0.000000in}}{\pgfqpoint{0.000000in}{0.000000in}}{%
\pgfpathmoveto{\pgfqpoint{0.000000in}{0.000000in}}%
\pgfpathlineto{\pgfqpoint{-0.048611in}{0.000000in}}%
\pgfusepath{stroke,fill}%
}%
\begin{pgfscope}%
\pgfsys@transformshift{0.688622in}{1.390778in}%
\pgfsys@useobject{currentmarker}{}%
\end{pgfscope}%
\end{pgfscope}%
\begin{pgfscope}%
\pgftext[x=0.248838in,y=1.294389in,left,base]{\rmfamily\fontsize{20.000000}{24.000000}\selectfont \(\displaystyle 0.2\)}%
\end{pgfscope}%
\begin{pgfscope}%
\pgfsetbuttcap%
\pgfsetroundjoin%
\definecolor{currentfill}{rgb}{0.000000,0.000000,0.000000}%
\pgfsetfillcolor{currentfill}%
\pgfsetlinewidth{0.803000pt}%
\definecolor{currentstroke}{rgb}{0.000000,0.000000,0.000000}%
\pgfsetstrokecolor{currentstroke}%
\pgfsetdash{}{0pt}%
\pgfsys@defobject{currentmarker}{\pgfqpoint{-0.048611in}{0.000000in}}{\pgfqpoint{0.000000in}{0.000000in}}{%
\pgfpathmoveto{\pgfqpoint{0.000000in}{0.000000in}}%
\pgfpathlineto{\pgfqpoint{-0.048611in}{0.000000in}}%
\pgfusepath{stroke,fill}%
}%
\begin{pgfscope}%
\pgfsys@transformshift{0.688622in}{2.296778in}%
\pgfsys@useobject{currentmarker}{}%
\end{pgfscope}%
\end{pgfscope}%
\begin{pgfscope}%
\pgftext[x=0.591400in,y=2.200389in,left,base]{\rmfamily\fontsize{20.000000}{24.000000}\selectfont  }%
\end{pgfscope}%
\begin{pgfscope}%
\pgfsetbuttcap%
\pgfsetroundjoin%
\definecolor{currentfill}{rgb}{0.000000,0.000000,0.000000}%
\pgfsetfillcolor{currentfill}%
\pgfsetlinewidth{0.803000pt}%
\definecolor{currentstroke}{rgb}{0.000000,0.000000,0.000000}%
\pgfsetstrokecolor{currentstroke}%
\pgfsetdash{}{0pt}%
\pgfsys@defobject{currentmarker}{\pgfqpoint{-0.048611in}{0.000000in}}{\pgfqpoint{0.000000in}{0.000000in}}{%
\pgfpathmoveto{\pgfqpoint{0.000000in}{0.000000in}}%
\pgfpathlineto{\pgfqpoint{-0.048611in}{0.000000in}}%
\pgfusepath{stroke,fill}%
}%
\begin{pgfscope}%
\pgfsys@transformshift{0.688622in}{3.202777in}%
\pgfsys@useobject{currentmarker}{}%
\end{pgfscope}%
\end{pgfscope}%
\begin{pgfscope}%
\pgftext[x=0.248838in,y=3.106389in,left,base]{\rmfamily\fontsize{20.000000}{24.000000}\selectfont \(\displaystyle 0.6\)}%
\end{pgfscope}%
\begin{pgfscope}%
\pgfsetbuttcap%
\pgfsetroundjoin%
\definecolor{currentfill}{rgb}{0.000000,0.000000,0.000000}%
\pgfsetfillcolor{currentfill}%
\pgfsetlinewidth{0.803000pt}%
\definecolor{currentstroke}{rgb}{0.000000,0.000000,0.000000}%
\pgfsetstrokecolor{currentstroke}%
\pgfsetdash{}{0pt}%
\pgfsys@defobject{currentmarker}{\pgfqpoint{-0.048611in}{0.000000in}}{\pgfqpoint{0.000000in}{0.000000in}}{%
\pgfpathmoveto{\pgfqpoint{0.000000in}{0.000000in}}%
\pgfpathlineto{\pgfqpoint{-0.048611in}{0.000000in}}%
\pgfusepath{stroke,fill}%
}%
\begin{pgfscope}%
\pgfsys@transformshift{0.688622in}{4.108778in}%
\pgfsys@useobject{currentmarker}{}%
\end{pgfscope}%
\end{pgfscope}%
\begin{pgfscope}%
\pgftext[x=0.248838in,y=4.012389in,left,base]{\rmfamily\fontsize{20.000000}{24.000000}\selectfont \(\displaystyle 0.8\)}%
\end{pgfscope}%
\begin{pgfscope}%
\pgfsetbuttcap%
\pgfsetroundjoin%
\definecolor{currentfill}{rgb}{0.000000,0.000000,0.000000}%
\pgfsetfillcolor{currentfill}%
\pgfsetlinewidth{0.803000pt}%
\definecolor{currentstroke}{rgb}{0.000000,0.000000,0.000000}%
\pgfsetstrokecolor{currentstroke}%
\pgfsetdash{}{0pt}%
\pgfsys@defobject{currentmarker}{\pgfqpoint{-0.048611in}{0.000000in}}{\pgfqpoint{0.000000in}{0.000000in}}{%
\pgfpathmoveto{\pgfqpoint{0.000000in}{0.000000in}}%
\pgfpathlineto{\pgfqpoint{-0.048611in}{0.000000in}}%
\pgfusepath{stroke,fill}%
}%
\begin{pgfscope}%
\pgfsys@transformshift{0.688622in}{5.014778in}%
\pgfsys@useobject{currentmarker}{}%
\end{pgfscope}%
\end{pgfscope}%
\begin{pgfscope}%
\pgftext[x=0.459293in,y=4.918389in,left,base]{\rmfamily\fontsize{20.000000}{24.000000}\selectfont \(\displaystyle 1\)}%
\end{pgfscope}%
\begin{pgfscope}%
\pgfsetbuttcap%
\pgfsetroundjoin%
\definecolor{currentfill}{rgb}{0.000000,0.000000,0.000000}%
\pgfsetfillcolor{currentfill}%
\pgfsetlinewidth{0.803000pt}%
\definecolor{currentstroke}{rgb}{0.000000,0.000000,0.000000}%
\pgfsetstrokecolor{currentstroke}%
\pgfsetdash{}{0pt}%
\pgfsys@defobject{currentmarker}{\pgfqpoint{-0.048611in}{0.000000in}}{\pgfqpoint{0.000000in}{0.000000in}}{%
\pgfpathmoveto{\pgfqpoint{0.000000in}{0.000000in}}%
\pgfpathlineto{\pgfqpoint{-0.048611in}{0.000000in}}%
\pgfusepath{stroke,fill}%
}%
\begin{pgfscope}%
\pgfsys@transformshift{0.688622in}{2.451832in}%
\pgfsys@useobject{currentmarker}{}%
\end{pgfscope}%
\end{pgfscope}%
\begin{pgfscope}%
\pgftext[x=0.100000in,y=2.327919in,left,base]{\rmfamily\fontsize{20.000000}{24.000000}\selectfont \(\displaystyle F_{\alpha,\beta}^1\)}%
\end{pgfscope}%
\begin{pgfscope}%
\pgfsetbuttcap%
\pgfsetroundjoin%
\definecolor{currentfill}{rgb}{0.000000,0.000000,0.000000}%
\pgfsetfillcolor{currentfill}%
\pgfsetlinewidth{0.803000pt}%
\definecolor{currentstroke}{rgb}{0.000000,0.000000,0.000000}%
\pgfsetstrokecolor{currentstroke}%
\pgfsetdash{}{0pt}%
\pgfsys@defobject{currentmarker}{\pgfqpoint{-0.048611in}{0.000000in}}{\pgfqpoint{0.000000in}{0.000000in}}{%
\pgfpathmoveto{\pgfqpoint{0.000000in}{0.000000in}}%
\pgfpathlineto{\pgfqpoint{-0.048611in}{0.000000in}}%
\pgfusepath{stroke,fill}%
}%
\begin{pgfscope}%
\pgfsys@transformshift{0.688622in}{2.083125in}%
\pgfsys@useobject{currentmarker}{}%
\end{pgfscope}%
\end{pgfscope}%
\begin{pgfscope}%
\pgftext[x=0.100000in,y=1.959212in,left,base]{\rmfamily\fontsize{20.000000}{24.000000}\selectfont \(\displaystyle F_{\alpha,\beta}^0\)}%
\end{pgfscope}%
\begin{pgfscope}%
\pgfpathrectangle{\pgfqpoint{0.688622in}{0.484778in}}{\pgfqpoint{6.200000in}{4.530000in}} %
\pgfusepath{clip}%
\pgfsetbuttcap%
\pgfsetroundjoin%
\pgfsetlinewidth{1.505625pt}%
\definecolor{currentstroke}{rgb}{0.501961,0.501961,0.501961}%
\pgfsetstrokecolor{currentstroke}%
\pgfsetdash{{5.550000pt}{2.400000pt}}{0.000000pt}%
\pgfpathmoveto{\pgfqpoint{4.964484in}{0.484778in}}%
\pgfpathlineto{\pgfqpoint{4.964484in}{5.028666in}}%
\pgfusepath{stroke}%
\end{pgfscope}%
\begin{pgfscope}%
\pgfpathrectangle{\pgfqpoint{0.688622in}{0.484778in}}{\pgfqpoint{6.200000in}{4.530000in}} %
\pgfusepath{clip}%
\pgfsetbuttcap%
\pgfsetroundjoin%
\pgfsetlinewidth{1.003750pt}%
\definecolor{currentstroke}{rgb}{0.000000,0.392157,0.000000}%
\pgfsetstrokecolor{currentstroke}%
\pgfsetdash{{3.700000pt}{1.600000pt}}{0.000000pt}%
\pgfpathmoveto{\pgfqpoint{0.688622in}{2.451832in}}%
\pgfpathlineto{\pgfqpoint{6.888622in}{2.451832in}}%
\pgfusepath{stroke}%
\end{pgfscope}%
\begin{pgfscope}%
\pgfpathrectangle{\pgfqpoint{0.688622in}{0.484778in}}{\pgfqpoint{6.200000in}{4.530000in}} %
\pgfusepath{clip}%
\pgfsetbuttcap%
\pgfsetroundjoin%
\pgfsetlinewidth{1.505625pt}%
\definecolor{currentstroke}{rgb}{0.501961,0.501961,0.501961}%
\pgfsetstrokecolor{currentstroke}%
\pgfsetdash{{5.550000pt}{2.400000pt}}{0.000000pt}%
\pgfpathmoveto{\pgfqpoint{3.885563in}{0.484778in}}%
\pgfpathlineto{\pgfqpoint{3.885563in}{2.451832in}}%
\pgfusepath{stroke}%
\end{pgfscope}%
\begin{pgfscope}%
\pgfpathrectangle{\pgfqpoint{0.688622in}{0.484778in}}{\pgfqpoint{6.200000in}{4.530000in}} %
\pgfusepath{clip}%
\pgfsetrectcap%
\pgfsetroundjoin%
\pgfsetlinewidth{2.007500pt}%
\definecolor{currentstroke}{rgb}{0.721569,0.525490,0.043137}%
\pgfsetstrokecolor{currentstroke}%
\pgfsetdash{}{0pt}%
\pgfpathmoveto{\pgfqpoint{0.688622in}{2.083125in}}%
\pgfpathlineto{\pgfqpoint{0.893222in}{2.084635in}}%
\pgfpathlineto{\pgfqpoint{1.097822in}{2.089165in}}%
\pgfpathlineto{\pgfqpoint{1.302422in}{2.096716in}}%
\pgfpathlineto{\pgfqpoint{1.507022in}{2.107287in}}%
\pgfpathlineto{\pgfqpoint{1.711622in}{2.120879in}}%
\pgfpathlineto{\pgfqpoint{1.916222in}{2.137490in}}%
\pgfpathlineto{\pgfqpoint{2.120822in}{2.157123in}}%
\pgfpathlineto{\pgfqpoint{2.325422in}{2.179775in}}%
\pgfpathlineto{\pgfqpoint{2.530022in}{2.205448in}}%
\pgfpathlineto{\pgfqpoint{2.740822in}{2.235057in}}%
\pgfpathlineto{\pgfqpoint{2.951622in}{2.267873in}}%
\pgfpathlineto{\pgfqpoint{3.162422in}{2.303895in}}%
\pgfpathlineto{\pgfqpoint{3.373222in}{2.343123in}}%
\pgfpathlineto{\pgfqpoint{3.584022in}{2.385558in}}%
\pgfpathlineto{\pgfqpoint{3.794822in}{2.431198in}}%
\pgfpathlineto{\pgfqpoint{4.005622in}{2.480045in}}%
\pgfpathlineto{\pgfqpoint{4.216422in}{2.532097in}}%
\pgfpathlineto{\pgfqpoint{4.427222in}{2.587356in}}%
\pgfpathlineto{\pgfqpoint{4.638022in}{2.645821in}}%
\pgfpathlineto{\pgfqpoint{4.848822in}{2.707493in}}%
\pgfpathlineto{\pgfqpoint{4.966622in}{2.743352in}}%
\pgfpathlineto{\pgfqpoint{4.964484in}{2.742693in}}%
\pgfusepath{stroke}%
\end{pgfscope}%
\begin{pgfscope}%
\pgfpathrectangle{\pgfqpoint{0.688622in}{0.484778in}}{\pgfqpoint{6.200000in}{4.530000in}} %
\pgfusepath{clip}%
\pgfsetbuttcap%
\pgfsetroundjoin%
\pgfsetlinewidth{3.011250pt}%
\definecolor{currentstroke}{rgb}{0.000000,0.000000,1.000000}%
\pgfsetstrokecolor{currentstroke}%
\pgfsetdash{{19.200000pt}{4.800000pt}{3.000000pt}{4.800000pt}}{0.000000pt}%
\pgfpathmoveto{\pgfqpoint{4.964484in}{2.742693in}}%
\pgfpathlineto{\pgfqpoint{4.972822in}{2.790588in}}%
\pgfpathlineto{\pgfqpoint{4.985222in}{2.821316in}}%
\pgfpathlineto{\pgfqpoint{5.003822in}{2.855608in}}%
\pgfpathlineto{\pgfqpoint{5.028622in}{2.893188in}}%
\pgfpathlineto{\pgfqpoint{5.059622in}{2.934165in}}%
\pgfpathlineto{\pgfqpoint{5.103022in}{2.985875in}}%
\pgfpathlineto{\pgfqpoint{5.158822in}{3.047283in}}%
\pgfpathlineto{\pgfqpoint{5.239422in}{3.130960in}}%
\pgfpathlineto{\pgfqpoint{5.369622in}{3.261036in}}%
\pgfpathlineto{\pgfqpoint{5.661022in}{3.551380in}}%
\pgfpathlineto{\pgfqpoint{5.809822in}{3.704651in}}%
\pgfpathlineto{\pgfqpoint{5.946222in}{3.849732in}}%
\pgfpathlineto{\pgfqpoint{6.070222in}{3.986012in}}%
\pgfpathlineto{\pgfqpoint{6.194222in}{4.126884in}}%
\pgfpathlineto{\pgfqpoint{6.312022in}{4.265272in}}%
\pgfpathlineto{\pgfqpoint{6.429822in}{4.408366in}}%
\pgfpathlineto{\pgfqpoint{6.547622in}{4.556401in}}%
\pgfpathlineto{\pgfqpoint{6.665422in}{4.709600in}}%
\pgfpathlineto{\pgfqpoint{6.777022in}{4.859690in}}%
\pgfpathlineto{\pgfqpoint{6.888622in}{5.014778in}}%
\pgfpathlineto{\pgfqpoint{6.888622in}{5.014778in}}%
\pgfusepath{stroke}%
\end{pgfscope}%
\begin{pgfscope}%
\pgfpathrectangle{\pgfqpoint{0.688622in}{0.484778in}}{\pgfqpoint{6.200000in}{4.530000in}} %
\pgfusepath{clip}%
\pgfsetbuttcap%
\pgfsetroundjoin%
\pgfsetlinewidth{2.509375pt}%
\definecolor{currentstroke}{rgb}{1.000000,0.000000,0.000000}%
\pgfsetstrokecolor{currentstroke}%
\pgfsetdash{{9.250000pt}{4.000000pt}}{0.000000pt}%
\pgfpathmoveto{\pgfqpoint{4.964484in}{2.742693in}}%
\pgfpathlineto{\pgfqpoint{4.972822in}{2.700837in}}%
\pgfpathlineto{\pgfqpoint{4.985222in}{2.679128in}}%
\pgfpathlineto{\pgfqpoint{4.997622in}{2.664474in}}%
\pgfpathlineto{\pgfqpoint{5.016222in}{2.648014in}}%
\pgfpathlineto{\pgfqpoint{5.041022in}{2.631385in}}%
\pgfpathlineto{\pgfqpoint{5.072022in}{2.615287in}}%
\pgfpathlineto{\pgfqpoint{5.109222in}{2.600008in}}%
\pgfpathlineto{\pgfqpoint{5.152622in}{2.585657in}}%
\pgfpathlineto{\pgfqpoint{5.208422in}{2.570753in}}%
\pgfpathlineto{\pgfqpoint{5.276622in}{2.556099in}}%
\pgfpathlineto{\pgfqpoint{5.357222in}{2.542139in}}%
\pgfpathlineto{\pgfqpoint{5.456422in}{2.528315in}}%
\pgfpathlineto{\pgfqpoint{5.574222in}{2.515161in}}%
\pgfpathlineto{\pgfqpoint{5.716822in}{2.502458in}}%
\pgfpathlineto{\pgfqpoint{5.890422in}{2.490230in}}%
\pgfpathlineto{\pgfqpoint{6.101222in}{2.478611in}}%
\pgfpathlineto{\pgfqpoint{6.361622in}{2.467529in}}%
\pgfpathlineto{\pgfqpoint{6.677822in}{2.457291in}}%
\pgfpathlineto{\pgfqpoint{6.888622in}{2.451832in}}%
\pgfpathlineto{\pgfqpoint{6.888622in}{2.451832in}}%
\pgfusepath{stroke}%
\end{pgfscope}%
\begin{pgfscope}%
\pgfpathrectangle{\pgfqpoint{0.688622in}{0.484778in}}{\pgfqpoint{6.200000in}{4.530000in}} %
\pgfusepath{clip}%
\pgfsetbuttcap%
\pgfsetroundjoin%
\definecolor{currentfill}{rgb}{0.721569,0.525490,0.043137}%
\pgfsetfillcolor{currentfill}%
\pgfsetlinewidth{1.003750pt}%
\definecolor{currentstroke}{rgb}{0.721569,0.525490,0.043137}%
\pgfsetstrokecolor{currentstroke}%
\pgfsetdash{}{0pt}%
\pgfsys@defobject{currentmarker}{\pgfqpoint{-0.027778in}{-0.027778in}}{\pgfqpoint{0.027778in}{0.027778in}}{%
\pgfpathmoveto{\pgfqpoint{0.000000in}{-0.027778in}}%
\pgfpathcurveto{\pgfqpoint{0.007367in}{-0.027778in}}{\pgfqpoint{0.014433in}{-0.024851in}}{\pgfqpoint{0.019642in}{-0.019642in}}%
\pgfpathcurveto{\pgfqpoint{0.024851in}{-0.014433in}}{\pgfqpoint{0.027778in}{-0.007367in}}{\pgfqpoint{0.027778in}{0.000000in}}%
\pgfpathcurveto{\pgfqpoint{0.027778in}{0.007367in}}{\pgfqpoint{0.024851in}{0.014433in}}{\pgfqpoint{0.019642in}{0.019642in}}%
\pgfpathcurveto{\pgfqpoint{0.014433in}{0.024851in}}{\pgfqpoint{0.007367in}{0.027778in}}{\pgfqpoint{0.000000in}{0.027778in}}%
\pgfpathcurveto{\pgfqpoint{-0.007367in}{0.027778in}}{\pgfqpoint{-0.014433in}{0.024851in}}{\pgfqpoint{-0.019642in}{0.019642in}}%
\pgfpathcurveto{\pgfqpoint{-0.024851in}{0.014433in}}{\pgfqpoint{-0.027778in}{0.007367in}}{\pgfqpoint{-0.027778in}{0.000000in}}%
\pgfpathcurveto{\pgfqpoint{-0.027778in}{-0.007367in}}{\pgfqpoint{-0.024851in}{-0.014433in}}{\pgfqpoint{-0.019642in}{-0.019642in}}%
\pgfpathcurveto{\pgfqpoint{-0.014433in}{-0.024851in}}{\pgfqpoint{-0.007367in}{-0.027778in}}{\pgfqpoint{0.000000in}{-0.027778in}}%
\pgfpathclose%
\pgfusepath{stroke,fill}%
}%
\begin{pgfscope}%
\pgfsys@transformshift{0.694822in}{2.083125in}%
\pgfsys@useobject{currentmarker}{}%
\end{pgfscope}%
\end{pgfscope}%
\begin{pgfscope}%
\pgfpathrectangle{\pgfqpoint{0.688622in}{0.484778in}}{\pgfqpoint{6.200000in}{4.530000in}} %
\pgfusepath{clip}%
\pgfsetbuttcap%
\pgfsetroundjoin%
\definecolor{currentfill}{rgb}{0.000000,0.000000,1.000000}%
\pgfsetfillcolor{currentfill}%
\pgfsetlinewidth{1.003750pt}%
\definecolor{currentstroke}{rgb}{0.000000,0.000000,1.000000}%
\pgfsetstrokecolor{currentstroke}%
\pgfsetdash{}{0pt}%
\pgfsys@defobject{currentmarker}{\pgfqpoint{-0.027778in}{-0.027778in}}{\pgfqpoint{0.027778in}{0.027778in}}{%
\pgfpathmoveto{\pgfqpoint{0.000000in}{-0.027778in}}%
\pgfpathcurveto{\pgfqpoint{0.007367in}{-0.027778in}}{\pgfqpoint{0.014433in}{-0.024851in}}{\pgfqpoint{0.019642in}{-0.019642in}}%
\pgfpathcurveto{\pgfqpoint{0.024851in}{-0.014433in}}{\pgfqpoint{0.027778in}{-0.007367in}}{\pgfqpoint{0.027778in}{0.000000in}}%
\pgfpathcurveto{\pgfqpoint{0.027778in}{0.007367in}}{\pgfqpoint{0.024851in}{0.014433in}}{\pgfqpoint{0.019642in}{0.019642in}}%
\pgfpathcurveto{\pgfqpoint{0.014433in}{0.024851in}}{\pgfqpoint{0.007367in}{0.027778in}}{\pgfqpoint{0.000000in}{0.027778in}}%
\pgfpathcurveto{\pgfqpoint{-0.007367in}{0.027778in}}{\pgfqpoint{-0.014433in}{0.024851in}}{\pgfqpoint{-0.019642in}{0.019642in}}%
\pgfpathcurveto{\pgfqpoint{-0.024851in}{0.014433in}}{\pgfqpoint{-0.027778in}{0.007367in}}{\pgfqpoint{-0.027778in}{0.000000in}}%
\pgfpathcurveto{\pgfqpoint{-0.027778in}{-0.007367in}}{\pgfqpoint{-0.024851in}{-0.014433in}}{\pgfqpoint{-0.019642in}{-0.019642in}}%
\pgfpathcurveto{\pgfqpoint{-0.014433in}{-0.024851in}}{\pgfqpoint{-0.007367in}{-0.027778in}}{\pgfqpoint{0.000000in}{-0.027778in}}%
\pgfpathclose%
\pgfusepath{stroke,fill}%
}%
\begin{pgfscope}%
\pgfsys@transformshift{6.882422in}{5.014778in}%
\pgfsys@useobject{currentmarker}{}%
\end{pgfscope}%
\end{pgfscope}%
\begin{pgfscope}%
\pgfpathrectangle{\pgfqpoint{0.688622in}{0.484778in}}{\pgfqpoint{6.200000in}{4.530000in}} %
\pgfusepath{clip}%
\pgfsetbuttcap%
\pgfsetroundjoin%
\definecolor{currentfill}{rgb}{1.000000,0.000000,0.000000}%
\pgfsetfillcolor{currentfill}%
\pgfsetlinewidth{1.003750pt}%
\definecolor{currentstroke}{rgb}{1.000000,0.000000,0.000000}%
\pgfsetstrokecolor{currentstroke}%
\pgfsetdash{}{0pt}%
\pgfsys@defobject{currentmarker}{\pgfqpoint{-0.027778in}{-0.027778in}}{\pgfqpoint{0.027778in}{0.027778in}}{%
\pgfpathmoveto{\pgfqpoint{0.000000in}{-0.027778in}}%
\pgfpathcurveto{\pgfqpoint{0.007367in}{-0.027778in}}{\pgfqpoint{0.014433in}{-0.024851in}}{\pgfqpoint{0.019642in}{-0.019642in}}%
\pgfpathcurveto{\pgfqpoint{0.024851in}{-0.014433in}}{\pgfqpoint{0.027778in}{-0.007367in}}{\pgfqpoint{0.027778in}{0.000000in}}%
\pgfpathcurveto{\pgfqpoint{0.027778in}{0.007367in}}{\pgfqpoint{0.024851in}{0.014433in}}{\pgfqpoint{0.019642in}{0.019642in}}%
\pgfpathcurveto{\pgfqpoint{0.014433in}{0.024851in}}{\pgfqpoint{0.007367in}{0.027778in}}{\pgfqpoint{0.000000in}{0.027778in}}%
\pgfpathcurveto{\pgfqpoint{-0.007367in}{0.027778in}}{\pgfqpoint{-0.014433in}{0.024851in}}{\pgfqpoint{-0.019642in}{0.019642in}}%
\pgfpathcurveto{\pgfqpoint{-0.024851in}{0.014433in}}{\pgfqpoint{-0.027778in}{0.007367in}}{\pgfqpoint{-0.027778in}{0.000000in}}%
\pgfpathcurveto{\pgfqpoint{-0.027778in}{-0.007367in}}{\pgfqpoint{-0.024851in}{-0.014433in}}{\pgfqpoint{-0.019642in}{-0.019642in}}%
\pgfpathcurveto{\pgfqpoint{-0.014433in}{-0.024851in}}{\pgfqpoint{-0.007367in}{-0.027778in}}{\pgfqpoint{0.000000in}{-0.027778in}}%
\pgfpathclose%
\pgfusepath{stroke,fill}%
}%
\begin{pgfscope}%
\pgfsys@transformshift{6.882422in}{2.451832in}%
\pgfsys@useobject{currentmarker}{}%
\end{pgfscope}%
\end{pgfscope}%
\begin{pgfscope}%
\pgfpathrectangle{\pgfqpoint{0.688622in}{0.484778in}}{\pgfqpoint{6.200000in}{4.530000in}} %
\pgfusepath{clip}%
\pgfsetbuttcap%
\pgfsetroundjoin%
\definecolor{currentfill}{rgb}{0.721569,0.525490,0.043137}%
\pgfsetfillcolor{currentfill}%
\pgfsetlinewidth{1.003750pt}%
\definecolor{currentstroke}{rgb}{0.721569,0.525490,0.043137}%
\pgfsetstrokecolor{currentstroke}%
\pgfsetdash{}{0pt}%
\pgfsys@defobject{currentmarker}{\pgfqpoint{-0.027778in}{-0.027778in}}{\pgfqpoint{0.027778in}{0.027778in}}{%
\pgfpathmoveto{\pgfqpoint{0.000000in}{-0.027778in}}%
\pgfpathcurveto{\pgfqpoint{0.007367in}{-0.027778in}}{\pgfqpoint{0.014433in}{-0.024851in}}{\pgfqpoint{0.019642in}{-0.019642in}}%
\pgfpathcurveto{\pgfqpoint{0.024851in}{-0.014433in}}{\pgfqpoint{0.027778in}{-0.007367in}}{\pgfqpoint{0.027778in}{0.000000in}}%
\pgfpathcurveto{\pgfqpoint{0.027778in}{0.007367in}}{\pgfqpoint{0.024851in}{0.014433in}}{\pgfqpoint{0.019642in}{0.019642in}}%
\pgfpathcurveto{\pgfqpoint{0.014433in}{0.024851in}}{\pgfqpoint{0.007367in}{0.027778in}}{\pgfqpoint{0.000000in}{0.027778in}}%
\pgfpathcurveto{\pgfqpoint{-0.007367in}{0.027778in}}{\pgfqpoint{-0.014433in}{0.024851in}}{\pgfqpoint{-0.019642in}{0.019642in}}%
\pgfpathcurveto{\pgfqpoint{-0.024851in}{0.014433in}}{\pgfqpoint{-0.027778in}{0.007367in}}{\pgfqpoint{-0.027778in}{0.000000in}}%
\pgfpathcurveto{\pgfqpoint{-0.027778in}{-0.007367in}}{\pgfqpoint{-0.024851in}{-0.014433in}}{\pgfqpoint{-0.019642in}{-0.019642in}}%
\pgfpathcurveto{\pgfqpoint{-0.014433in}{-0.024851in}}{\pgfqpoint{-0.007367in}{-0.027778in}}{\pgfqpoint{0.000000in}{-0.027778in}}%
\pgfpathclose%
\pgfusepath{stroke,fill}%
}%
\begin{pgfscope}%
\pgfsys@transformshift{3.885563in}{2.451832in}%
\pgfsys@useobject{currentmarker}{}%
\end{pgfscope}%
\end{pgfscope}%
\begin{pgfscope}%
\pgfsetrectcap%
\pgfsetmiterjoin%
\pgfsetlinewidth{0.803000pt}%
\definecolor{currentstroke}{rgb}{0.000000,0.000000,0.000000}%
\pgfsetstrokecolor{currentstroke}%
\pgfsetdash{}{0pt}%
\pgfpathmoveto{\pgfqpoint{0.688622in}{0.484778in}}%
\pgfpathlineto{\pgfqpoint{0.688622in}{5.014778in}}%
\pgfusepath{stroke}%
\end{pgfscope}%
\begin{pgfscope}%
\pgfsetrectcap%
\pgfsetmiterjoin%
\pgfsetlinewidth{0.803000pt}%
\definecolor{currentstroke}{rgb}{0.000000,0.000000,0.000000}%
\pgfsetstrokecolor{currentstroke}%
\pgfsetdash{}{0pt}%
\pgfpathmoveto{\pgfqpoint{6.888622in}{0.484778in}}%
\pgfpathlineto{\pgfqpoint{6.888622in}{5.014778in}}%
\pgfusepath{stroke}%
\end{pgfscope}%
\begin{pgfscope}%
\pgfsetrectcap%
\pgfsetmiterjoin%
\pgfsetlinewidth{0.803000pt}%
\definecolor{currentstroke}{rgb}{0.000000,0.000000,0.000000}%
\pgfsetstrokecolor{currentstroke}%
\pgfsetdash{}{0pt}%
\pgfpathmoveto{\pgfqpoint{0.688622in}{0.484778in}}%
\pgfpathlineto{\pgfqpoint{6.888622in}{0.484778in}}%
\pgfusepath{stroke}%
\end{pgfscope}%
\begin{pgfscope}%
\pgfsetrectcap%
\pgfsetmiterjoin%
\pgfsetlinewidth{0.803000pt}%
\definecolor{currentstroke}{rgb}{0.000000,0.000000,0.000000}%
\pgfsetstrokecolor{currentstroke}%
\pgfsetdash{}{0pt}%
\pgfpathmoveto{\pgfqpoint{0.688622in}{5.014778in}}%
\pgfpathlineto{\pgfqpoint{6.888622in}{5.014778in}}%
\pgfusepath{stroke}%
\end{pgfscope}%
\begin{pgfscope}%
\definecolor{textcolor}{rgb}{0.000000,0.392157,0.000000}%
\pgfsetstrokecolor{textcolor}%
\pgfsetfillcolor{textcolor}%
\pgftext[x=1.308622in,y=2.542432in,left,base]{\color{textcolor}\rmfamily\fontsize{20.000000}{24.000000}\selectfont \(\displaystyle y=f_{\alpha,\beta,1}(1)\)}%
\end{pgfscope}%
\begin{pgfscope}%
\pgfsetbuttcap%
\pgfsetmiterjoin%
\definecolor{currentfill}{rgb}{1.000000,1.000000,1.000000}%
\pgfsetfillcolor{currentfill}%
\pgfsetlinewidth{1.003750pt}%
\definecolor{currentstroke}{rgb}{0.000000,0.000000,0.000000}%
\pgfsetstrokecolor{currentstroke}%
\pgfsetdash{}{0pt}%
\pgfpathmoveto{\pgfqpoint{2.855155in}{5.166235in}}%
\pgfpathlineto{\pgfqpoint{4.722089in}{5.166235in}}%
\pgfpathlineto{\pgfqpoint{4.722089in}{5.490679in}}%
\pgfpathlineto{\pgfqpoint{2.855155in}{5.490679in}}%
\pgfpathclose%
\pgfusepath{stroke,fill}%
\end{pgfscope}%
\begin{pgfscope}%
\pgftext[x=3.788622in,y=5.268458in,,base]{\rmfamily\fontsize{16.000000}{19.200000}\selectfont \(\displaystyle \alpha=\) 0.35, \(\displaystyle \beta=\) 0.58}%
\end{pgfscope}%
\begin{pgfscope}%
\pgfsetbuttcap%
\pgfsetmiterjoin%
\definecolor{currentfill}{rgb}{1.000000,1.000000,1.000000}%
\pgfsetfillcolor{currentfill}%
\pgfsetfillopacity{0.800000}%
\pgfsetlinewidth{1.003750pt}%
\definecolor{currentstroke}{rgb}{0.800000,0.800000,0.800000}%
\pgfsetstrokecolor{currentstroke}%
\pgfsetstrokeopacity{0.800000}%
\pgfsetdash{}{0pt}%
\pgfpathmoveto{\pgfqpoint{0.868178in}{3.429277in}}%
\pgfpathlineto{\pgfqpoint{3.402511in}{3.429277in}}%
\pgfpathquadraticcurveto{\pgfqpoint{3.458067in}{3.429277in}}{\pgfqpoint{3.458067in}{3.484833in}}%
\pgfpathlineto{\pgfqpoint{3.458067in}{4.767360in}}%
\pgfpathquadraticcurveto{\pgfqpoint{3.458067in}{4.822915in}}{\pgfqpoint{3.402511in}{4.822915in}}%
\pgfpathlineto{\pgfqpoint{0.868178in}{4.822915in}}%
\pgfpathquadraticcurveto{\pgfqpoint{0.812622in}{4.822915in}}{\pgfqpoint{0.812622in}{4.767360in}}%
\pgfpathlineto{\pgfqpoint{0.812622in}{3.484833in}}%
\pgfpathquadraticcurveto{\pgfqpoint{0.812622in}{3.429277in}}{\pgfqpoint{0.868178in}{3.429277in}}%
\pgfpathclose%
\pgfusepath{stroke,fill}%
\end{pgfscope}%
\begin{pgfscope}%
\pgfsetrectcap%
\pgfsetroundjoin%
\pgfsetlinewidth{2.007500pt}%
\definecolor{currentstroke}{rgb}{0.721569,0.525490,0.043137}%
\pgfsetstrokecolor{currentstroke}%
\pgfsetdash{}{0pt}%
\pgfpathmoveto{\pgfqpoint{0.923733in}{4.593244in}}%
\pgfpathlineto{\pgfqpoint{1.290400in}{4.593244in}}%
\pgfusepath{stroke}%
\end{pgfscope}%
\begin{pgfscope}%
\pgftext[x=1.512622in,y=4.496022in,left,base]{\rmfamily\fontsize{20.000000}{24.000000}\selectfont \(\displaystyle f_{\alpha,\beta,r}(c)\), \(\displaystyle c< \widehat{c}_{\beta}\)}%
\end{pgfscope}%
\begin{pgfscope}%
\pgfsetbuttcap%
\pgfsetroundjoin%
\pgfsetlinewidth{3.011250pt}%
\definecolor{currentstroke}{rgb}{0.000000,0.000000,1.000000}%
\pgfsetstrokecolor{currentstroke}%
\pgfsetdash{{19.200000pt}{4.800000pt}{3.000000pt}{4.800000pt}}{0.000000pt}%
\pgfpathmoveto{\pgfqpoint{0.923733in}{4.156476in}}%
\pgfpathlineto{\pgfqpoint{1.290400in}{4.156476in}}%
\pgfusepath{stroke}%
\end{pgfscope}%
\begin{pgfscope}%
\pgftext[x=1.512622in,y=4.059254in,left,base]{\rmfamily\fontsize{20.000000}{24.000000}\selectfont \(\displaystyle f_{\alpha,\beta,2}(c)\), \(\displaystyle c\geq \widehat{c}_{\beta}\)}%
\end{pgfscope}%
\begin{pgfscope}%
\pgfsetbuttcap%
\pgfsetroundjoin%
\pgfsetlinewidth{2.509375pt}%
\definecolor{currentstroke}{rgb}{1.000000,0.000000,0.000000}%
\pgfsetstrokecolor{currentstroke}%
\pgfsetdash{{9.250000pt}{4.000000pt}}{0.000000pt}%
\pgfpathmoveto{\pgfqpoint{0.923733in}{3.719708in}}%
\pgfpathlineto{\pgfqpoint{1.290400in}{3.719708in}}%
\pgfusepath{stroke}%
\end{pgfscope}%
\begin{pgfscope}%
\pgftext[x=1.512622in,y=3.622486in,left,base]{\rmfamily\fontsize{20.000000}{24.000000}\selectfont \(\displaystyle f_{\alpha,\beta,1}(c)\), \(\displaystyle c\geq \widehat{c}_{\beta}\)}%
\end{pgfscope}%
\end{pgfpicture}%
\makeatother%
\endgroup%

%% file: Fig2d.pgf
\begingroup%
\makeatletter%
\begin{pgfpicture}%
\pgfpathrectangle{\pgfpointorigin}{\pgfqpoint{13.842056in}{0.852511in}}%
\pgfusepath{use as bounding box, clip}%
\begin{pgfscope}%
\pgfsetbuttcap%
\pgfsetmiterjoin%
\definecolor{currentfill}{rgb}{1.000000,1.000000,1.000000}%
\pgfsetfillcolor{currentfill}%
\pgfsetlinewidth{0.000000pt}%
\definecolor{currentstroke}{rgb}{1.000000,1.000000,1.000000}%
\pgfsetstrokecolor{currentstroke}%
\pgfsetdash{}{0pt}%
\pgfpathmoveto{\pgfqpoint{0.000000in}{0.000000in}}%
\pgfpathlineto{\pgfqpoint{13.842056in}{0.000000in}}%
\pgfpathlineto{\pgfqpoint{13.842056in}{0.852511in}}%
\pgfpathlineto{\pgfqpoint{0.000000in}{0.852511in}}%
\pgfpathclose%
\pgfusepath{fill}%
\end{pgfscope}%
\begin{pgfscope}%
\pgfsetbuttcap%
\pgfsetmiterjoin%
\definecolor{currentfill}{rgb}{1.000000,1.000000,1.000000}%
\pgfsetfillcolor{currentfill}%
\pgfsetfillopacity{0.800000}%
\pgfsetlinewidth{1.003750pt}%
\definecolor{currentstroke}{rgb}{0.800000,0.800000,0.800000}%
\pgfsetstrokecolor{currentstroke}%
\pgfsetstrokeopacity{0.800000}%
\pgfsetdash{}{0pt}%
\pgfpathmoveto{\pgfqpoint{0.720778in}{0.000000in}}%
\pgfpathlineto{\pgfqpoint{13.508278in}{0.000000in}}%
\pgfpathlineto{\pgfqpoint{13.508278in}{0.775399in}}%
\pgfpathlineto{\pgfqpoint{0.720778in}{0.775399in}}%
\pgfpathclose%
\pgfusepath{stroke,fill}%
\end{pgfscope}%
\begin{pgfscope}%
\pgfsetbuttcap%
\pgfsetmiterjoin%
\definecolor{currentfill}{rgb}{0.000000,0.000000,1.000000}%
\pgfsetfillcolor{currentfill}%
\pgfsetfillopacity{0.450000}%
\pgfsetlinewidth{0.000000pt}%
\definecolor{currentstroke}{rgb}{0.000000,0.000000,1.000000}%
\pgfsetstrokecolor{currentstroke}%
\pgfsetstrokeopacity{0.450000}%
\pgfsetdash{}{0pt}%
\pgfsys@defobject{currentmarker}{\pgfqpoint{-0.111111in}{-0.111111in}}{\pgfqpoint{0.111111in}{0.111111in}}{%
\pgfpathmoveto{\pgfqpoint{-0.111111in}{-0.111111in}}%
\pgfpathlineto{\pgfqpoint{0.111111in}{-0.111111in}}%
\pgfpathlineto{\pgfqpoint{0.111111in}{0.111111in}}%
\pgfpathlineto{\pgfqpoint{-0.111111in}{0.111111in}}%
\pgfpathclose%
\pgfusepath{fill}%
}%
\begin{pgfscope}%
\pgfsys@transformshift{1.037444in}{0.445137in}%
\pgfsys@useobject{currentmarker}{}%
\end{pgfscope}%
\end{pgfscope}%
\begin{pgfscope}%
\pgftext[x=1.327722in,y=0.352776in,left,base]{\rmfamily\fontsize{19.000000}{22.800000}\selectfont \(\displaystyle D_1\)}%
\end{pgfscope}%
\begin{pgfscope}%
\pgfsetbuttcap%
\pgfsetmiterjoin%
\definecolor{currentfill}{rgb}{0.854902,0.647059,0.125490}%
\pgfsetfillcolor{currentfill}%
\pgfsetfillopacity{0.650000}%
\pgfsetlinewidth{0.000000pt}%
\definecolor{currentstroke}{rgb}{0.854902,0.647059,0.125490}%
\pgfsetstrokecolor{currentstroke}%
\pgfsetstrokeopacity{0.650000}%
\pgfsetdash{}{0pt}%
\pgfsys@defobject{currentmarker}{\pgfqpoint{-0.111111in}{-0.111111in}}{\pgfqpoint{0.111111in}{0.111111in}}{%
\pgfpathmoveto{\pgfqpoint{-0.111111in}{-0.111111in}}%
\pgfpathlineto{\pgfqpoint{0.111111in}{-0.111111in}}%
\pgfpathlineto{\pgfqpoint{0.111111in}{0.111111in}}%
\pgfpathlineto{\pgfqpoint{-0.111111in}{0.111111in}}%
\pgfpathclose%
\pgfusepath{fill}%
}%
\begin{pgfscope}%
\pgfsys@transformshift{2.302281in}{0.445137in}%
\pgfsys@useobject{currentmarker}{}%
\end{pgfscope}%
\end{pgfscope}%
\begin{pgfscope}%
\pgftext[x=2.592559in,y=0.352776in,left,base]{\rmfamily\fontsize{19.000000}{22.800000}\selectfont \(\displaystyle D_2\)}%
\end{pgfscope}%
\begin{pgfscope}%
\pgfsetbuttcap%
\pgfsetmiterjoin%
\definecolor{currentfill}{rgb}{1.000000,0.000000,0.000000}%
\pgfsetfillcolor{currentfill}%
\pgfsetfillopacity{0.550000}%
\pgfsetlinewidth{0.000000pt}%
\definecolor{currentstroke}{rgb}{1.000000,0.000000,0.000000}%
\pgfsetstrokecolor{currentstroke}%
\pgfsetstrokeopacity{0.550000}%
\pgfsetdash{}{0pt}%
\pgfsys@defobject{currentmarker}{\pgfqpoint{-0.111111in}{-0.111111in}}{\pgfqpoint{0.111111in}{0.111111in}}{%
\pgfpathmoveto{\pgfqpoint{-0.111111in}{-0.111111in}}%
\pgfpathlineto{\pgfqpoint{0.111111in}{-0.111111in}}%
\pgfpathlineto{\pgfqpoint{0.111111in}{0.111111in}}%
\pgfpathlineto{\pgfqpoint{-0.111111in}{0.111111in}}%
\pgfpathclose%
\pgfusepath{fill}%
}%
\begin{pgfscope}%
\pgfsys@transformshift{3.567118in}{0.445137in}%
\pgfsys@useobject{currentmarker}{}%
\end{pgfscope}%
\end{pgfscope}%
\begin{pgfscope}%
\pgftext[x=3.857396in,y=0.352776in,left,base]{\rmfamily\fontsize{19.000000}{22.800000}\selectfont \(\displaystyle D_3\)}%
\end{pgfscope}%
\begin{pgfscope}%
\pgfsetbuttcap%
\pgfsetroundjoin%
\pgfsetlinewidth{3.011250pt}%
\definecolor{currentstroke}{rgb}{0.000000,0.000000,0.545098}%
\pgfsetstrokecolor{currentstroke}%
\pgfsetdash{{11.100000pt}{4.800000pt}}{0.000000pt}%
\pgfpathmoveto{\pgfqpoint{4.620843in}{0.445137in}}%
\pgfpathlineto{\pgfqpoint{5.043066in}{0.445137in}}%
\pgfusepath{stroke}%
\end{pgfscope}%
\begin{pgfscope}%
\pgftext[x=5.122232in,y=0.352776in,left,base]{\rmfamily\fontsize{19.000000}{22.800000}\selectfont \(\displaystyle \beta=\frac{1}{1+s_F}\)}%
\end{pgfscope}%
\begin{pgfscope}%
\pgfsetbuttcap%
\pgfsetroundjoin%
\pgfsetlinewidth{3.011250pt}%
\definecolor{currentstroke}{rgb}{0.501961,0.000000,0.000000}%
\pgfsetstrokecolor{currentstroke}%
\pgfsetdash{{3.000000pt}{4.950000pt}}{0.000000pt}%
\pgfpathmoveto{\pgfqpoint{6.732203in}{0.445137in}}%
\pgfpathlineto{\pgfqpoint{7.154426in}{0.445137in}}%
\pgfusepath{stroke}%
\end{pgfscope}%
\begin{pgfscope}%
\pgftext[x=7.233592in,y=0.352776in,left,base]{\rmfamily\fontsize{19.000000}{22.800000}\selectfont \(\displaystyle \beta=\frac{1}{1+s_F^2}\)}%
\end{pgfscope}%
\begin{pgfscope}%
\pgfsetbuttcap%
\pgfsetroundjoin%
\pgfsetlinewidth{3.011250pt}%
\definecolor{currentstroke}{rgb}{0.000000,0.392157,0.000000}%
\pgfsetstrokecolor{currentstroke}%
\pgfsetdash{{19.200000pt}{4.800000pt}{3.000000pt}{4.800000pt}}{0.000000pt}%
\pgfpathmoveto{\pgfqpoint{8.843563in}{0.445137in}}%
\pgfpathlineto{\pgfqpoint{9.265786in}{0.445137in}}%
\pgfusepath{stroke}%
\end{pgfscope}%
\begin{pgfscope}%
\pgftext[x=9.344952in,y=0.352776in,left,base]{\rmfamily\fontsize{19.000000}{22.800000}\selectfont \(\displaystyle \alpha{\beta\left(4(1-\beta)^2-s_F^2\right)}={1-\beta(1+s_F^2)}\;\quad\)}%
\end{pgfscope}%
\end{pgfpicture}%
\makeatother%
\endgroup%

%% file: Fig3.pgf
\begingroup%
\makeatletter%
\begin{pgfpicture}%
\pgfpathrectangle{\pgfpointorigin}{\pgfqpoint{5.648213in}{3.931928in}}%
\pgfusepath{use as bounding box, clip}%
\begin{pgfscope}%
\pgfsetbuttcap%
\pgfsetmiterjoin%
\definecolor{currentfill}{rgb}{1.000000,1.000000,1.000000}%
\pgfsetfillcolor{currentfill}%
\pgfsetlinewidth{0.000000pt}%
\definecolor{currentstroke}{rgb}{1.000000,1.000000,1.000000}%
\pgfsetstrokecolor{currentstroke}%
\pgfsetdash{}{0pt}%
\pgfpathmoveto{\pgfqpoint{0.000000in}{0.000000in}}%
\pgfpathlineto{\pgfqpoint{5.648213in}{0.000000in}}%
\pgfpathlineto{\pgfqpoint{5.648213in}{3.931928in}}%
\pgfpathlineto{\pgfqpoint{0.000000in}{3.931928in}}%
\pgfpathclose%
\pgfusepath{fill}%
\end{pgfscope}%
\begin{pgfscope}%
\pgfsetbuttcap%
\pgfsetmiterjoin%
\definecolor{currentfill}{rgb}{1.000000,1.000000,1.000000}%
\pgfsetfillcolor{currentfill}%
\pgfsetlinewidth{0.000000pt}%
\definecolor{currentstroke}{rgb}{0.000000,0.000000,0.000000}%
\pgfsetstrokecolor{currentstroke}%
\pgfsetstrokeopacity{0.000000}%
\pgfsetdash{}{0pt}%
\pgfpathmoveto{\pgfqpoint{0.612802in}{0.754094in}}%
\pgfpathlineto{\pgfqpoint{5.482270in}{0.754094in}}%
\pgfpathlineto{\pgfqpoint{5.482270in}{3.774094in}}%
\pgfpathlineto{\pgfqpoint{0.612802in}{3.774094in}}%
\pgfpathclose%
\pgfusepath{fill}%
\end{pgfscope}%
\begin{pgfscope}%
\pgfsetbuttcap%
\pgfsetroundjoin%
\definecolor{currentfill}{rgb}{0.000000,0.000000,0.000000}%
\pgfsetfillcolor{currentfill}%
\pgfsetlinewidth{0.803000pt}%
\definecolor{currentstroke}{rgb}{0.000000,0.000000,0.000000}%
\pgfsetstrokecolor{currentstroke}%
\pgfsetdash{}{0pt}%
\pgfsys@defobject{currentmarker}{\pgfqpoint{0.000000in}{-0.048611in}}{\pgfqpoint{0.000000in}{0.000000in}}{%
\pgfpathmoveto{\pgfqpoint{0.000000in}{0.000000in}}%
\pgfpathlineto{\pgfqpoint{0.000000in}{-0.048611in}}%
\pgfusepath{stroke,fill}%
}%
\begin{pgfscope}%
\pgfsys@transformshift{0.612802in}{0.754094in}%
\pgfsys@useobject{currentmarker}{}%
\end{pgfscope}%
\end{pgfscope}%
\begin{pgfscope}%
\pgftext[x=0.612802in,y=0.656872in,,top]{\rmfamily\fontsize{12.000000}{14.400000}\selectfont \(\displaystyle 0\)}%
\end{pgfscope}%
\begin{pgfscope}%
\pgfsetbuttcap%
\pgfsetroundjoin%
\definecolor{currentfill}{rgb}{0.000000,0.000000,0.000000}%
\pgfsetfillcolor{currentfill}%
\pgfsetlinewidth{0.803000pt}%
\definecolor{currentstroke}{rgb}{0.000000,0.000000,0.000000}%
\pgfsetstrokecolor{currentstroke}%
\pgfsetdash{}{0pt}%
\pgfsys@defobject{currentmarker}{\pgfqpoint{0.000000in}{-0.048611in}}{\pgfqpoint{0.000000in}{0.000000in}}{%
\pgfpathmoveto{\pgfqpoint{0.000000in}{0.000000in}}%
\pgfpathlineto{\pgfqpoint{0.000000in}{-0.048611in}}%
\pgfusepath{stroke,fill}%
}%
\begin{pgfscope}%
\pgfsys@transformshift{1.830169in}{0.754094in}%
\pgfsys@useobject{currentmarker}{}%
\end{pgfscope}%
\end{pgfscope}%
\begin{pgfscope}%
\pgftext[x=1.830169in,y=0.656872in,,top]{\rmfamily\fontsize{12.000000}{14.400000}\selectfont \(\displaystyle \frac{\pi}{8}\)}%
\end{pgfscope}%
\begin{pgfscope}%
\pgfsetbuttcap%
\pgfsetroundjoin%
\definecolor{currentfill}{rgb}{0.000000,0.000000,0.000000}%
\pgfsetfillcolor{currentfill}%
\pgfsetlinewidth{0.803000pt}%
\definecolor{currentstroke}{rgb}{0.000000,0.000000,0.000000}%
\pgfsetstrokecolor{currentstroke}%
\pgfsetdash{}{0pt}%
\pgfsys@defobject{currentmarker}{\pgfqpoint{0.000000in}{-0.048611in}}{\pgfqpoint{0.000000in}{0.000000in}}{%
\pgfpathmoveto{\pgfqpoint{0.000000in}{0.000000in}}%
\pgfpathlineto{\pgfqpoint{0.000000in}{-0.048611in}}%
\pgfusepath{stroke,fill}%
}%
\begin{pgfscope}%
\pgfsys@transformshift{3.047536in}{0.754094in}%
\pgfsys@useobject{currentmarker}{}%
\end{pgfscope}%
\end{pgfscope}%
\begin{pgfscope}%
\pgftext[x=3.047536in,y=0.656872in,,top]{\rmfamily\fontsize{12.000000}{14.400000}\selectfont \(\displaystyle \frac{\pi}{4}\)}%
\end{pgfscope}%
\begin{pgfscope}%
\pgfsetbuttcap%
\pgfsetroundjoin%
\definecolor{currentfill}{rgb}{0.000000,0.000000,0.000000}%
\pgfsetfillcolor{currentfill}%
\pgfsetlinewidth{0.803000pt}%
\definecolor{currentstroke}{rgb}{0.000000,0.000000,0.000000}%
\pgfsetstrokecolor{currentstroke}%
\pgfsetdash{}{0pt}%
\pgfsys@defobject{currentmarker}{\pgfqpoint{0.000000in}{-0.048611in}}{\pgfqpoint{0.000000in}{0.000000in}}{%
\pgfpathmoveto{\pgfqpoint{0.000000in}{0.000000in}}%
\pgfpathlineto{\pgfqpoint{0.000000in}{-0.048611in}}%
\pgfusepath{stroke,fill}%
}%
\begin{pgfscope}%
\pgfsys@transformshift{4.264903in}{0.754094in}%
\pgfsys@useobject{currentmarker}{}%
\end{pgfscope}%
\end{pgfscope}%
\begin{pgfscope}%
\pgftext[x=4.264903in,y=0.656872in,,top]{\rmfamily\fontsize{12.000000}{14.400000}\selectfont \(\displaystyle \frac{3\pi}{8}\)}%
\end{pgfscope}%
\begin{pgfscope}%
\pgfsetbuttcap%
\pgfsetroundjoin%
\definecolor{currentfill}{rgb}{0.000000,0.000000,0.000000}%
\pgfsetfillcolor{currentfill}%
\pgfsetlinewidth{0.803000pt}%
\definecolor{currentstroke}{rgb}{0.000000,0.000000,0.000000}%
\pgfsetstrokecolor{currentstroke}%
\pgfsetdash{}{0pt}%
\pgfsys@defobject{currentmarker}{\pgfqpoint{0.000000in}{-0.048611in}}{\pgfqpoint{0.000000in}{0.000000in}}{%
\pgfpathmoveto{\pgfqpoint{0.000000in}{0.000000in}}%
\pgfpathlineto{\pgfqpoint{0.000000in}{-0.048611in}}%
\pgfusepath{stroke,fill}%
}%
\begin{pgfscope}%
\pgfsys@transformshift{5.482270in}{0.754094in}%
\pgfsys@useobject{currentmarker}{}%
\end{pgfscope}%
\end{pgfscope}%
\begin{pgfscope}%
\pgftext[x=5.482270in,y=0.656872in,,top]{\rmfamily\fontsize{12.000000}{14.400000}\selectfont \(\displaystyle \frac{\pi}{2}\)}%
\end{pgfscope}%
\begin{pgfscope}%
\pgftext[x=3.047536in,y=0.266833in,,top]{\rmfamily\fontsize{12.000000}{14.400000}\selectfont Friedrichs angle (in radians)}%
\end{pgfscope}%
\begin{pgfscope}%
\pgfsetbuttcap%
\pgfsetroundjoin%
\definecolor{currentfill}{rgb}{0.000000,0.000000,0.000000}%
\pgfsetfillcolor{currentfill}%
\pgfsetlinewidth{0.803000pt}%
\definecolor{currentstroke}{rgb}{0.000000,0.000000,0.000000}%
\pgfsetstrokecolor{currentstroke}%
\pgfsetdash{}{0pt}%
\pgfsys@defobject{currentmarker}{\pgfqpoint{-0.048611in}{0.000000in}}{\pgfqpoint{0.000000in}{0.000000in}}{%
\pgfpathmoveto{\pgfqpoint{0.000000in}{0.000000in}}%
\pgfpathlineto{\pgfqpoint{-0.048611in}{0.000000in}}%
\pgfusepath{stroke,fill}%
}%
\begin{pgfscope}%
\pgfsys@transformshift{0.612802in}{0.754094in}%
\pgfsys@useobject{currentmarker}{}%
\end{pgfscope}%
\end{pgfscope}%
\begin{pgfscope}%
\pgftext[x=0.307055in,y=0.696261in,left,base]{\rmfamily\fontsize{12.000000}{14.400000}\selectfont \(\displaystyle 0.0\)}%
\end{pgfscope}%
\begin{pgfscope}%
\pgfsetbuttcap%
\pgfsetroundjoin%
\definecolor{currentfill}{rgb}{0.000000,0.000000,0.000000}%
\pgfsetfillcolor{currentfill}%
\pgfsetlinewidth{0.803000pt}%
\definecolor{currentstroke}{rgb}{0.000000,0.000000,0.000000}%
\pgfsetstrokecolor{currentstroke}%
\pgfsetdash{}{0pt}%
\pgfsys@defobject{currentmarker}{\pgfqpoint{-0.048611in}{0.000000in}}{\pgfqpoint{0.000000in}{0.000000in}}{%
\pgfpathmoveto{\pgfqpoint{0.000000in}{0.000000in}}%
\pgfpathlineto{\pgfqpoint{-0.048611in}{0.000000in}}%
\pgfusepath{stroke,fill}%
}%
\begin{pgfscope}%
\pgfsys@transformshift{0.612802in}{1.358094in}%
\pgfsys@useobject{currentmarker}{}%
\end{pgfscope}%
\end{pgfscope}%
\begin{pgfscope}%
\pgftext[x=0.307055in,y=1.300261in,left,base]{\rmfamily\fontsize{12.000000}{14.400000}\selectfont \(\displaystyle 0.2\)}%
\end{pgfscope}%
\begin{pgfscope}%
\pgfsetbuttcap%
\pgfsetroundjoin%
\definecolor{currentfill}{rgb}{0.000000,0.000000,0.000000}%
\pgfsetfillcolor{currentfill}%
\pgfsetlinewidth{0.803000pt}%
\definecolor{currentstroke}{rgb}{0.000000,0.000000,0.000000}%
\pgfsetstrokecolor{currentstroke}%
\pgfsetdash{}{0pt}%
\pgfsys@defobject{currentmarker}{\pgfqpoint{-0.048611in}{0.000000in}}{\pgfqpoint{0.000000in}{0.000000in}}{%
\pgfpathmoveto{\pgfqpoint{0.000000in}{0.000000in}}%
\pgfpathlineto{\pgfqpoint{-0.048611in}{0.000000in}}%
\pgfusepath{stroke,fill}%
}%
\begin{pgfscope}%
\pgfsys@transformshift{0.612802in}{1.962094in}%
\pgfsys@useobject{currentmarker}{}%
\end{pgfscope}%
\end{pgfscope}%
\begin{pgfscope}%
\pgftext[x=0.307055in,y=1.904261in,left,base]{\rmfamily\fontsize{12.000000}{14.400000}\selectfont \(\displaystyle 0.4\)}%
\end{pgfscope}%
\begin{pgfscope}%
\pgfsetbuttcap%
\pgfsetroundjoin%
\definecolor{currentfill}{rgb}{0.000000,0.000000,0.000000}%
\pgfsetfillcolor{currentfill}%
\pgfsetlinewidth{0.803000pt}%
\definecolor{currentstroke}{rgb}{0.000000,0.000000,0.000000}%
\pgfsetstrokecolor{currentstroke}%
\pgfsetdash{}{0pt}%
\pgfsys@defobject{currentmarker}{\pgfqpoint{-0.048611in}{0.000000in}}{\pgfqpoint{0.000000in}{0.000000in}}{%
\pgfpathmoveto{\pgfqpoint{0.000000in}{0.000000in}}%
\pgfpathlineto{\pgfqpoint{-0.048611in}{0.000000in}}%
\pgfusepath{stroke,fill}%
}%
\begin{pgfscope}%
\pgfsys@transformshift{0.612802in}{2.566094in}%
\pgfsys@useobject{currentmarker}{}%
\end{pgfscope}%
\end{pgfscope}%
\begin{pgfscope}%
\pgftext[x=0.307055in,y=2.508261in,left,base]{\rmfamily\fontsize{12.000000}{14.400000}\selectfont \(\displaystyle 0.6\)}%
\end{pgfscope}%
\begin{pgfscope}%
\pgfsetbuttcap%
\pgfsetroundjoin%
\definecolor{currentfill}{rgb}{0.000000,0.000000,0.000000}%
\pgfsetfillcolor{currentfill}%
\pgfsetlinewidth{0.803000pt}%
\definecolor{currentstroke}{rgb}{0.000000,0.000000,0.000000}%
\pgfsetstrokecolor{currentstroke}%
\pgfsetdash{}{0pt}%
\pgfsys@defobject{currentmarker}{\pgfqpoint{-0.048611in}{0.000000in}}{\pgfqpoint{0.000000in}{0.000000in}}{%
\pgfpathmoveto{\pgfqpoint{0.000000in}{0.000000in}}%
\pgfpathlineto{\pgfqpoint{-0.048611in}{0.000000in}}%
\pgfusepath{stroke,fill}%
}%
\begin{pgfscope}%
\pgfsys@transformshift{0.612802in}{3.170094in}%
\pgfsys@useobject{currentmarker}{}%
\end{pgfscope}%
\end{pgfscope}%
\begin{pgfscope}%
\pgftext[x=0.307055in,y=3.112261in,left,base]{\rmfamily\fontsize{12.000000}{14.400000}\selectfont \(\displaystyle 0.8\)}%
\end{pgfscope}%
\begin{pgfscope}%
\pgfsetbuttcap%
\pgfsetroundjoin%
\definecolor{currentfill}{rgb}{0.000000,0.000000,0.000000}%
\pgfsetfillcolor{currentfill}%
\pgfsetlinewidth{0.803000pt}%
\definecolor{currentstroke}{rgb}{0.000000,0.000000,0.000000}%
\pgfsetstrokecolor{currentstroke}%
\pgfsetdash{}{0pt}%
\pgfsys@defobject{currentmarker}{\pgfqpoint{-0.048611in}{0.000000in}}{\pgfqpoint{0.000000in}{0.000000in}}{%
\pgfpathmoveto{\pgfqpoint{0.000000in}{0.000000in}}%
\pgfpathlineto{\pgfqpoint{-0.048611in}{0.000000in}}%
\pgfusepath{stroke,fill}%
}%
\begin{pgfscope}%
\pgfsys@transformshift{0.612802in}{3.774094in}%
\pgfsys@useobject{currentmarker}{}%
\end{pgfscope}%
\end{pgfscope}%
\begin{pgfscope}%
\pgftext[x=0.307055in,y=3.716261in,left,base]{\rmfamily\fontsize{12.000000}{14.400000}\selectfont \(\displaystyle 1.0\)}%
\end{pgfscope}%
\begin{pgfscope}%
\pgftext[x=0.251500in,y=2.264094in,,bottom,rotate=90.000000]{\rmfamily\fontsize{12.000000}{14.400000}\selectfont Rate of convergence}%
\end{pgfscope}%
\begin{pgfscope}%
\pgfpathrectangle{\pgfqpoint{0.612802in}{0.754094in}}{\pgfqpoint{4.869469in}{3.020000in}} %
\pgfusepath{clip}%
\pgfsetbuttcap%
\pgfsetroundjoin%
\pgfsetlinewidth{3.011250pt}%
\definecolor{currentstroke}{rgb}{1.000000,0.000000,0.000000}%
\pgfsetstrokecolor{currentstroke}%
\pgfsetdash{{11.100000pt}{4.800000pt}}{0.000000pt}%
\pgfpathmoveto{\pgfqpoint{0.612802in}{3.774094in}}%
\pgfpathlineto{\pgfqpoint{0.661988in}{3.773334in}}%
\pgfpathlineto{\pgfqpoint{0.711175in}{3.771054in}}%
\pgfpathlineto{\pgfqpoint{0.760361in}{3.767257in}}%
\pgfpathlineto{\pgfqpoint{0.809548in}{3.761946in}}%
\pgfpathlineto{\pgfqpoint{0.858734in}{3.755127in}}%
\pgfpathlineto{\pgfqpoint{0.907921in}{3.746807in}}%
\pgfpathlineto{\pgfqpoint{0.957108in}{3.736993in}}%
\pgfpathlineto{\pgfqpoint{1.006294in}{3.725697in}}%
\pgfpathlineto{\pgfqpoint{1.055481in}{3.712929in}}%
\pgfpathlineto{\pgfqpoint{1.104667in}{3.698702in}}%
\pgfpathlineto{\pgfqpoint{1.153854in}{3.683030in}}%
\pgfpathlineto{\pgfqpoint{1.203040in}{3.665930in}}%
\pgfpathlineto{\pgfqpoint{1.252227in}{3.647418in}}%
\pgfpathlineto{\pgfqpoint{1.301413in}{3.627514in}}%
\pgfpathlineto{\pgfqpoint{1.350600in}{3.606236in}}%
\pgfpathlineto{\pgfqpoint{1.399786in}{3.583607in}}%
\pgfpathlineto{\pgfqpoint{1.448973in}{3.559649in}}%
\pgfpathlineto{\pgfqpoint{1.498160in}{3.534387in}}%
\pgfpathlineto{\pgfqpoint{1.547346in}{3.507846in}}%
\pgfpathlineto{\pgfqpoint{1.596533in}{3.480053in}}%
\pgfpathlineto{\pgfqpoint{1.645719in}{3.451035in}}%
\pgfpathlineto{\pgfqpoint{1.694906in}{3.420822in}}%
\pgfpathlineto{\pgfqpoint{1.744092in}{3.389444in}}%
\pgfpathlineto{\pgfqpoint{1.793279in}{3.356933in}}%
\pgfpathlineto{\pgfqpoint{1.842465in}{3.323322in}}%
\pgfpathlineto{\pgfqpoint{1.891652in}{3.288644in}}%
\pgfpathlineto{\pgfqpoint{1.940839in}{3.252934in}}%
\pgfpathlineto{\pgfqpoint{1.990025in}{3.216229in}}%
\pgfpathlineto{\pgfqpoint{2.039212in}{3.178565in}}%
\pgfpathlineto{\pgfqpoint{2.088398in}{3.139980in}}%
\pgfpathlineto{\pgfqpoint{2.137585in}{3.100514in}}%
\pgfpathlineto{\pgfqpoint{2.186771in}{3.060205in}}%
\pgfpathlineto{\pgfqpoint{2.235958in}{3.019094in}}%
\pgfpathlineto{\pgfqpoint{2.285144in}{2.977224in}}%
\pgfpathlineto{\pgfqpoint{2.334331in}{2.934635in}}%
\pgfpathlineto{\pgfqpoint{2.383518in}{2.891371in}}%
\pgfpathlineto{\pgfqpoint{2.432704in}{2.847476in}}%
\pgfpathlineto{\pgfqpoint{2.481891in}{2.802993in}}%
\pgfpathlineto{\pgfqpoint{2.531077in}{2.757967in}}%
\pgfpathlineto{\pgfqpoint{2.580264in}{2.712444in}}%
\pgfpathlineto{\pgfqpoint{2.629450in}{2.666470in}}%
\pgfpathlineto{\pgfqpoint{2.678637in}{2.620090in}}%
\pgfpathlineto{\pgfqpoint{2.727823in}{2.573353in}}%
\pgfpathlineto{\pgfqpoint{2.777010in}{2.526303in}}%
\pgfpathlineto{\pgfqpoint{2.826196in}{2.478990in}}%
\pgfpathlineto{\pgfqpoint{2.875383in}{2.431460in}}%
\pgfpathlineto{\pgfqpoint{2.924570in}{2.383762in}}%
\pgfpathlineto{\pgfqpoint{2.973756in}{2.335943in}}%
\pgfpathlineto{\pgfqpoint{3.022943in}{2.288052in}}%
\pgfpathlineto{\pgfqpoint{3.072129in}{2.240137in}}%
\pgfpathlineto{\pgfqpoint{3.121316in}{2.192246in}}%
\pgfpathlineto{\pgfqpoint{3.170502in}{2.144427in}}%
\pgfpathlineto{\pgfqpoint{3.219689in}{2.096729in}}%
\pgfpathlineto{\pgfqpoint{3.268875in}{2.049199in}}%
\pgfpathlineto{\pgfqpoint{3.318062in}{2.001886in}}%
\pgfpathlineto{\pgfqpoint{3.367249in}{1.954836in}}%
\pgfpathlineto{\pgfqpoint{3.416435in}{1.908098in}}%
\pgfpathlineto{\pgfqpoint{3.465622in}{1.861719in}}%
\pgfpathlineto{\pgfqpoint{3.514808in}{1.815745in}}%
\pgfpathlineto{\pgfqpoint{3.563995in}{1.770222in}}%
\pgfpathlineto{\pgfqpoint{3.613181in}{1.725196in}}%
\pgfpathlineto{\pgfqpoint{3.662368in}{1.680713in}}%
\pgfpathlineto{\pgfqpoint{3.711554in}{1.636818in}}%
\pgfpathlineto{\pgfqpoint{3.760741in}{1.593554in}}%
\pgfpathlineto{\pgfqpoint{3.809928in}{1.550965in}}%
\pgfpathlineto{\pgfqpoint{3.859114in}{1.509094in}}%
\pgfpathlineto{\pgfqpoint{3.908301in}{1.467984in}}%
\pgfpathlineto{\pgfqpoint{3.957487in}{1.427675in}}%
\pgfpathlineto{\pgfqpoint{4.006674in}{1.388209in}}%
\pgfpathlineto{\pgfqpoint{4.055860in}{1.349624in}}%
\pgfpathlineto{\pgfqpoint{4.105047in}{1.311960in}}%
\pgfpathlineto{\pgfqpoint{4.154233in}{1.275255in}}%
\pgfpathlineto{\pgfqpoint{4.203420in}{1.239545in}}%
\pgfpathlineto{\pgfqpoint{4.252606in}{1.204867in}}%
\pgfpathlineto{\pgfqpoint{4.301793in}{1.171256in}}%
\pgfpathlineto{\pgfqpoint{4.350980in}{1.138745in}}%
\pgfpathlineto{\pgfqpoint{4.400166in}{1.107367in}}%
\pgfpathlineto{\pgfqpoint{4.449353in}{1.077154in}}%
\pgfpathlineto{\pgfqpoint{4.498539in}{1.048136in}}%
\pgfpathlineto{\pgfqpoint{4.547726in}{1.020343in}}%
\pgfpathlineto{\pgfqpoint{4.596912in}{0.993802in}}%
\pgfpathlineto{\pgfqpoint{4.646099in}{0.968539in}}%
\pgfpathlineto{\pgfqpoint{4.695285in}{0.944582in}}%
\pgfpathlineto{\pgfqpoint{4.744472in}{0.921953in}}%
\pgfpathlineto{\pgfqpoint{4.793659in}{0.900675in}}%
\pgfpathlineto{\pgfqpoint{4.842845in}{0.880771in}}%
\pgfpathlineto{\pgfqpoint{4.892032in}{0.862259in}}%
\pgfpathlineto{\pgfqpoint{4.941218in}{0.845159in}}%
\pgfpathlineto{\pgfqpoint{4.990405in}{0.829487in}}%
\pgfpathlineto{\pgfqpoint{5.039591in}{0.815260in}}%
\pgfpathlineto{\pgfqpoint{5.088778in}{0.802492in}}%
\pgfpathlineto{\pgfqpoint{5.137964in}{0.791195in}}%
\pgfpathlineto{\pgfqpoint{5.187151in}{0.781382in}}%
\pgfpathlineto{\pgfqpoint{5.236338in}{0.773062in}}%
\pgfpathlineto{\pgfqpoint{5.285524in}{0.766243in}}%
\pgfpathlineto{\pgfqpoint{5.334711in}{0.760932in}}%
\pgfpathlineto{\pgfqpoint{5.383897in}{0.757135in}}%
\pgfpathlineto{\pgfqpoint{5.433084in}{0.754855in}}%
\pgfpathlineto{\pgfqpoint{5.482270in}{0.754094in}}%
\pgfusepath{stroke}%
\end{pgfscope}%
\begin{pgfscope}%
\pgfpathrectangle{\pgfqpoint{0.612802in}{0.754094in}}{\pgfqpoint{4.869469in}{3.020000in}} %
\pgfusepath{clip}%
\pgfsetbuttcap%
\pgfsetroundjoin%
\pgfsetlinewidth{3.011250pt}%
\definecolor{currentstroke}{rgb}{0.800000,0.600000,0.000000}%
\pgfsetstrokecolor{currentstroke}%
\pgfsetdash{{19.200000pt}{4.800000pt}{3.000000pt}{4.800000pt}}{0.000000pt}%
\pgfpathmoveto{\pgfqpoint{0.612802in}{3.774094in}}%
\pgfpathlineto{\pgfqpoint{0.661988in}{3.772574in}}%
\pgfpathlineto{\pgfqpoint{0.711175in}{3.768020in}}%
\pgfpathlineto{\pgfqpoint{0.760361in}{3.760451in}}%
\pgfpathlineto{\pgfqpoint{0.809548in}{3.749895in}}%
\pgfpathlineto{\pgfqpoint{0.858734in}{3.736397in}}%
\pgfpathlineto{\pgfqpoint{0.907921in}{3.720008in}}%
\pgfpathlineto{\pgfqpoint{0.957108in}{3.700793in}}%
\pgfpathlineto{\pgfqpoint{1.006294in}{3.678826in}}%
\pgfpathlineto{\pgfqpoint{1.055481in}{3.654192in}}%
\pgfpathlineto{\pgfqpoint{1.104667in}{3.626982in}}%
\pgfpathlineto{\pgfqpoint{1.153854in}{3.597297in}}%
\pgfpathlineto{\pgfqpoint{1.203040in}{3.565246in}}%
\pgfpathlineto{\pgfqpoint{1.252227in}{3.530941in}}%
\pgfpathlineto{\pgfqpoint{1.301413in}{3.494503in}}%
\pgfpathlineto{\pgfqpoint{1.350600in}{3.456055in}}%
\pgfpathlineto{\pgfqpoint{1.399786in}{3.415724in}}%
\pgfpathlineto{\pgfqpoint{1.448973in}{3.373640in}}%
\pgfpathlineto{\pgfqpoint{1.498160in}{3.329935in}}%
\pgfpathlineto{\pgfqpoint{1.547346in}{3.284740in}}%
\pgfpathlineto{\pgfqpoint{1.596533in}{3.238189in}}%
\pgfpathlineto{\pgfqpoint{1.645719in}{3.190413in}}%
\pgfpathlineto{\pgfqpoint{1.694906in}{3.141543in}}%
\pgfpathlineto{\pgfqpoint{1.744092in}{3.091707in}}%
\pgfpathlineto{\pgfqpoint{1.793279in}{3.041031in}}%
\pgfpathlineto{\pgfqpoint{1.842465in}{2.989639in}}%
\pgfpathlineto{\pgfqpoint{1.891652in}{2.937648in}}%
\pgfpathlineto{\pgfqpoint{1.940839in}{2.885174in}}%
\pgfpathlineto{\pgfqpoint{1.990025in}{2.832330in}}%
\pgfpathlineto{\pgfqpoint{2.039212in}{2.779220in}}%
\pgfpathlineto{\pgfqpoint{2.088398in}{2.725947in}}%
\pgfpathlineto{\pgfqpoint{2.137585in}{2.672608in}}%
\pgfpathlineto{\pgfqpoint{2.186771in}{2.619295in}}%
\pgfpathlineto{\pgfqpoint{2.235958in}{2.566094in}}%
\pgfpathlineto{\pgfqpoint{2.285144in}{2.513088in}}%
\pgfpathlineto{\pgfqpoint{2.334331in}{2.460352in}}%
\pgfpathlineto{\pgfqpoint{2.383518in}{2.407959in}}%
\pgfpathlineto{\pgfqpoint{2.432704in}{2.355975in}}%
\pgfpathlineto{\pgfqpoint{2.481891in}{2.304461in}}%
\pgfpathlineto{\pgfqpoint{2.531077in}{2.253476in}}%
\pgfpathlineto{\pgfqpoint{2.580264in}{2.203071in}}%
\pgfpathlineto{\pgfqpoint{2.629450in}{2.153295in}}%
\pgfpathlineto{\pgfqpoint{2.678637in}{2.104191in}}%
\pgfpathlineto{\pgfqpoint{2.727823in}{2.055799in}}%
\pgfpathlineto{\pgfqpoint{2.777010in}{2.008155in}}%
\pgfpathlineto{\pgfqpoint{2.826196in}{1.961292in}}%
\pgfpathlineto{\pgfqpoint{2.875383in}{1.915238in}}%
\pgfpathlineto{\pgfqpoint{2.924570in}{1.870018in}}%
\pgfpathlineto{\pgfqpoint{2.973756in}{1.825656in}}%
\pgfpathlineto{\pgfqpoint{3.022943in}{1.782170in}}%
\pgfpathlineto{\pgfqpoint{3.072129in}{1.739577in}}%
\pgfpathlineto{\pgfqpoint{3.121316in}{1.697893in}}%
\pgfpathlineto{\pgfqpoint{3.170502in}{1.657128in}}%
\pgfpathlineto{\pgfqpoint{3.219689in}{1.617292in}}%
\pgfpathlineto{\pgfqpoint{3.268875in}{1.578394in}}%
\pgfpathlineto{\pgfqpoint{3.318062in}{1.540439in}}%
\pgfpathlineto{\pgfqpoint{3.367249in}{1.503433in}}%
\pgfpathlineto{\pgfqpoint{3.416435in}{1.467376in}}%
\pgfpathlineto{\pgfqpoint{3.465622in}{1.432272in}}%
\pgfpathlineto{\pgfqpoint{3.514808in}{1.398120in}}%
\pgfpathlineto{\pgfqpoint{3.563995in}{1.364919in}}%
\pgfpathlineto{\pgfqpoint{3.613181in}{1.332667in}}%
\pgfpathlineto{\pgfqpoint{3.662368in}{1.301362in}}%
\pgfpathlineto{\pgfqpoint{3.711554in}{1.271000in}}%
\pgfpathlineto{\pgfqpoint{3.760741in}{1.241576in}}%
\pgfpathlineto{\pgfqpoint{3.809928in}{1.213085in}}%
\pgfpathlineto{\pgfqpoint{3.859114in}{1.185523in}}%
\pgfpathlineto{\pgfqpoint{3.908301in}{1.158883in}}%
\pgfpathlineto{\pgfqpoint{3.957487in}{1.133158in}}%
\pgfpathlineto{\pgfqpoint{4.006674in}{1.108342in}}%
\pgfpathlineto{\pgfqpoint{4.055860in}{1.084429in}}%
\pgfpathlineto{\pgfqpoint{4.105047in}{1.061412in}}%
\pgfpathlineto{\pgfqpoint{4.154233in}{1.039282in}}%
\pgfpathlineto{\pgfqpoint{4.203420in}{1.018033in}}%
\pgfpathlineto{\pgfqpoint{4.252606in}{0.997658in}}%
\pgfpathlineto{\pgfqpoint{4.301793in}{0.978150in}}%
\pgfpathlineto{\pgfqpoint{4.350980in}{0.959501in}}%
\pgfpathlineto{\pgfqpoint{4.400166in}{0.941704in}}%
\pgfpathlineto{\pgfqpoint{4.449353in}{0.924752in}}%
\pgfpathlineto{\pgfqpoint{4.498539in}{0.908639in}}%
\pgfpathlineto{\pgfqpoint{4.547726in}{0.893357in}}%
\pgfpathlineto{\pgfqpoint{4.596912in}{0.878901in}}%
\pgfpathlineto{\pgfqpoint{4.646099in}{0.865264in}}%
\pgfpathlineto{\pgfqpoint{4.695285in}{0.852440in}}%
\pgfpathlineto{\pgfqpoint{4.744472in}{0.840423in}}%
\pgfpathlineto{\pgfqpoint{4.793659in}{0.829208in}}%
\pgfpathlineto{\pgfqpoint{4.842845in}{0.818789in}}%
\pgfpathlineto{\pgfqpoint{4.892032in}{0.809163in}}%
\pgfpathlineto{\pgfqpoint{4.941218in}{0.800324in}}%
\pgfpathlineto{\pgfqpoint{4.990405in}{0.792267in}}%
\pgfpathlineto{\pgfqpoint{5.039591in}{0.784990in}}%
\pgfpathlineto{\pgfqpoint{5.088778in}{0.778489in}}%
\pgfpathlineto{\pgfqpoint{5.137964in}{0.772760in}}%
\pgfpathlineto{\pgfqpoint{5.187151in}{0.767800in}}%
\pgfpathlineto{\pgfqpoint{5.236338in}{0.763608in}}%
\pgfpathlineto{\pgfqpoint{5.285524in}{0.760181in}}%
\pgfpathlineto{\pgfqpoint{5.334711in}{0.757517in}}%
\pgfpathlineto{\pgfqpoint{5.383897in}{0.755615in}}%
\pgfpathlineto{\pgfqpoint{5.433084in}{0.754475in}}%
\pgfpathlineto{\pgfqpoint{5.482270in}{0.754094in}}%
\pgfusepath{stroke}%
\end{pgfscope}%
\begin{pgfscope}%
\pgfpathrectangle{\pgfqpoint{0.612802in}{0.754094in}}{\pgfqpoint{4.869469in}{3.020000in}} %
\pgfusepath{clip}%
\pgfsetbuttcap%
\pgfsetroundjoin%
\pgfsetlinewidth{3.011250pt}%
\definecolor{currentstroke}{rgb}{0.000000,0.000000,1.000000}%
\pgfsetstrokecolor{currentstroke}%
\pgfsetdash{{3.000000pt}{4.950000pt}}{0.000000pt}%
\pgfpathmoveto{\pgfqpoint{0.612802in}{3.774094in}}%
\pgfpathlineto{\pgfqpoint{0.661988in}{3.773714in}}%
\pgfpathlineto{\pgfqpoint{0.711175in}{3.772574in}}%
\pgfpathlineto{\pgfqpoint{0.760361in}{3.770674in}}%
\pgfpathlineto{\pgfqpoint{0.809548in}{3.768014in}}%
\pgfpathlineto{\pgfqpoint{0.858734in}{3.764596in}}%
\pgfpathlineto{\pgfqpoint{0.907921in}{3.760420in}}%
\pgfpathlineto{\pgfqpoint{0.957108in}{3.755487in}}%
\pgfpathlineto{\pgfqpoint{1.006294in}{3.749798in}}%
\pgfpathlineto{\pgfqpoint{1.055481in}{3.743355in}}%
\pgfpathlineto{\pgfqpoint{1.104667in}{3.736160in}}%
\pgfpathlineto{\pgfqpoint{1.153854in}{3.728214in}}%
\pgfpathlineto{\pgfqpoint{1.203040in}{3.719519in}}%
\pgfpathlineto{\pgfqpoint{1.252227in}{3.710078in}}%
\pgfpathlineto{\pgfqpoint{1.301413in}{3.699892in}}%
\pgfpathlineto{\pgfqpoint{1.350600in}{3.688965in}}%
\pgfpathlineto{\pgfqpoint{1.399786in}{3.677300in}}%
\pgfpathlineto{\pgfqpoint{1.448973in}{3.664898in}}%
\pgfpathlineto{\pgfqpoint{1.498160in}{3.651763in}}%
\pgfpathlineto{\pgfqpoint{1.547346in}{3.637899in}}%
\pgfpathlineto{\pgfqpoint{1.596533in}{3.623309in}}%
\pgfpathlineto{\pgfqpoint{1.645719in}{3.607997in}}%
\pgfpathlineto{\pgfqpoint{1.694906in}{3.591966in}}%
\pgfpathlineto{\pgfqpoint{1.744092in}{3.575221in}}%
\pgfpathlineto{\pgfqpoint{1.793279in}{3.557766in}}%
\pgfpathlineto{\pgfqpoint{1.842465in}{3.539604in}}%
\pgfpathlineto{\pgfqpoint{1.891652in}{3.520742in}}%
\pgfpathlineto{\pgfqpoint{1.940839in}{3.501183in}}%
\pgfpathlineto{\pgfqpoint{1.990025in}{3.480933in}}%
\pgfpathlineto{\pgfqpoint{2.039212in}{3.459996in}}%
\pgfpathlineto{\pgfqpoint{2.088398in}{3.438377in}}%
\pgfpathlineto{\pgfqpoint{2.137585in}{3.416084in}}%
\pgfpathlineto{\pgfqpoint{2.186771in}{3.393120in}}%
\pgfpathlineto{\pgfqpoint{2.235958in}{3.369491in}}%
\pgfpathlineto{\pgfqpoint{2.285144in}{3.345204in}}%
\pgfpathlineto{\pgfqpoint{2.334331in}{3.320265in}}%
\pgfpathlineto{\pgfqpoint{2.383518in}{3.294680in}}%
\pgfpathlineto{\pgfqpoint{2.432704in}{3.268455in}}%
\pgfpathlineto{\pgfqpoint{2.481891in}{3.241598in}}%
\pgfpathlineto{\pgfqpoint{2.531077in}{3.214114in}}%
\pgfpathlineto{\pgfqpoint{2.580264in}{3.186011in}}%
\pgfpathlineto{\pgfqpoint{2.629450in}{3.157295in}}%
\pgfpathlineto{\pgfqpoint{2.678637in}{3.127975in}}%
\pgfpathlineto{\pgfqpoint{2.727823in}{3.098057in}}%
\pgfpathlineto{\pgfqpoint{2.777010in}{3.067549in}}%
\pgfpathlineto{\pgfqpoint{2.826196in}{3.036458in}}%
\pgfpathlineto{\pgfqpoint{2.875383in}{3.004793in}}%
\pgfpathlineto{\pgfqpoint{2.924570in}{2.972561in}}%
\pgfpathlineto{\pgfqpoint{2.973756in}{2.939771in}}%
\pgfpathlineto{\pgfqpoint{3.022943in}{2.906431in}}%
\pgfpathlineto{\pgfqpoint{3.072129in}{2.872549in}}%
\pgfpathlineto{\pgfqpoint{3.121316in}{2.838133in}}%
\pgfpathlineto{\pgfqpoint{3.170502in}{2.803193in}}%
\pgfpathlineto{\pgfqpoint{3.219689in}{2.767737in}}%
\pgfpathlineto{\pgfqpoint{3.268875in}{2.731774in}}%
\pgfpathlineto{\pgfqpoint{3.318062in}{2.695313in}}%
\pgfpathlineto{\pgfqpoint{3.367249in}{2.658363in}}%
\pgfpathlineto{\pgfqpoint{3.416435in}{2.620935in}}%
\pgfpathlineto{\pgfqpoint{3.465622in}{2.583036in}}%
\pgfpathlineto{\pgfqpoint{3.514808in}{2.544676in}}%
\pgfpathlineto{\pgfqpoint{3.563995in}{2.505866in}}%
\pgfpathlineto{\pgfqpoint{3.613181in}{2.466615in}}%
\pgfpathlineto{\pgfqpoint{3.662368in}{2.426933in}}%
\pgfpathlineto{\pgfqpoint{3.711554in}{2.386830in}}%
\pgfpathlineto{\pgfqpoint{3.760741in}{2.346315in}}%
\pgfpathlineto{\pgfqpoint{3.809928in}{2.305400in}}%
\pgfpathlineto{\pgfqpoint{3.859114in}{2.264094in}}%
\pgfpathlineto{\pgfqpoint{3.908301in}{2.222409in}}%
\pgfpathlineto{\pgfqpoint{3.957487in}{2.180353in}}%
\pgfpathlineto{\pgfqpoint{4.006674in}{2.137939in}}%
\pgfpathlineto{\pgfqpoint{4.055860in}{2.095176in}}%
\pgfpathlineto{\pgfqpoint{4.105047in}{2.052075in}}%
\pgfpathlineto{\pgfqpoint{4.154233in}{2.008648in}}%
\pgfpathlineto{\pgfqpoint{4.203420in}{1.964905in}}%
\pgfpathlineto{\pgfqpoint{4.252606in}{1.920857in}}%
\pgfpathlineto{\pgfqpoint{4.301793in}{1.876515in}}%
\pgfpathlineto{\pgfqpoint{4.350980in}{1.831891in}}%
\pgfpathlineto{\pgfqpoint{4.400166in}{1.786995in}}%
\pgfpathlineto{\pgfqpoint{4.449353in}{1.741840in}}%
\pgfpathlineto{\pgfqpoint{4.498539in}{1.696435in}}%
\pgfpathlineto{\pgfqpoint{4.547726in}{1.650794in}}%
\pgfpathlineto{\pgfqpoint{4.596912in}{1.604927in}}%
\pgfpathlineto{\pgfqpoint{4.646099in}{1.558845in}}%
\pgfpathlineto{\pgfqpoint{4.695285in}{1.512561in}}%
\pgfpathlineto{\pgfqpoint{4.744472in}{1.466086in}}%
\pgfpathlineto{\pgfqpoint{4.793659in}{1.419432in}}%
\pgfpathlineto{\pgfqpoint{4.842845in}{1.372611in}}%
\pgfpathlineto{\pgfqpoint{4.892032in}{1.325633in}}%
\pgfpathlineto{\pgfqpoint{4.941218in}{1.278512in}}%
\pgfpathlineto{\pgfqpoint{4.990405in}{1.231259in}}%
\pgfpathlineto{\pgfqpoint{5.039591in}{1.183885in}}%
\pgfpathlineto{\pgfqpoint{5.088778in}{1.136404in}}%
\pgfpathlineto{\pgfqpoint{5.137964in}{1.088826in}}%
\pgfpathlineto{\pgfqpoint{5.187151in}{1.041164in}}%
\pgfpathlineto{\pgfqpoint{5.236338in}{0.993429in}}%
\pgfpathlineto{\pgfqpoint{5.285524in}{0.945635in}}%
\pgfpathlineto{\pgfqpoint{5.334711in}{0.897792in}}%
\pgfpathlineto{\pgfqpoint{5.383897in}{0.849913in}}%
\pgfpathlineto{\pgfqpoint{5.433084in}{0.802010in}}%
\pgfpathlineto{\pgfqpoint{5.482270in}{0.754094in}}%
\pgfusepath{stroke}%
\end{pgfscope}%
\begin{pgfscope}%
\pgfpathrectangle{\pgfqpoint{0.612802in}{0.754094in}}{\pgfqpoint{4.869469in}{3.020000in}} %
\pgfusepath{clip}%
\pgfsetrectcap%
\pgfsetroundjoin%
\pgfsetlinewidth{3.011250pt}%
\definecolor{currentstroke}{rgb}{0.000000,0.501961,0.000000}%
\pgfsetstrokecolor{currentstroke}%
\pgfsetdash{}{0pt}%
\pgfpathmoveto{\pgfqpoint{0.612802in}{3.774094in}}%
\pgfpathlineto{\pgfqpoint{0.661988in}{3.679761in}}%
\pgfpathlineto{\pgfqpoint{0.711175in}{3.588351in}}%
\pgfpathlineto{\pgfqpoint{0.760361in}{3.499753in}}%
\pgfpathlineto{\pgfqpoint{0.809548in}{3.413861in}}%
\pgfpathlineto{\pgfqpoint{0.858734in}{3.330574in}}%
\pgfpathlineto{\pgfqpoint{0.907921in}{3.249794in}}%
\pgfpathlineto{\pgfqpoint{0.957108in}{3.171430in}}%
\pgfpathlineto{\pgfqpoint{1.006294in}{3.095394in}}%
\pgfpathlineto{\pgfqpoint{1.055481in}{3.021603in}}%
\pgfpathlineto{\pgfqpoint{1.104667in}{2.949978in}}%
\pgfpathlineto{\pgfqpoint{1.153854in}{2.880441in}}%
\pgfpathlineto{\pgfqpoint{1.203040in}{2.812920in}}%
\pgfpathlineto{\pgfqpoint{1.252227in}{2.747347in}}%
\pgfpathlineto{\pgfqpoint{1.301413in}{2.683654in}}%
\pgfpathlineto{\pgfqpoint{1.350600in}{2.621779in}}%
\pgfpathlineto{\pgfqpoint{1.399786in}{2.561661in}}%
\pgfpathlineto{\pgfqpoint{1.448973in}{2.503242in}}%
\pgfpathlineto{\pgfqpoint{1.498160in}{2.446466in}}%
\pgfpathlineto{\pgfqpoint{1.547346in}{2.391281in}}%
\pgfpathlineto{\pgfqpoint{1.596533in}{2.337636in}}%
\pgfpathlineto{\pgfqpoint{1.645719in}{2.285482in}}%
\pgfpathlineto{\pgfqpoint{1.694906in}{2.234772in}}%
\pgfpathlineto{\pgfqpoint{1.744092in}{2.185463in}}%
\pgfpathlineto{\pgfqpoint{1.793279in}{2.137510in}}%
\pgfpathlineto{\pgfqpoint{1.842465in}{2.090874in}}%
\pgfpathlineto{\pgfqpoint{1.891652in}{2.045514in}}%
\pgfpathlineto{\pgfqpoint{1.940839in}{2.001394in}}%
\pgfpathlineto{\pgfqpoint{1.990025in}{1.958477in}}%
\pgfpathlineto{\pgfqpoint{2.039212in}{1.916727in}}%
\pgfpathlineto{\pgfqpoint{2.088398in}{1.876112in}}%
\pgfpathlineto{\pgfqpoint{2.137585in}{1.836600in}}%
\pgfpathlineto{\pgfqpoint{2.186771in}{1.798159in}}%
\pgfpathlineto{\pgfqpoint{2.235958in}{1.760761in}}%
\pgfpathlineto{\pgfqpoint{2.285144in}{1.724377in}}%
\pgfpathlineto{\pgfqpoint{2.334331in}{1.688979in}}%
\pgfpathlineto{\pgfqpoint{2.383518in}{1.654541in}}%
\pgfpathlineto{\pgfqpoint{2.432704in}{1.621038in}}%
\pgfpathlineto{\pgfqpoint{2.481891in}{1.588446in}}%
\pgfpathlineto{\pgfqpoint{2.531077in}{1.556742in}}%
\pgfpathlineto{\pgfqpoint{2.580264in}{1.525902in}}%
\pgfpathlineto{\pgfqpoint{2.629450in}{1.495905in}}%
\pgfpathlineto{\pgfqpoint{2.678637in}{1.466731in}}%
\pgfpathlineto{\pgfqpoint{2.727823in}{1.438360in}}%
\pgfpathlineto{\pgfqpoint{2.777010in}{1.410772in}}%
\pgfpathlineto{\pgfqpoint{2.826196in}{1.383948in}}%
\pgfpathlineto{\pgfqpoint{2.875383in}{1.357872in}}%
\pgfpathlineto{\pgfqpoint{2.924570in}{1.332525in}}%
\pgfpathlineto{\pgfqpoint{2.973756in}{1.307892in}}%
\pgfpathlineto{\pgfqpoint{3.022943in}{1.283956in}}%
\pgfpathlineto{\pgfqpoint{3.072129in}{1.260702in}}%
\pgfpathlineto{\pgfqpoint{3.121316in}{1.238115in}}%
\pgfpathlineto{\pgfqpoint{3.170502in}{1.216182in}}%
\pgfpathlineto{\pgfqpoint{3.219689in}{1.194888in}}%
\pgfpathlineto{\pgfqpoint{3.268875in}{1.174221in}}%
\pgfpathlineto{\pgfqpoint{3.318062in}{1.154167in}}%
\pgfpathlineto{\pgfqpoint{3.367249in}{1.134715in}}%
\pgfpathlineto{\pgfqpoint{3.416435in}{1.115853in}}%
\pgfpathlineto{\pgfqpoint{3.465622in}{1.097569in}}%
\pgfpathlineto{\pgfqpoint{3.514808in}{1.079854in}}%
\pgfpathlineto{\pgfqpoint{3.563995in}{1.062696in}}%
\pgfpathlineto{\pgfqpoint{3.613181in}{1.046085in}}%
\pgfpathlineto{\pgfqpoint{3.662368in}{1.030012in}}%
\pgfpathlineto{\pgfqpoint{3.711554in}{1.014468in}}%
\pgfpathlineto{\pgfqpoint{3.760741in}{0.999444in}}%
\pgfpathlineto{\pgfqpoint{3.809928in}{0.984931in}}%
\pgfpathlineto{\pgfqpoint{3.859114in}{0.970921in}}%
\pgfpathlineto{\pgfqpoint{3.908301in}{0.957406in}}%
\pgfpathlineto{\pgfqpoint{3.957487in}{0.944379in}}%
\pgfpathlineto{\pgfqpoint{4.006674in}{0.931832in}}%
\pgfpathlineto{\pgfqpoint{4.055860in}{0.919759in}}%
\pgfpathlineto{\pgfqpoint{4.105047in}{0.908153in}}%
\pgfpathlineto{\pgfqpoint{4.154233in}{0.897008in}}%
\pgfpathlineto{\pgfqpoint{4.203420in}{0.886317in}}%
\pgfpathlineto{\pgfqpoint{4.252606in}{0.876075in}}%
\pgfpathlineto{\pgfqpoint{4.301793in}{0.866277in}}%
\pgfpathlineto{\pgfqpoint{4.350980in}{0.856917in}}%
\pgfpathlineto{\pgfqpoint{4.400166in}{0.847990in}}%
\pgfpathlineto{\pgfqpoint{4.449353in}{0.839492in}}%
\pgfpathlineto{\pgfqpoint{4.498539in}{0.831417in}}%
\pgfpathlineto{\pgfqpoint{4.547726in}{0.823763in}}%
\pgfpathlineto{\pgfqpoint{4.596912in}{0.816524in}}%
\pgfpathlineto{\pgfqpoint{4.646099in}{0.809698in}}%
\pgfpathlineto{\pgfqpoint{4.695285in}{0.803280in}}%
\pgfpathlineto{\pgfqpoint{4.744472in}{0.797267in}}%
\pgfpathlineto{\pgfqpoint{4.793659in}{0.791657in}}%
\pgfpathlineto{\pgfqpoint{4.842845in}{0.786446in}}%
\pgfpathlineto{\pgfqpoint{4.892032in}{0.781631in}}%
\pgfpathlineto{\pgfqpoint{4.941218in}{0.777210in}}%
\pgfpathlineto{\pgfqpoint{4.990405in}{0.773182in}}%
\pgfpathlineto{\pgfqpoint{5.039591in}{0.769543in}}%
\pgfpathlineto{\pgfqpoint{5.088778in}{0.766292in}}%
\pgfpathlineto{\pgfqpoint{5.137964in}{0.763427in}}%
\pgfpathlineto{\pgfqpoint{5.187151in}{0.760947in}}%
\pgfpathlineto{\pgfqpoint{5.236338in}{0.758851in}}%
\pgfpathlineto{\pgfqpoint{5.285524in}{0.757138in}}%
\pgfpathlineto{\pgfqpoint{5.334711in}{0.755806in}}%
\pgfpathlineto{\pgfqpoint{5.383897in}{0.754855in}}%
\pgfpathlineto{\pgfqpoint{5.433084in}{0.754285in}}%
\pgfpathlineto{\pgfqpoint{5.482270in}{0.754094in}}%
\pgfusepath{stroke}%
\end{pgfscope}%
\begin{pgfscope}%
\pgfsetrectcap%
\pgfsetmiterjoin%
\pgfsetlinewidth{0.803000pt}%
\definecolor{currentstroke}{rgb}{0.000000,0.000000,0.000000}%
\pgfsetstrokecolor{currentstroke}%
\pgfsetdash{}{0pt}%
\pgfpathmoveto{\pgfqpoint{0.612802in}{0.754094in}}%
\pgfpathlineto{\pgfqpoint{0.612802in}{3.774094in}}%
\pgfusepath{stroke}%
\end{pgfscope}%
\begin{pgfscope}%
\pgfsetrectcap%
\pgfsetmiterjoin%
\pgfsetlinewidth{0.803000pt}%
\definecolor{currentstroke}{rgb}{0.000000,0.000000,0.000000}%
\pgfsetstrokecolor{currentstroke}%
\pgfsetdash{}{0pt}%
\pgfpathmoveto{\pgfqpoint{5.482270in}{0.754094in}}%
\pgfpathlineto{\pgfqpoint{5.482270in}{3.774094in}}%
\pgfusepath{stroke}%
\end{pgfscope}%
\begin{pgfscope}%
\pgfsetrectcap%
\pgfsetmiterjoin%
\pgfsetlinewidth{0.803000pt}%
\definecolor{currentstroke}{rgb}{0.000000,0.000000,0.000000}%
\pgfsetstrokecolor{currentstroke}%
\pgfsetdash{}{0pt}%
\pgfpathmoveto{\pgfqpoint{0.612802in}{0.754094in}}%
\pgfpathlineto{\pgfqpoint{5.482270in}{0.754094in}}%
\pgfusepath{stroke}%
\end{pgfscope}%
\begin{pgfscope}%
\pgfsetrectcap%
\pgfsetmiterjoin%
\pgfsetlinewidth{0.803000pt}%
\definecolor{currentstroke}{rgb}{0.000000,0.000000,0.000000}%
\pgfsetstrokecolor{currentstroke}%
\pgfsetdash{}{0pt}%
\pgfpathmoveto{\pgfqpoint{0.612802in}{3.774094in}}%
\pgfpathlineto{\pgfqpoint{5.482270in}{3.774094in}}%
\pgfusepath{stroke}%
\end{pgfscope}%
\begin{pgfscope}%
\pgfsetbuttcap%
\pgfsetmiterjoin%
\definecolor{currentfill}{rgb}{1.000000,1.000000,1.000000}%
\pgfsetfillcolor{currentfill}%
\pgfsetfillopacity{0.800000}%
\pgfsetlinewidth{1.003750pt}%
\definecolor{currentstroke}{rgb}{0.800000,0.800000,0.800000}%
\pgfsetstrokecolor{currentstroke}%
\pgfsetstrokeopacity{0.800000}%
\pgfsetdash{}{0pt}%
\pgfpathmoveto{\pgfqpoint{3.651640in}{2.711428in}}%
\pgfpathlineto{\pgfqpoint{5.365604in}{2.711428in}}%
\pgfpathquadraticcurveto{\pgfqpoint{5.398937in}{2.711428in}}{\pgfqpoint{5.398937in}{2.744762in}}%
\pgfpathlineto{\pgfqpoint{5.398937in}{3.657428in}}%
\pgfpathquadraticcurveto{\pgfqpoint{5.398937in}{3.690761in}}{\pgfqpoint{5.365604in}{3.690761in}}%
\pgfpathlineto{\pgfqpoint{3.651640in}{3.690761in}}%
\pgfpathquadraticcurveto{\pgfqpoint{3.618307in}{3.690761in}}{\pgfqpoint{3.618307in}{3.657428in}}%
\pgfpathlineto{\pgfqpoint{3.618307in}{2.744762in}}%
\pgfpathquadraticcurveto{\pgfqpoint{3.618307in}{2.711428in}}{\pgfqpoint{3.651640in}{2.711428in}}%
\pgfpathclose%
\pgfusepath{stroke,fill}%
\end{pgfscope}%
\begin{pgfscope}%
\pgfsetbuttcap%
\pgfsetroundjoin%
\pgfsetlinewidth{3.011250pt}%
\definecolor{currentstroke}{rgb}{1.000000,0.000000,0.000000}%
\pgfsetstrokecolor{currentstroke}%
\pgfsetdash{{11.100000pt}{4.800000pt}}{0.000000pt}%
\pgfpathmoveto{\pgfqpoint{3.684973in}{3.565761in}}%
\pgfpathlineto{\pgfqpoint{4.084973in}{3.565761in}}%
\pgfusepath{stroke}%
\end{pgfscope}%
\begin{pgfscope}%
\pgftext[x=4.218307in,y=3.507428in,left,base]{\rmfamily\fontsize{12.000000}{14.400000}\selectfont AP}%
\end{pgfscope}%
\begin{pgfscope}%
\pgfsetbuttcap%
\pgfsetroundjoin%
\pgfsetlinewidth{3.011250pt}%
\definecolor{currentstroke}{rgb}{0.800000,0.600000,0.000000}%
\pgfsetstrokecolor{currentstroke}%
\pgfsetdash{{19.200000pt}{4.800000pt}{3.000000pt}{4.800000pt}}{0.000000pt}%
\pgfpathmoveto{\pgfqpoint{3.684973in}{3.333428in}}%
\pgfpathlineto{\pgfqpoint{4.084973in}{3.333428in}}%
\pgfusepath{stroke}%
\end{pgfscope}%
\begin{pgfscope}%
\pgftext[x=4.218307in,y=3.275095in,left,base]{\rmfamily\fontsize{12.000000}{14.400000}\selectfont RAP}%
\end{pgfscope}%
\begin{pgfscope}%
\pgfsetbuttcap%
\pgfsetroundjoin%
\pgfsetlinewidth{3.011250pt}%
\definecolor{currentstroke}{rgb}{0.000000,0.000000,1.000000}%
\pgfsetstrokecolor{currentstroke}%
\pgfsetdash{{3.000000pt}{4.950000pt}}{0.000000pt}%
\pgfpathmoveto{\pgfqpoint{3.684973in}{3.101095in}}%
\pgfpathlineto{\pgfqpoint{4.084973in}{3.101095in}}%
\pgfusepath{stroke}%
\end{pgfscope}%
\begin{pgfscope}%
\pgftext[x=4.218307in,y=3.042761in,left,base]{\rmfamily\fontsize{12.000000}{14.400000}\selectfont DR}%
\end{pgfscope}%
\begin{pgfscope}%
\pgfsetrectcap%
\pgfsetroundjoin%
\pgfsetlinewidth{3.011250pt}%
\definecolor{currentstroke}{rgb}{0.000000,0.501961,0.000000}%
\pgfsetstrokecolor{currentstroke}%
\pgfsetdash{}{0pt}%
\pgfpathmoveto{\pgfqpoint{3.684973in}{2.868762in}}%
\pgfpathlineto{\pgfqpoint{4.084973in}{2.868762in}}%
\pgfusepath{stroke}%
\end{pgfscope}%
\begin{pgfscope}%
\pgftext[x=4.218307in,y=2.810428in,left,base]{\rmfamily\fontsize{12.000000}{14.400000}\selectfont AAMR \(\displaystyle \equiv\) GAP}%
\end{pgfscope}%
\end{pgfpicture}%
\makeatother%
\endgroup%